\documentclass[psamsfonts]{article}

\usepackage{amsfonts}
\usepackage{amsmath}
\usepackage{amssymb}
\usepackage{amscd}
\usepackage{amstext}
\usepackage{amsthm}

 \title{
 Simple nuclear $C^*$-algebras of tracial topological rank one\,
 \thanks{Research partially supported by NSF grant DMS 0097903.
     AMS 1991 Subject Classification Numbers:
             Primary 46L05,
                46L35.
            Key words:  TAI $C^*$-algebras, Classification
            of nuclear $C^*$-algebras, Tracial rank one
                \protect\\}}

\author{Huaxin Lin\\
 Department of Mathematics\\
University of Oregon\\
Eugene, Oregon 97403-1222}
\date{}
\begin{document}
\maketitle

\newcommand{\CA}{$C^*$-algebra}
\newcommand{\aue}{approximate unitary equivalence}
\newcommand{\ayue}{approximately unitarily equivalent}
\newcommand{\mops}{mutually orthogonal projections}
\newcommand{\hm}{homomorphism}
\newcommand{\pisca}{purely infinite simple \CA}
\newcommand{\andeqn}{\,\,\,\,\,\, {\rm and} \,\,\,\,\,\,}
\newcommand{\QED}{\rule{1.5mm}{3mm}}
\newcommand{\morp}{contractive completely
positive linear map}
\newcommand{\asmorp}{asymptotic morphism}
\newcommand{\arrow}{\rightarrow}
\newcommand{\tdsum}{\widetilde{\oplus}}
\newcommand{\pa}{\|}  
\newcommand{\ep}{\varepsilon}
\newcommand{\id}{{\rm id}}
\newcommand{\aueeps}[1]{\stackrel{#1}{\sim}}
\newcommand{\aeps}[1]{\stackrel{#1}{\approx}}
\newcommand{\dt}{\delta}
\newcommand{\yu}{\fang}
\newcommand{\ca}{{\cal C}_1}
\newcommand{\Co}{unital separable simple $C^*$-algebra of real rank zero and
stable rank one with  weakly unperforated $K_0$-group}
\newcommand{\n}{\diamond}
\newcommand{\SCA}{$C^*$-subalgebra}
\newcommand{\rforal}{\,\,\,{\rm for\,\,\,all}\,\,\,}
\begin{abstract}

We give a classification theorem
for unital separable nuclear simple \CA s with tracial
rank no more than one. Let $A$ and $B$ be two
unital separable simple nuclear \CA s with $TR(A), \,TR(B)\le 1$
which satisfy the universal coefficient theorem. We show that
$A\cong B$ if and only if there is an order and unit preserving
isomorphism
$$
\gamma=(\gamma_0, \gamma_1, \gamma_2): (K_0(A), K_0(A)_+, [1_A], K_1(A), T(A))
\cong
(K_0(B), K_0(B)_+, [1_B], K_1(B), T(B)),
$$
where $\gamma_2^{-1}(\tau)(x)=\tau(\gamma_0(x))$ for each $x\in K_0(A)$
and $\tau\in T(B).$

\end{abstract}

\section{Introduction}

This paper is a part of the program to classify nuclear
simple \CA s initiated by George A. Elliott.
By a classification theorem for  a class of nuclear \CA s, one usually
means the following: Two \CA s in the class with the same $K$-theoretical
data are isomorphic (as \CA s) and given a set of $K$-theoretical
data there is a \CA\, in the class which possesses the $K$-theoretical
data. One should note that
 $K$-theoretical data
is certainly not a complete
invariance for commutative \CA s in general.
 This is one of the reasons that we study \CA s with ``lower rank".
Moreover, we are only interested in non-commutative \CA s.
By \CA s of ``lower rank", one often means that the \CA s have real rank zero,
or stable rank one. It has been known (but not surprisingly)
 many important \CA s arising from
applications are of real rank zero or stable rank one. Notably,
all purely infinite simple \CA s have real rank zero and many \CA
s arising from dynamical systems are of stable rank one. One of
the classical results of this kind states that all irrational
rotation \CA s are simple nuclear \CA s with real rank zero and
stable rank one (see \cite{EE} and \cite{P}). This is another
reason that we study \CA s of ``lower rank". One may view (simple)
\CA s of real rank zero and stable rank one as some kind of
generalization of AF-algebras. A more suitable generalization of
AF-algebras has been demonstrated to be \CA s with tracial
topological rank zero. Simple \CA s with tracial topological rank
zero have real rank zero, stable rank one, with weakly
unperforated $K_0$ and are  quasidiagonal. All simple AH-algebras
with slow dimension growth and with real rank zero have tracial
topological rank zero. This shows that simple \CA s with zero
tracial rank could have rich $K$-theory.
Simple AH-algebras with slow dimension growth and
with real rank zero have been classified in \cite{EG} ( together
with
\cite{D}, \cite{G4}  and \cite{G5}). A classification theorem for
unital nuclear separable simple \CA s with tracial topological
rank zero which satisfy the UCT was given in \cite{Lnmsri} (see
also \cite{Lnan}, \cite{LnC} and \cite{DE1} for earlier
references). Simple \CA s with tracial topological rank zero are
also called TAF (tracially AF) \CA s.

This paper studies \CA s of ``rank one". A standard example of a
\CA\, with stable rank one is of course $M_k(C([0,1])).$ A notion
of ``tracially approximately interval" \CA s (TAI \CA s) are
introduced in this paper----see definition 2.2 below. It turns out
that simple TAI \CA s are the same as simple \CA s with tracial
topological rank no more than one. Roughly speaking, TAI \CA s are
those \CA s whose finite subsets can be approximated by \SCA s
which are finite direct sums of finite dimensional \CA s and
matrix algebras over $C([0,1])$ in ``measure" or rather in trace.
It is proved here that simple TAI \CA s have stable rank one. From
a result of G. Gong (\cite{G2}) we observe that all simple
AH-algebra with very slow dimension growth are in fact TAI \CA s.
It is also shown here that simple TAI \CA s are quasidiagonal,
their ordered $K_0$-groups are weakly unperforated and satisfy the
Riesz interpolation property, and these \CA s also satisfy the
Fundamental Comparison Property of Blackadar. All these suggest
that the notion of TAI \CA s is a suitable notion for ``rank one
(or less)" \CA s. As mentioned above, these \CA s have tracial
topological rank one or zero.

Elliott, Gong and Li in \cite{EGL} (also \cite{G2})
give a complete classification (up to isomorphism)
for simple AH-algebras with
very slow dimension growth by their $K$-theoretical data
(An important special case can be found in
K. Thomsen's work \cite{T2}).
These \CA s are nuclear separable simple \CA s of stable
rank one. Their work is
a significant advance in classifying
finite simple \CA s
after the remarkable result of
\cite{EG} which classifies simple AH-algebras of real rank zero
(with slow dimension growth). Therefore, it is  the time
to classify nuclear simple separable finite \CA s with
real rank other than zero
{\it without assuming} that they are inductive limits (AH-algebras are
inductive limits of finite direct sums of some
standard homogeneous \CA s) of certain special ``building blocks".
The main purpose of this paper is to present such a result.
(In a recent paper \cite{LO}, using Kishmoto's work for crossed products,
it is shown that many simple crossed products can be constructed
to have tracial topological rank one).

Section 1-6 and section 8 and most of section 9 were written in 1998. Together
with later sections, the original preprint has two parts.
A preliminary report on the results in the two-part
preprints
was reported in EU conference
on Operator Algebra in Copenhagen in the August 1998.
Since then a great deal of progress on the subject
has been made. The present paper absorbs both parts of the original
preprint and reflects the new development. But it is significantly
shorter than the original preprint.
More importantly,
the main result of the paper has been greatly improved and
a technical condition in original preprint has been removed:

Let $A$ and $B$ be two unital separable nuclear simple \CA s with
$TR(A)\le 1,$ $TR(B)\le 1$ and satisfying the UCT. Then $A\cong B$
if and only if they have the same $K$-theoretical data ---see
Theorem 10.10.

Consider two \CA s $A$ and $B$ as above.
As in \cite{EGL}, we will construct the following approximately commutative
diagram:
\vspace{-0.05in}
$$
\begin{array}{ccccc}
B & {\stackrel{{\rm id}_B}{\rightarrow}} & B\\
\hspace{0.1in} \downarrow_{L_1} & \nearrow_{\Phi_1} & \hspace{0.1in}\downarrow_{L_2}\\
A & {\stackrel{{\rm id}_A}{\rightarrow}} & A.
\end{array}
$$
We then apply an approximately intertwining argument of Elliott
to obtain an isomorphism.

It requires two types of results: an ``existence theorem" and
a ``uniqueness theorem".
The existence theorem will state that
given the same $K$-theoretical data of $A$ and $B,$ one should
be able to produce maps from $A$ to $B$ and $B$ to $A$
which carry the $K$-theoretical data. So one has at least
maps between these \CA s. In order to get the ``approximately commutative"
diagram, one also needs the uniqueness theorem:
two maps with the same $K$-theoretical data are approximately
unitarily equivalent.

Denote by ${\cal I}$ the class of unital \CA s
which are finite direct sums of \CA s with the form $M_k,$
or $M_k(C([0,1]).$
An important fact is the following classification
of monomorphisms from $\oplus_{k=1}^n M_{r(k)}(C([0,1])$ to a simple TAI-algebra:
Let $A=\oplus_{k=1}^n M_{r(k)}(C([0,1])$ and $B$ be a unital simple TAI \CA.
Suppose that $\phi_i: A\to B$ ($i=1,2$) are two unital monomorphisms
such that
$(\phi_1)_*=(\phi_2)_*$ on $K_0(A)$ and
$$
\tau\circ\phi_1(a)=\tau\circ\phi_2(a)
$$
for all $a\in A$ and all tracial states $\tau.$ Then
there exists a sequence of unitaries $u_n\in B$ such that
$$
\lim_{n\to\infty}u_n^*\phi_1(a)u_n=\phi_2(a)\,\,\,\,{\rm for\,\,\,all}
\,\,\, a\in A.
$$
Combining this with the more general uniqueness theorem of \cite{LnA},
we are able to obtain a uniqueness theorem
for  nuclear separable simple
TAI \CA s.
As in \cite{EGL}, the $K$-theoretical data mentioned above
contains not only $(K_0(A), K_0(A)_+, [1_A], K_1(A))$
but an additional tracial data, namely the tracial state space
 $T(A)$ together
with a paring of $T(A)$ and $K_0(A).$ Since traces become part of
the invariant as in \cite{Ell2}, \cite{Li} and \cite{EGL}, we also
use some earlier results of J. Cuntz and G. K. Pedersen. However,
we are able to avoid some difficult topological techniques
involving higher dimensional CW-complexes in \cite{EGL}. We also
show that the set of ``$K$-theoretical data" for unital separable
simple \CA s with $TR(A)\le 1$ which satisfy the UCT is the same
as that of simple AH-algebras with no dimension growth as
described by J. Villadsen (\cite{V}), namely the paring of $T(A)$
and $K_0(A)$ maps extremal traces to extremal states on $(K_0(A),
[1_A]).$ The uniqueness theorem has also to be adjusted to deal
with other complication caused by the fact that \CA s are no
longer assumed to have real rank zero. A careful treatment on
exponential length is needed. Existence theorem also needs to be
improved from that in \cite{Lnmsri}. The existence theorem should
also control the exponential length. It turns out that when \CA s
are assumed to have only torsion $K_1,$ the proof can be made much
shorter. This was done without using de la Harpe and Skandalis
determinants as in \cite{EGL}.

The paper is organized as follows. Section 2 gives the definition
of TAI \CA s. Section 3 gives some elementary properties of simple
TAI \CA s. In section 4, we show that simple TAI \CA s have stable
rank one, weakly unperforated $K_0,$ the fundamental comparison
property, and are MF. Starting from section 5, we will use term ``
simple \CA s with $TR(A)\le 1$ instead of ``TAI \CA\, $A$". Even
though that the term ``$TR(A)\le 1$" has appeared in \cite{Lntr1},
the term ``TAI" has been used and results in the first  4 sections
have been quoted in a number of  places including \cite{Lntr1}. We
feel that we can keep the literature consistent by keeping the
term ``TAI" in the  first 4 sections here. In section 5, we show
that every simple (nonelementary) \CA\, $A$ with $TR(A)\le 1$ is
tracially approximately divisible. We also give a classification
theorem for monomorphisms from $M_r(C([0,1]))$ to a unital simple
TAI \CA\, mentioned earlier. In section 6, we study the unitary
group of a simple \CA\, $A$ with $TR(A)\le 1.$ Exponential rank of
a simple \CA\, $A$ with $TR(A)\le 1$ is proved to be no more than
$3+\ep.$ Let $CU(A)$ be the closure of commutator group of $U(A).$
We show that $U_0(A)/CU(A)$ is always divisible and if $A$ is
simple and $TR(A)\le 1$ then $U_0(A)/CU(A)$ is torsion free. In
section 7, we present some  results concerning \hm\, from
$U(C)/CU(C)$ to $U(B)/CU(B),$ where $B$ is a unital simple \CA\,
with $TR(B)\le 1$ and $C$ is some very special unital \CA s. These
results may be viewed as part of the existence theorem which
controls the exponential length of unitaries under certain maps.
In section 8, we present a uniqueness theorem suitable to use in
the proof of 10.4  which is based on results in \cite{LnA}. One
immediate consequence of that is the following: Let $A$ be a
unital separable simple nuclear \CA\, with $TR(A)\le 1$ and with
torsion $K_1(A).$ Then an automorphism $\alpha: A\to A$ is
approximately inner if and only if $[\alpha]=[{\rm id}_A]$ in
$KL(A)$ and $\tau\circ \alpha(x)=\tau(x)$ for all $x\in A_{s.a}$
and $\tau\in T(A).$ In section 9, we present  several versions of
existence theorem. The purpose is to establish a map from $A$ to
$B$ if $TR(A)\le 1,$ $TR(B)\le 1$ and $(K_0(A), K_0(A)_+, [1_A],
K_1(A), T(A)) =(K_0(B), K_0(B)_+, [1_B], K_1(B), T(B)),$ where
$T(A)$ and $T(B)$ are tracial state spaces of $A$ and $B,$
respectively. Finally, in section 10, we give the proof of the
main theorem--Theorem 10.4 and Theorem 10.10.

\vspace{0.1in}

The following terminology and notation will
be used throughout this  paper.

Let $A$ be a \CA.

(i) Two projections in $A$ are said to be equivalent if they are
Murray-von Neumann equivalent. We write $p\preceq q$ if
$p$ is equivalent to a projection in $qAq.$ We use $[p]$ for
the equivalence class of projections equivalent to $p.$
Let $a\in A_+.$ We write $p\preceq a$ if $p\preceq q$ for some projection
$q\in \overline{aAa}.$

(ii) An element in $A$ is said to be {\it full} if
the (closed) ideal generated by the element is $A$ itself. Every
nonzero element in a simple \CA \, is full.

(iii) Let $\ep>0,$ ${\cal F}$ and $S$ be a subset of $A.$
We write $x\in_{\ep} S$ if there exists $y\in S$ such that
$\|x-y\|<\ep,$ and write ${\cal F}\subset_{\ep} S,$ if
$x\in_{\ep} S$ for all $x\in {\cal F}.$

(iv) Denote by $T(A)$ the compact convex set of the
tracial state space of a unital separable $A.$

(v) Let ${\cal G}\subset A$ and $\dt>0.$
A  \morp\, $L: A\to B$ is said to be
${\cal G}$-$\dt$-multiplicative if
$$
\|L(ab)-L(a)L(b)\|<\dt\,\,\,{\rm for\,\,\,all}\,\,\, a,b\in {\cal G}.
$$

(vi) Let $X$ be a compact metric space and $h: PM_r(C(X))P\to A,$
where $P\in M_r(C(X))$ is a projection, be
a \hm. We say $h$ is {\it homotopically trivial}, if $h$
is homotopy to a point-evaluation.
A \morp\, $L: PM_r(C(X))P\to A$ is said to be
{\it homotopically trivial}, if $L$ factors through
to a homotopically trivial \hm, i.e., $L=L'\circ h,$ where
$h$ is homotopically trivial.


\vspace{0.1in}

{\bf Acknowledgements} This work is partially
supported by a grant from NSF.
We wish to thank
Guihua Gong and Chris Phillips for helpful  correspondences and
G. A. Elliott, Guihua Gong  and Lianqing Li for sharing their manuscripts.

\section{ Definition of tracially AI \CA s}

{\bf 2.1. Definition}
We denote by ${\cal I}$ the class of all unital \CA s
with the form $\oplus_{i=1}^n B_i,$ where each $B_i\cong M_k$
for some integer $k$ or
$B_i\cong M_k(C([0,1])).$
Let $A\in {\cal I}.$ Then we have the following well known
facts.

(i) Every \CA \, in ${\cal I}$ is of stable rank one;

(ii) Two projections $p$ and $q$ in a \CA\, $A\in {\cal I}$
are equivalent if and only if
$\tau(p)=\tau(q)$ for all $\tau\in T(A);$

(iii) For any $\ep>0$ and any finite subset ${\cal F}\subset
A,$ there exist $\dt>0$ and a finite subset ${\cal G}\subset A$
satisfying the following:
if $L: A\to B$ is a ${\cal G}$-$\dt$-multiplicative \morp,
where $B$ is a \CA,
then there exists a \hm\, $h: A\to B$
such that
$$
\|h(a)-L(a)\|<\ep
\,\,\,\,\,\rforal a\in {\cal F}.
$$
These facts will be used throughout the paper
without further notice.

{\bf 2.2.  Definition} A unital \CA \, $A$ is said to be
{\it tracially AI} (TAI) if for any finite subset ${\cal F}\subset A$
containing a nonzero element $b,$
$\ep>0,$ integer $n>0$ and any full element $a\in A_+,$
there exists a nonzero projection $p\in A$ and
a $C^*$-subalgebra $I\in {\cal I}$ with $1_{I}=p,$
such that

(1) $\|[x,p]\|<\ep$ for all $x\in {\cal F},$

(2) $pxp\in_{\ep} I$ for all $x\in {\cal F}$
and $\|pbp\|\ge \|b\|-\ep,$

(3) $n[1-p]\le [p]$ and $1-p\preceq a.$

A non-unital \CA \, $A$ is said to be TAI if ${\tilde A}$
is TAI.

In 4.10, we show that, if $A$ is simple,
condition (3) can replaced by

(3') $1-p$ is unitarily equivalent to a projection in
$eAe$ for any previously given nonzero projection $e\in A.$

If $A$ has the Fundamental Comparability (see \cite{Bl2}),
condition (3) can be replaced by

(3'') $\tau(1-p)<\sigma$ for any prescribed $\sigma>0$ and for
all normalized quasi-traces of $A.$

From the definition, one sees that the ``part " of $A$ which
may not be approximated by \CA s in ${\cal I}$ has
small 'measure" or trace.
Note in the above, if ${\cal I}$ is replaced by
finite dimensional \CA s, then it is precisely the definition
of TAF \CA s ( see \cite{LnB}).

\vspace{0.1in}

{\bf 2.3. Examples}
Every AF-algebra is TAI.
Every TAF \CA \, introduced in \cite{LnB} is
a TAI \CA. However, in general, TAI \CA s have real rank
other than zero. In 4.5 we will show that every
simple TAI \CA \, has stable rank one, which implies that
simple TAI \CA s have real rank one or zero.
It is obvious that every direct limit of \CA s in ${\cal I}$ is
a TAI \CA. These \CA s provide many examples of TAI \CA s
that have real rank one. However, TAI \CA s may not be
inductive limits of \CA s in ${\cal I}.$

Let $A=\lim_{n\to\infty} (A_n, \phi_{n,m}),$ where
$A_n=\oplus_{i=1}^{s(n)} P_{n,i}M_{(n,i)}(C(X_{n,i}))P_{n,i},$
$X_{n,i}$ is a finite dimensional compact metric space
and $P_{n,i}\in M_{(n,i)}(C(X_{n,i}))$ is a projection
for all $n$ and $i.$ Such a \CA\, is called an AH-algebra.
Suppose that $A$ is unital. Following \cite{G2},
$A$ is said to have {\it very slow dimension growth} if
$$
\lim_{n\to\infty} \min_{i}{{\rm rank}(P_{n,i})\over{
({\rm dim}X_{n,i}+1)^3}}=\infty.
$$
$A$ is said to have {\it no dimension growth} if
there is an integer $m>0$ such that
${\rm dim}X_{n,i}\le m.$
Note these \CA s may not be of real rank zero. Since these
\CA s could have non-trivial $K_1$-groups (see 10.1), they
are not inductive limits of \CA s in ${\cal I}.$
In \cite{LnB}, example of simple TAF \CA s which
are non-nuclear were given. In particular, there are simple
TAI \CA s that are not even nuclear.

\vspace{0.1in}

{\bf 2.4. Lemma} {\it Let $a$ be a positive element in a unital
\CA\, $A$ with $sp(a)\subset [0,1].$ Then for any $\ep>0,$
there exists $b\in A_+$ such that $sp(b)$ is a union of
finitely many mutually disjoint closed intervals and  finitely many points
and
$$
\|a-b\|<\ep.
$$}

{\it Proof:} Fix $\ep>0.$
Let $I_1,I_2,...,I_k$ be all disjoint closed intervals in $sp(a)$
with length at least  $\ep/8$ such that if $I\supset I_j,$ is an interval,
then $I\not\subset sp(a).$
Let $d'=\min\{{\rm dist}(I_i, I_j), i\not=j\}$ and
$d=\min(d'/2, \ep/16).$

Choose  $J_i=\{\xi\in [0,1]: {\rm dist}(\xi, I_j)<d_i\},$ $i=1,2,...,k$
with $d_i\le d$ and the endpoints of $J_i$ are not in $sp(a).$
Since the endpoints of $J_i$ are not in $sp(a),$ there are
open intervals $J_i'\subset J_i$ such that ${\bar J_i}'\subset J_i$
and $I_i\subset J_i'.$
Set $Y=sp(a)\setminus (\bigcup_{i=1}^k J_i).$
Then $Y=sp(a)\setminus (\bigcup_{i=1}^k J_i')=sp(a)\setminus (\bigcup_{i=1}^k
{\bar J_i}').$
Since $Y$ is compact and $Y$ contains no intervals with length
more than $\ep/8,$ it is routine to show that there are finitely
many disjoint closed intervals $K_1,K_2,...,K_n$ in $[0,1]\setminus
(\bigcup_{i=1}^kJ_i)$ with
length no more than $9\ep/64$ such that
$Y\subset \bigcup_{j=1}^mK_j.$
Note that $\{{\bar J_1}', {\bar J_2'},..., {\bar J}_k',
K_1, K_2,...,K_n\}$ are disjoint closed intervals.
Fix a point $\xi_j\in K_j,$ $j=1,2,...,m.$
One can define a continuous function $f: (\bigcup_{i=1}^k {\bar J}_i')
\bigcup (\bigcup_{s=1}^n K_s)
\to [0,1]$ which
maps each $J_i$ onto $I_i,$ $i=1,2,...,k$ and
maps $K_j$ to a single point $\xi_j$ such that
$$
|f(\xi)-\xi|<\ep/2\,\,\, {\rm for\,\,\, all}\,\,\, \xi\in [0,1].
$$
Define $b=f(a).$ We see that $b$ meets the requirements
of the lemma.
\QED

{\bf 2.5. Theorem} {\it Let $A$ be a unital simple AH-algebra\,
with very slow dimension growth. Then $A$ is TAI.}

{\it Proof:} By 6.31 in \cite{G2}, 2.1 in \cite{EGL} and
6.39 in \cite{G2},
we may write $A=\lim_{n\to\infty} (A_n, \phi_{n,m}),$
where $A_n=\oplus_{i=1}^{i(n)}M_{n(i)}(C(X_{n,i})),$
$X_{n,i}$ are simplicial complexes with ${\rm dim}X_{n,i}\le 3$
and $\phi_{n,m}$ are injective.
Furthermore, for any finite subset ${\cal F}\subset A_n$ and $\ep>0,$
if $m$ is large enough, then there are two mutually orthogonal
projections $P,\, Q\in A_m$ and two \hm s $\phi: A_n\to PA_mP$ and
$\psi: A_n\to QA_mQ$ such that

(1) $\|\phi_{n,n+1}(f)-(\phi(f)\oplus \psi(f))\|<\ep/2$ for all $f
\in {\cal F};$

(2) $\|SPV(\phi(f))\|<\ep/2$ \,(see \cite{G2} for the definition
of $SPV$) and

(3) $\psi$ factors through a matrix algebra over $C([0,1]).$

It follows from Lemma 2.4 that there is a unital $C^*$-subalgebra
$B_1\in {\cal I}$ of $QA_mQ$ such that
$$
\psi(f)\in_{\ep} B_1\,\,\,\,\,\,{\rm for}\,\,\,f\in {\cal F}.
$$
By applying 2.21 in \cite{EG}, we obtain a projection $P_1\le P$ and a
finite dimensional $C^*$-subalgebra $B_2$ with $1_{B_2}=P_1$ such that
$$
\|[\phi(f), P_1]\|<\ep
$$
for all $f\in {\cal F},$ $P_1\phi(f)P_1\subset B_2$ and
$N[1-P_1-Q]\le [P_1+Q]$ and $[1-P_1-Q]\le [e]$ for any given integer $N>0$
and any prescribed projection $e\in A.$
It follows that $A$ is TAI.
\QED

\vspace{0.1in}

{\bf 2.6. Remark} It is not necessary to employ the full strength
of \cite{G2} to prove the above theorem. In fact, it follows from 4.34 and
4.35 in
\cite{G2}.

\section{Elementary properties of simple TAI \CA s}

{\bf 3.1. Lemma} {\it For any $d>0$ there
are $f_1,f_2,...,f_m\in C([0,1])_+$ satisfying
the following.
For any $n,$ and any positive element $x\in B=M_n(C([0,1]))$
with $\|x\|\le 1$ if there exist $a_{ij}\in B,$ $i=1,2,...,n(j),$
$j=1,2,...,m$
with
$$
\|\sum_{i=1}^{n(j)}a_{i,j}f_j(x)a_{i,j}^*-1_A\|<1/2\,\,\,\,\,\,j=1,2,...,m,
$$
then, for any subinterval $J$ of $[0,1]$ with $\mu(J)\ge d$
($\mu$ is the Lebesgue measure),

$sp(\pi_t(x))\cap J\not=\emptyset $ for all $t\in [0,1],$
where $\pi_t: B\to M_n$ is the point evaluation at $t.$
Moreover, denote by $N=max\{n(j): j=1,2,...,m\},$
$$
|sp(\pi_t(x))\cap J|\ge 1/N |sp(\pi_t(x))|,
$$
where $|S|$ means the number of elements in the finite set $S$
(counting multiplicities).
}

{\it Proof:}
Divide $[0,1]$ into $m$ closed subintervals $\{J_j\}$
each of which has the same length $<d/4.$
Let $f_j\in C([0,1])$ such that $0\le f_j\le 1,$ $f_j(t)=1$ for $t\in J_j$
and $f_j(t)=0$ for ${\rm dist}(t, J_j)\ge \mu(J_j).$
Note that, for any subinterval $J$ with $\mu(J)\ge d,$
there exists $j$ such that $J_j\subset J.$
For any $t\in [0,1],$ set
$$
I=\{g\in B: g(t)=0\}.
$$
Then $I$ is a (closed) ideal of $A.$
If $sp(\pi_t(x))\cap J=\emptyset,$
there would be $j$ such that $\pi_t(f_j(x))=0.$
Therefore $f_j\in I.$ But this is impossible, since there is an
element  $z\in B$ with
$$
z(\sum_{i=1}a_{i,j}f_j(x)a_{i,j}^*)=1_B.
$$

For the last part of the lemma, fix $t\in [0,1]$ and
an interval $J$ with $t\in J$ and $\mu(J)\ge d.$
Let $\pi_t(B)=M_{l(t)}.$
Then $|sp(\pi_t(x))|=l(t).$
Suppose that $J_j\subset J$ so that
$f_j(t)=0$ for all $t\not\in J.$
Let $q_t$ be the spectral projection
of $\pi_t(x)$ in $M_{l(t)}$ corresponding to $J.$  Then
$q_t\ge f_j(\pi_t(x)).$ An elementary linear
algebra argument shows that
${\rm rank}q_t\ge (1/N)l(i).$
\QED

\vspace{0.1in}

{\bf 3.2. Theorem} {\it Every  unital
simple \CA\,
satisfying (1) and (2) in 2.2
has property (SP), i.e.,
every hereditary $C^*$-subalgebra contains a nonzero projection.}

{\it Proof:}
Let $A$ be a unital simple \CA\, satisfying (1) and (2) and $B\subset A$ be a hereditary
$C^*$-subalgebra.
We may assume that $A$ is not elementary. Thus $B$ is not elementary.
By p.61 (item 4) in \cite{AS}, there is $a\in B_+$ such that $sp(a)=[0,1].$
It suffices to show that $B_1=\overline{aAa}$ has a nonzero
projection.

Let $d$ be a positive number with $0<d<1/16.$
Let $f_1,f_2,...,f_m$ be as in 3.1.
Since $A$ is simple, there are $a_{ij}\in A$ such that
$$
\|\sum_{i=1}a_{i,j}f_j(a)a_{i,j}^*-1_A\|<1/8\,\,\,\,\,\,j=1,2,...,m.
$$
Let $g\in C_0((0,1))$ with $0\le g \le 1,$
$g(t)=1$ if $g\in [1/4,3/4].$

Denote by ${\cal F}$ the subset
$$
\{a, g(a), f_j(a), a_{i,j}, a_{ij}^*: i,j\}.
$$
For any $\ep>0,$ there exists a projection $p\in A$ and
a $C^*$-subalgebra $C\in {\cal I}$ with $1_{I}=p$
such that

(1) $\|[b,p]\|<\ep$ for all $b\in {\cal F},$

(2) $pbp\in_{\ep} C$ and $\|pap\|\ge 1-\ep.$

A standard perturbation argument shows that,
for any $\eta>0,$  with sufficiently
small $\ep,$
there is a \hm\, $\phi: C^*(pap)\to C$
(where $C^*(pap)$ is the $C^*$-subalgebra generated by $pap$)
and there are $b_{ij}\in C$
such that
$$
\|pap-\phi(pap)\|<\eta,\,\,\,\,\,\|\phi(f_j(a))-pf_j(a)p\|<\eta,
$$
$$
\|\phi(g(a))-pg(a)p\|<\eta,\,\,\,
\andeqn \|\sum_{i=1}b_{ij}\phi(f_j(pap))b_{ij}^*-p\|<1/4
$$
for $j=1,2,..,m.$

Write $C=\oplus_k M_{m(k)}(C([0,1]))$ and $C_k=M_{m(k)}(C([0,1])).$
By Lemma 3.1, with sufficiently small $d,$ we may assume
that
$$
|sp(\pi_t(\phi_k(pap)))\cap (0,1]|\ge 4
$$
for each $t$ and $k,$ where $\pi_t$ is the point-evaluation at $t\in [0,1],$
$\phi_k$ is the map from $C([0,1])$ to $M_{m(k)}(C([0,1]))$
induced by $\phi$ and where $``|sp(\pi_t(\phi_k(pap)))|"$
is the number of eigenvalues of $\pi_t(\phi_k(pap))$
counting multiplicities.

Fix $k.$  For $t\in [0,1],$ let $J_t\subset [1/4,3/4]$
be an open interval with $\mu(J_t)\ge d$ whose endpoints are not eigenvalues
of $\pi_t(\phi_k(pap)).$
Let $V_t$ be an open neighborhood of $t$ such that
the end points of $J_t$ are not eigenvalues of $\pi_y(\phi_k(pap))$
for $y\in V_t.$ Let $\xi_{J_t}$ be the characteristic function
on $J_t.$
Then using the continuous functional calculus we can define
a continuous projection valued function $q_t: V_t\to M_{m(k)}$ by
$q_t(y)=\xi_{J_t}(\phi_k(pap)).$ From the previous paragraph,
${\rm rank} (q_t(y))\ge 4.$
It follows from Proposition 3.2 in \cite{DNP} that
there exists a nonzero projection $q\in M_{m(k)}(C([0,1]))$ such that
$q(t)\le f(\phi(pap))(t)$ for all $t\in [0,1],$ where
$f\in C_0((0,1))$ with $0\le f \le 1,$
$f(t)=1$ if $t\in [1/4,3/4]$ (see also the proof of 1.4 of \cite{Ph2}).
In particular, $qf(\phi(pap))=q.$

Let $b=f(a).$
We estimate that
$$
\|bq-q\|<2\eta.
$$
It is standard that if $\eta<1/8,$ there is a projection
$q'\subset \overline{bAb}$ such that
$$
\|q'-q\|<1/2.
$$
This implies that $\overline{bAb}\subset \overline{aAa}\subset B$ contains
a nonzero projection (\,$q'$\,).
\QED

\vspace{0.1in}

{\bf 3.3. Corollary} {\it Let $A$ be a unital simple
\CA \, satisfying {\rm (1)} and {\rm (2)} in 2.2.
Then,  for any integer $N,$
we may assume that $I=\oplus_{i=1}^k M_{m_i}(C([0,1]))\oplus_{j=1}^l
M_{n_i}$ where $m_i, n_i\ge N.$}

{\it Proof:} In the proof of 3.2, we see that if $1/2d\ge N,$
since
$sp(\pi_t(pap))\cap J_j\not=\emptyset$ for each $j,$
then $\pi_t(pap)$ has at least
$N$ distinct eigenvalues (see also the proof of 3.1).
Therefore,
each summand $C$ in the proof 3.2 has rank at least $N.$
\QED

\vspace{0.1in}

{\bf 3.4. Proposition} {\it Let $A$ be a unital  TAI \CA\,
and $e\in A$ be a full projection. Then $eAe$
satisfies
{\rm (1)} and {\rm (2)} in {\rm 2.2},
and for any full positive element $a\in eAe,$
we can have

{\rm (3')} $1-p\preceq a$.

If  $A$ is also simple, $eAe$ is TAI.
}

{\it Proof:} Fix $\ep>0,$ a finite subset ${\cal F}\subset eAe,$
an integer $n>0$ and  nonzero elements $a,b\in eAe$ with $a\ge 0$
and $b\in {\cal F}.$
Let ${\cal F}_1=\{e\}\bigcup {\cal F}.$
Since $A$ is TAI,
there exists $q\in A$ and
a $C^*$-subalgebra $C\in {\cal I}$ with $1_{C}=q$
such that

(i) $\|[x,q]\|<\ep/64$ for all $x\in {\cal F},$

(ii) $qxq\in_{\ep/64} C$ for all  $x\in {\cal F},$
and $\|qbq\|\ge \|b\|-\ep/64;$ and,

(iii) $n[1-q]\le [q]$ and $1-q\preceq a.$

Note that, by the second part of (ii),
$qeq\not=0.$

We estimate that
$$
\|(eqe)^2-eqe\|<\ep/64 \andeqn \|eqe-qeq\|<\ep/32.
$$

Therefore there is a projection $p\in eAe$ such that
$$
\|p-eqe\|<\ep/16.
$$

Consequently,  there is a projection $d\in C$ such that
$$
\|d-p\|<\ep/8.
$$
Note that
$$
\|qp-pq\|<\ep/8+\|qeqe-eqeq\|<\ep/8+\ep/32=5\ep/32,
$$
and
$B=dCd\in {\cal I}.$ With $\ep/2<1/2,$ we obtain a unitary $u\in A$ such that
$$
\|u-1\|<\ep/4
\andeqn u^*du=p.
$$
Set $C_1=u^*Bu.$ Then $C_1\in I$ and $C_1\subset eAe.$
Now $1_{C_1}=p,$

(1) $\|[x, p]\|<\ep/2$ for all $x\in {\cal F},$

(2) $pxp\in_{\ep/2} C_1$ for all $x\in {\cal F}$
and $\|pbp\|\ge \|b\|-\ep/2.$

We also have
$$
\|(e-p)-(1-q)(e-p)(1-q)\|\le
\|(e-p)-(e-p)(1-q)+q(e-p)(1-q)\|
$$
$$
<\|(e-p)q\|+\ep/16< 5\ep/32+\ep/64+\|qeq-qpq\|+\ep/16
$$
$$
<5\ep/32+\ep/64+\ep/16+\|qeq-qeqeq\|+\ep/16<9\ep/32.
$$

We have (with $\ep<1$)

(3') $(e-p)\preceq (1-q)\preceq a.$

Finally, if we assume that $A$ is simple, by 3.2 and 3.3 there is
a nonzero projection $p_1\le p$ such that $n[p_1]\le [p].$ There
is a nonzero projection $p_1\in \overline{aAa}.$ By  applying 3.3,
we obtain  a nonzero projection $q_1\le p_1$ such that $n[q_1]\le
[p_1].$ Applying the first part of the proof to $(e-p){\cal
F}(e-p),$ we obtain a projection $p'\le (e-p)$ and a unital
$C^*$-subalgebra $C_2\in {\cal I}$ with $1_{C_2}=p'$ such that

(1'') $\|[(e-p)x(e-p),p']\|<\ep/2$ for all $x\in {\cal F},$

(2'') $p'xp'\in_{\ep/2}C_2$ for all $x\in {\cal F},$ and

(3'') $(e-p-p')\preceq p_1.$

Now since $n[p_1]\le [p]\le [p+p']$ and
$(e-p-p')\preceq (e-p)\preceq a,$ we obtain

(3) $n[(e-p-p')]\le [p+p']$ and $(e-p-p')\preceq a.$

We also have $\|[x,(p+p')]\|<\ep$ and
$(p+p')x(p+p')\in_{\ep}C_1\oplus C_2$ for all $x\in {\cal F}.$
Hence  $eAe$ is TAI. \QED

\vspace{0.1in}

{\bf 3.5. Corollary} {\it If $A$ is a unital simple TAI \CA,
then condition {\rm (2)} can be strengthened to

{\rm (2')} $pxp\in_{\ep} B$ and
$\|pxp\|\ge \|x\|-\ep$ for all $x\in {\cal F}.$}

We omit the proof.

%
%
%
%
%
%

{\bf 3.6. Theorem} {\it Let $A$ be a unital simple \CA.
Then $A$ is TAI if and only if $M_n(A)$ is TAI for
all $n$ (or for some $n>0$).}

{\it Proof:} If $M_n(A)$ is TAI, then by identifying $A$ with
a unital hereditary $C^*$-subalgebra of $M_n(A)$ and by using 3.4,
we know $A$ is TAI. It remains to prove the `` only if" part.

We prove this in two steps. The first step is
to prove that $M_n(A)$ satisfies (1) and (2) in 2.2.
To do this, we let $\ep>0$ and ${\cal F}$ be a finite
subset of the unit ball of $M_n(A).$
Set ${\cal G}=\{ f_{ij}\in A: (f_{ij})_{n\times n}\in {\cal F}\}.$
Note that ${\cal G}\subset A.$ Since $A$ is TAI, there exists
a projection $p\in A$ and a unital $C^*$-subalgebra $B\in {\cal I}$
such that

(1) $\|[x,p]\|<\ep/2n^2,$

(2) $pxp\in_{\ep/2n^2} B$ for all $x\in {\cal G}$
and for some $x_1\in {\cal G},$ $\|px_1p\|\ge \|x_1\|-\ep/2n^2.$

Put $P=diag(p,p,...,p)\in M_n(A)$ and $D=M_n(B).$
Then, it is easy to check that

(i) $\|[f,P]\|<\ep$ and

(ii) $PfP\in_{\ep} D$ for all $f\in {\cal F}$ and
$\|pf_1p\|\ge \|f_1\|-\ep$ (if $f_1$ is prescribed).

This completes the first step. Now we also know by 3.2 that
$M_n(A)$ has (SP). Let $a\in M_n(A)$ be given. Choose any nonzero
projection $e\in \overline{aM_n(A)a}.$ Since $M_n(A)$ is simple
and has (SP), by 3.1 in \cite{LnB}, there is a nonzero projection
$q\le e$ and $[q]\le [1_A].$ Applying 3.2 in \cite{LnB}, there
exists a nonzero projection $q_1\le q$ such that $(n+1)[q_1]\le
[q].$ In the first step, we can also require, for any integer
$N>0,$  that

(3)  $N[1_A-p]\le [p]$ and $1_A-p\preceq q_1.$

This implies that

(iii) $N[1_{M_n(A)}-P]\le [P]$ and $(1_{M_n(A)}-P)\preceq q\preceq e.$

Therefore  $M_n(A)$ is TAI.
\QED

\vspace{0.1in}

Next we show that every simple TAI \CA \, has the property
introduced by Popa (\cite{Po}).

\vspace{0.1in}

{\bf 3.7. Proposition} {\it Let $A$ be a unital simple
TAI \CA. Then for any finite subset ${\cal F}\subset A$
and $\ep>0,$
there exists a projection $p\in A$ and a finite dimensional
\CA\, $F\subset A$ with $1_F=p$ such that

{\rm (P1)} $\|[x,p]\|<\ep$ and

{\rm (P2)} $pxp\in_{\ep} F$ for all $x\in {\cal F}$ and
$\|pxp\|\ge \|x\|-\ep$ for all $x\in {\cal F}.$

}

{\it Proof:} By 3.5 it is clear that it suffices to prove
the following claim:

For any unital $C^*$-subalgebra $B\in {\cal I},$
the proposition holds for any finite subset
${\cal F}\subset B\subset A.$

This can be further reduced to the case that $B=C([0,1])\otimes
M_k.$ Moreover, it suffices to prove the claim for the case in
which $B=C([0,1]).$ In this reduced case, we only need to consider
the case in which  ${\cal F}$ contains a single element $x\in B,$
where $x$ is the identity function on $[0,1].$ Now, for any
$\ep>0,$ let $\xi_1, \xi_2,...,\xi_n$ be $\ep/4$-dense in $[0,1]$
and ${\rm dist}(\xi_i, \xi_j)>\ep/8$ if $i\not=j.$ Denote by $f_i$
continuous functions with $0\le f_i\le 1,$ which are  one on
$(\xi_i-\ep/32, \xi_i+\ep/32)$ and zero on $[0,1]\setminus
(\xi_i-\ep/16,\xi_i+\ep/16).$ By 3.2, there is  a nonzero
projection $e_i\in \overline{f_iAf_i}.$ Note that $e_ie_j=0$ if
$i\not=j.$ Set $p=\sum_{i=1}^ne_i.$ We estimate that (see Lemma 2
in \cite{Ln1})
$$
\|x-[(1-p)x(1-p)+\sum_{i=1}^n\xi_ie_i\|<\ep/2
\andeqn
$$
(P1) $\|[p,x]\|<\ep.$

Let $F_1$ be the finite dimensional $C^*$-subalgebra generated by
$e_1,e_2,...,e_n.$ Then

(P2) $pxp\in_{\ep} F_1$ and $\|pxp\|\ge \|x\|-\ep.$
\QED

\section{The structure of simple TAI \CA s}

{\bf 4.1. Theorem} {\it Every unital separable simple TAI \CA\, is
MF} (\cite{BK1}).

{\it Proof:} Let $A$ be such a \CA\, and
let $\{x_n\}$ be a dense sequence in the unit ball
of $A.$
By 3.7, there are projections $p_n\in A$ and finite dimensional
$C^*$-subalgebras $B_n$ with $1_{B_n}=p_n$ such that

(1) $\|[p_n,x_i]\|<1/n,$ and

(2) $px_ip\in_{1/n} B_n$ and $\|p_nx_ip_n\|\ge \|x_i\|-1/n$
for $i=1,2,...,n.$

Let ${\rm id}_n: B_n \to B_n$ be the identity map and let
$j: B_n\to M_{K(n)}$ be a unital embedding. We note that
such $j$ exists provided that $K(n)$ is large enough.
By 5.2 in \cite{Pa}, there exists a completely positive
map $L'_n: p_nAp_n\to M_{K(n)}$ such that
$L_n'|_{B_n}=j\circ {\rm id}_n.$ Since $L'_n$ is unital, by 5.9 and 5.10
in \cite{Pa},
$L'_n$ is a contraction. We define $L_n: A\to M_{K(n)}$ by
$L_n(a)=L_n'(p_nap_n).$
Let $y_{i,n}\in B_n$ such that
$\|p_nx_ip_n-y_{i,n}\|<1/n,$ $n=1,2,...\,.$ Then
$$
\|L_n(x_i)-p_nx_ip_n\|\le \|L_n(x_i-y_{i,n})-(y_{i,n}-p_nx_ip_n)\|
<2/n\to 0
$$
as $n\to \infty.$ Combining this with (1) above, we see  that
$$
\|L_n(ab)-L_n(a)L_n(b)\|\to 0
$$
as $n\to\infty.$
Define
$\Phi: A\to \prod_{n=1}^{\infty} M_{m(n)}$ by
sending $a$ to $\{ L_n(a)\}.$
Then $\Phi$ is a completely positive map.
Denote by $\pi: \prod_{n=1}^{\infty} M_{m(n)}\to
\prod_{n=1}^{\infty} M_{m(n)}/\oplus_{n=1}^{\infty}M_{m(n)}.$
Then
$$
\pi\circ \Phi: A\to \prod_{n=1}^{\infty} M_{m(n)}
/\oplus_{n=1}^{\infty}M_{m(n)}
$$
is a (nonzero) \hm.
Since $A$ is simple, $\pi\circ \Phi$ is injective.
It follows from 3.22 in \cite{BK1} that $A$ is an MF-algebra.
\QED

\vspace{0.1in}

{\bf 4.2. Corollary} {\it Every separable unital
\CA\, satisfying {\rm (P1)} and {\rm (P2)} is MF.}

{\it Proof:} We actually proved this above.
Note, simplicity is not needed for injectivity since
$\|p_nxp_n\|\to \|x\|.$
\QED

\vspace{0.1in}

{\bf 4.3. Proposition} {\it Every nuclear
separable
simple TAI \CA \, is quasidiagonal.}

{\it Proof:} As in \cite{BK1}, a separable nuclear MF \CA\, is NF,
and it is quasidiagonal. In fact it is strong NF (see \cite{BK2}).
\QED

\vspace{0.1in}

{\bf 4.4. Corollary} {\it Every unital separable simple TAI \CA\,
has at least one tracial state.}

{\it Proof:} It is well known that
$\prod_{n=1}^{\infty}M_{m(n)}/\oplus _{n=1}^{\infty}M_{m(n)}$
has tracial states. Tracial states are defined by
weak limits of tracial states on each $M_{m(n)}.$
Let $\tau$ be such a tracial state. Then,
in the proof of 4.2, let $t(a)=\tau\circ \pi\circ \Phi(a).$
\QED

\vspace{0.1in}

{\bf 4.5. Theorem} {\it A unital simple TAI \CA \, has stable rank one.}

{\it Proof:} Let $A$ be a unital simple \CA.
Take a nonzero element $a\in A.$ We will show that
$a$ is a norm limit of invertible elements in $A.$
So we may assume that $a$ is not invertible and $\|a\|=1.$ Since
$A$ is finite, $a$ is not one-sided invertible.
For any $\ep>0,$ by 3.2 in \cite{R1}, there is a zero divisor $b\in A$
such that
$\|a-b\|<\ep/2.$
We further assume that $\|b\|\le 1.$
Therefore, by \cite{R1}, there is a unitary $u\in A$ such that
$ub$ is orthogonal to a nonzero positive element $c\in A.$
Set $d=ub.$
Since $A$ has (SP) (by 3.2), there exists a nonzero projection
$e\in A$ such that
$de=ed=0.$
Since $A$ is simple and has (SP) (by 3.2),
we may write $e=e_1\oplus e_2$ with $e_2\lesssim e_1.$
Note that $d\le (1-e)\le (1-e_1).$
Moreover, $(1-e_1)A(1-e_1)$ is TAI.

Let $\eta>0$ be a positive number. There is a projection
$p\in (1-e_1)A(1-e_1)$ and a unital $C^*$-subalgebra $B\in {\cal I}$ with
$1_B=p$ such that

(1) $\|[x, p]\|<\eta,$

(2) $pxp\in_{\eta} B$ for all $x\in {\cal F},$ and

(3) $[1-e_1-p]\le  [e_2],$

\noindent
where ${\cal F}$ contains $d.$
Thus,
with sufficiently small $\eta,$ we may assume that
$$
\|d-(d_1+d_2)\|<\ep/16,
$$
where $d_1\in B$ and $d_2\in (1-e_1-p)A(1-e_1-p).$

Since \CA s in ${\cal I}$ have stable
rank one and $B\in {\cal I},$ there is an invertible
$d_1'\in B$
such that
$$
\|d_1-d_1'\|<\ep/8.
$$
Let $v$ be a partial isometry such that
$v^*v=(1-e_1-p)$ and $vv^*\le e_1.$ Set $e_1'=vv^*$ and
$d_2'=\ep/8(e_1-e_1')+(\ep/8)v+(\ep/8)v^*+d_2.$
Note that $(\ep/8)v+(\ep/8)v^*+d_2$ has matrix
decomposition
$$
\begin{pmatrix}
0 & \ep/8\\
     \ep/8 & d_2\end{pmatrix}
     $$
Therefore $d_2'$ is invertible in $(1-p)A(1-p).$
This implies that $d'=d_1'+d_2'$ is invertible in $A.$
We also have
$$
\|d_2'-d_2\|<\ep/8,
$$
whence
$$
\|d-d'\|<\|d-(d_1+d_2)\|+\|(d_1+d_2)-(d_1'+d_2')\|<
\ep/16+\ep/8+\ep/8<3\ep/8.
$$
We have
$$
\|b-u^*d'\|\le \|u^*u(b-u^*d')\|=\|ub-d'\|<3\ep/8.
$$
Finally,
$$
\|a-u^*d'\|\le \|a-b\|+\|b-u^*d'\|<\ep/2+3\ep/8<\ep.
$$
Note that $u^*d'$ is invertible.
\QED

\vspace{0.1in}

{\bf 4.6. Corollary} {\it Every unital simple TAI \CA \,
has the cancellation of projections, i.e.,
if $p\oplus e\sim q\oplus e$ then $p\sim q.$ }

\vspace{0.1in}

{\bf 4.7. Theorem} {\it Every unital simple TAI \CA \,
has the following Fundamental Comparability }(\cite{Bl2}):
{\it if $p,\,q\in A$ are two projections with $\tau(p)<\tau(q)$ for
all tracial states $\tau$ on $A,$ then $p\preceq q.$}

{\it Proof:} Denote by $T(A)$ the space of all normalized
traces. It is compact. There is $d>0$ such that
$\tau(q-p)>d$ for all $\tau\in T(A).$
It follows from 3.2 in \cite{LnB} that there exists
a nonzero projection $e\le q$ such that $\tau(e)<d/2$
for all $\tau\in T(A).$ Set $q'=q-e.$
Then $\tau(q'-p)>d/2$ for all $\tau\in T(A).$

It follows from
6.4 in \cite{CP} that there exists a nonzero $a\in A_+$
such that $q'-p-(d/4)=a +z$ and there is a sequence $\{u_n\}$ in $A$
$$
z=\sum_nu_n^*u_n-\sum_nu_nu_n^*.
$$

Choose an integer $N>0$ such that
$$
\|\sum_nu_n^*u_n-\sum_{n=1}^Nu_n^*u_n\|<d/128
\andeqn
\|\sum_nu_nu_n^*-\sum_{n=1}^Nu_nu_n^*\|<d/128.
$$

Let ${\cal F}=\{p, q, q', e, z, u_n, u_n^*, n=1,2,...,N\}$ and
let  $0<\ep<1.$
Since $A$ is TAI, there exists a projection $P\in A$
and a $C^*$-subalgebra $B\in {\cal I}$ with $1_B=P$ such that

(1) $\|[x,P]\|<\ep/2N,$

(2) $PxP\in_{\ep/2N} B$ for all $x\in {\cal F}$ and

(3) $(1-P)\preceq e.$

With sufficiently small $\ep,$ using a standard perturbation
argument, we obtain projections $q''=q_1+q_2,$ $p'=p_1+p_2,$
where $q_1, q_2, p_1, p_2$ are projections, $p_1, q_1\in B$
and $q_2, p_2\in (1-P)A(1-p)$ such that
$$
\|q''-q'\|<d/32\andeqn \|p'-p\|<d/32.
$$
Furthermore (with sufficiently small $\ep$),
we obtain $v_1,v_2,...,v_N\in B$ such that
$$
\|(q_1-p_1-(d/4)P)-(b+\sum_{n=1}^Nv_n^*v_n-\sum_{n=1}^Nv_nv_n^*)\|<d/16,
$$
where $b\in B_+$ and $\|PaP-b\|<\ep/2N.$
Denote by $T(B)$ the space of all normalized traces on $B.$
Then
$$
\tau(q_1-p_1-(d/4)P-b)>-d/16
$$
for all $\tau\in T(B).$
Therefore
$$
\tau(q_1-p_1)>d/4-d/16=3d/16
$$
for all $\tau\in T(B).$
This implies that $p_1\preceq q_1$ in $B,$ whence also in $A.$
Since $p_2\preceq (1-P)\preceq e,$ we conclude that
$$
[p]=[p_1+p_2]\le [q_1]+[e]\le [q].
$$
\QED

{\bf 4.8. Theorem} {\it Let $A$ be a unital simple
TAI \CA. Then $K_0(A)$ is weakly unperforated and satisfies the
Riesz interpolation property.}

{\it Proof:} First we note that, by 3.6, $M_n(A)$ is a unital simple TAI \CA.
To show that $K_0(A)$ is weakly unperforated, it suffices to show that
if $k[p]>k[q]$ for any projections in $M_n(A),$ then
$[p]\ge [q],$ where $k>0$ is an integer.
But $k[p]>k[q]$ implies that
$\tau(p)>\tau(q)$ for all traces. This implies
that $[p]\ge [q]$ by 4.7.
So $K_0(A)$ is  weakly unperforated.

Since $A$ has cancellation, to show that $K_0(A)$ has the Riesz
interpolation property, it suffices to show the following.
If $p\le q$ are two projections in $A$ and $q=q_1+q_2,$
where $q_1$ and $q_2$ are two mutually orthogonal projections,
then $p=p_1+p_2$ with $p_1\le q_1$ and $p_2\le q_2.$
Without loss of generality, we may assume $q-p\not=0.$
Since $q_1Aq_1$ is a TAI \CA, there is a nonzero projection $q_1'\le q_1$
such that  $[p]\le (q_1-q_1')+q_2$ and $q_1-q_1'\not=0.$
Let $q'=(q_1-q_1')+q_2.$
Then $p\preceq q'.$ Without loss of generality, we may assume
that $p\le q'.$

Let ${\cal F}$ be a finite subset containing
$p, q_1-q_1', q_2.$  For any $\ep>0,$ there exists
 a unital $C^*$-subalgebra $B\in {\cal I}$ and
a projection $P\in A$ with $1_B=P$ such that

(1) $\|[x, P]\|<\ep,$

(2) $PxP\in_{\ep} B$ for all $x\in {\cal F}$ and

(3) $(1-P)\preceq q_1'.$

With sufficiently small $\ep,$ without loss of generality, we may assume that
$[p,P]=[q_1-q_1',P]=[q_2,P]=0.$
Write $p''=PpP,$ $q_1''=P(q_1-q_1')P$ and $q_2''=Pq_2P.$
We have $p''\le q_1''+q_2''.$
Note that $M_n$ and $M_n(C([0,1]))$ have the Riesz interpolation
property.
So $B$ has the Riesz property.
There are $p_1'\preceq q_1''$ and $p_2'\preceq q_2''$ such that
$p_1'+p_2'=p''.$
Since $p-p''\le (1-P)\le q_1',$ we let $p_1''=p-p'+p_1'.$
Then $p_1''\preceq q_1'+q_1''\le q_1.$ Now $p=p_1''+p_2'.$
\QED

{\bf 4.9.} Let $A$ be a unital separable simple TAI \CA.
We summarize some of its properties:

(i) $A$ has stable rank one;
(ii) $A$ has at least one tracial state;
(iii) $A$ has Fundamental Comparison property;
(iv) $A$ has weakly unperforated $K_0(A)$ and
satisfies the Riesz interpolation property;
(v) $A$ has property (SP);
(vi) $A$ is MF;
(vii) if $A$ is nuclear, $A$ is also quasidiagonal;
(viii) $M_n(A)$ is TAI;
(ix) Every quasitrace on $A$ is a trace and $T(A)$
is a (metrizable) Choquet simplex;
(x) $A\otimes F$ is TAI for all AF-algebras $F$;
(xi) direct limits of TAI \CA s are TAI and, in fact,
locally TAI \CA s are TAI.

We have not shown (ix). The only thing that one needs to note is
that every quasitrace on \CA s in ${\cal I}$ is in fact a trace.
Then, from condition (3) of Definition 2.2, it is easy to see
that every quasitrace is a trace. Note that it was proved in
\cite{BH} that set of quasitraces on a unital \CA\, is a
Choquet simplex.

We end this section with the following
necessary and sufficient condition for a unital simple
\CA s to be TAI. For the simple case, one could use
it as the definition.

{\bf 4.10. Theorem} {\it Let $A$ be a unital simple \CA.
Then $A$ is TAI if and only if the following hold:

For any finite subset ${\cal F}\subset A$
containing a nonzero element $b,$
$\ep>0,$ integers $n>0$ and $N>0,$ and any nonzero projection $e\in A,$
there exist  a nonzero projection $p\in A$ and
a $C^*$-subalgebra $I=\oplus_{i=1}^k M_{n_i}(C([0,1]),$
with $1_{I}=p$ and $min\{n_i: 1\le i \le k\}\ge N,$
such that

{\rm (1)} $\|[x,p]\|<\ep$ for all $x\in {\cal F},$

{\rm (2)} $pxp\in_{\ep} I$ for all $x\in {\cal F}$
and $\|pbp\|\ge \|b\|-\ep,$ and

{\rm (3')} $1-p$ is unitarily equivalent to a projection in $eAe.$}

{\it Proof:} To show that the above is sufficient for $A$ being
TAI we note that $A$ has property (SP) by 3.2. Then, by 3.2 in
\cite{LnB}, a result of Cuntz, there exists a projection $q\in
eAe$ such that $(n+1)[q]\le [e].$ Then it is clear that  the above
(3') implies (3) in 2.2 (if we use the projection $q$ instead of
$e$).

To see it is also necessary, we use the fact that simple
TAI \CA s have stable rank one (so they have cancellation).
It remains to show that we can make each summand of $I$
have large rank.
But this follows from (the proof of) 3.3.
\QED

\section{Tracially approximately divisibility and \hm s
from \CA\, in {\cal I}}

The main purpose of this section is to prove 5.8 and 5.9.
Theorem 5.8 will be used to prove
theorem 8.6.
Theorem 5.9 classifies monomorphisms from a
\CA\, in ${\cal I}$ to a unital simple TAI \CA.

{\bf 5.1.} Section 1-6 and 8 and most of 9
were written in a 1998 preprint in titled
"Classification of simple TAI \CA s, part I" which was reported
in EU Operator Algebra Conference at Copenhagen in August 1998.
The author later introduced the notation of tracial topological
rank. When $A$ is a unital  simple \CA,
$A$ is an TAI \CA\, if and only if $A$ has tracial topological
rank no more than 1 (see 7.1 in \cite{Lntr1}).

The following is the definition of tracial topological rank no
more than one for simple \CA s:

{\bf 5.2. Definition}
Let $A$ be a unital simple \CA.
Then $A$ has tracial topological rank no more than one
and we will write $TR(A)\le 1$ if the following holds:
For any $\ep>0,$ and any finite subset ${\cal F}\subset A$
containing a nonzero element $a\in A_+,$
there is a \SCA\, $C$ in $A$ where $C=\oplus_{i=1}^k M_{n_i}(C(X_i)),$
where each $X_i$ is a finite CW complex with dimension no more
than one such that $1_C=p$ satisfying the following:

(i) $\|px-xp\|<\ep\,\,\,{\rm for}\,\,\, x\in {\cal F},$

(ii) $pxp\in_{\ep} C \,\,\,{\rm for}\,\,\, x\in {\cal F}$ and

(iii) $1-p$ is equivalent to a projection in ${\overline{ aAa}}.$

In the above definition, if $C$ can be chosen to be a finite dimensional
\SCA\, then we write $TR(A)=0$ (see \cite{Lntr1}).
If $TR(A)\le 1$ but $TR(A)\not=0$ (see \cite{Lntr1})
 then we will write $TR(A)=1.$

In the light of Theorem 7.1 in \cite{Lntr1}, {\it in what follows,
we will replace unital simple TAI \CA s by unital simple \CA s
with tracial topological rank no more than one and write
$TR(A)\le 1.$ }

\vspace{0.1in}

{\bf 5.3. Definition}
A unital simple \CA\, $A$ is said to be {\it tracially approximately
divisible} if for any $\ep>0,$ any projection $e\in A,$
any integer $N>0$ and
any finite subset ${\cal F}\subset A,$ there exists a
projection $q\in A$ and there exists a finite dimensional
$C^*$-subalgebra $B$ with each simple summand having rank at least $N$
such that

(1) $\|qx-xq\|<\ep$ for all $x\in {\cal F},$

(2) $\|y(1-q)x(1-q)-(1-q)x(1-q) y]\|<\ep$ for all $x\in {\cal F}$ and
all $y\in B$ with $\|y\|\le 1,$ and

(3) $q$ is unitarily equivalent to a projection of $eAe.$

Of course if $A$ is approximately divisible, then $A$ is
tracially approximately divisible (see \cite{BKR}).

\vspace{0.1in}

{\bf  5.4. Theorem} {\it Every nonelementary unital simple
\CA\, with $TR(A)\le 1$  is tracially approximately divisible.}

{\it Proof:} Let $A$ be a unital simple \CA\, with $TR(A)\le 1.$
Fix $\ep>0,$ $\sigma>0,$ $N>0$  and a finite subset ${\cal F}\subset A.$
Let $b\in A$ with $\|b\|=1$ and assume that $b\in {\cal F}.$
There exist a projection $p\in A$ and a $C^*$-subalgebra $C\in {\cal I}$
with $1_C=p$ such that

(1) $\|px-xp\|<\ep/4$ for all $x\in {\cal F},$

(2) $pxp\in_{\ep/4} C$ and $\|pbp\|\ge \|b\|-\ep/2,$ and

(3) $\tau(1-p)<\sigma/2$ for all traces $\tau$ on $A.$

Write $C=\oplus C_i,$ where $C_i=M_{l(i)}(C[0,1]),$ or
$C_i=M_{l(i)}.$
It will become clear that, without loss of generality, to simplify notation,
we may assume that $C=C_i$ (i.e., there is only one summand).
If  $C=M_l,$ let $\{e_{ij}\}$
be matrix units for $M_l.$ Since $A$ is not elementary,
there is a positive element $a\in e_{11}Ae_{11}$ such that
$sp(a)=[0,1]$ (see p.6.1 in \cite{AS}). This implies that
$C\subset M_l(C([0,1])).$ So, we may assume that
$C=M_l(C([0,1])).$
Let ${\cal G}_1\subset C$ be a finite subset such that
$$
{\rm dist}(pxp, {\cal G}_1)<\ep/4
$$
for all $x\in {\cal F}.$
Let ${\cal G}$ be a finite subset of $C$
containing $\{e_{ij}\}$ and $e_{ij}ge_{ij}^*$ for all
$g\in {\cal G}_1.$

Let $\eta>0.$ Denote by  $\dt$  the positive number in Theorem 4.3 of
\cite{Li} corresponding to $\eta$ (instead of $\ep$).
Let $\{f_1, f_2,...,f_m\}\subset C([0,1])$ be as in 3.1 with
respect to $\dt$ ($=d$). We identify $C([0,1])$ with $e_{11}Ce_{11}.$
Since $e_{11}Ae_{11}$ is simple, there are $b_{ij}\in e_{11}Ae_{11}$
such that
$$
\|\sum_{i=1}b_{ij}f_jb_{ij}^*-e_{11}\|<1/16,
$$
$j=1,2,...,m.$
Let  ${\cal G}_2$ be a finite subset containing
$\{f_j, b_{ij}, b_{ij}^*\}\bigcup\{ a_{ij}\in e_{11}Ae_{11}: (a_{ij})_{l\times l}\in
{\cal G}\}.$

By 3.4,  $TR(e_{11}Ae_{11})\le 1.$
So for any $0<\sigma<\eta/2$ and any finite subset
${\cal G}_3\supset {\cal G}_2,$
there exist a projection $q\in e_{11}Ae_{11}$ and
a $C^*$-subalgebra $C_1\subset e_{11}Ae_{11}$ with
$1_{C_1}=q$ and $C_1\in {\cal I}$ satisfying the following:

(a) $\|qx-xq\|<\sigma,$

(b) $qxq\in_{\sigma} C_1$
for all $x\in {\cal G}_3,$

(c) $\tau(e_{11}-q)<\sigma/2l$ for all traces $\tau.$

With sufficiently small $\sigma$ and sufficiently large ${\cal G}_2,$
we may assume that there exists a \hm\, $\phi: C([0,1])\to C_1$ such that

(b') $\|\phi(x)-qxq\|<\eta/2$ for all $x\in {\cal G}_2\cap C([0,1]).$

Note that we also have $c_{ij}\subset C_1$ such that
$$
\|\sum_{i=1}c_{ij}\phi(f_j)c_{ij}^*-q\|<1/8, \,\,\,\,j=1,2,...,m.
$$

We are now applying Theorem 4.3 in \cite{Li}.
It follows from 3.1 that $Sp(\phi_t)$ is $\dt$-dense in $[0,1].$
By applying 4.3 in \cite{Li}, there is a \hm\,
$\psi: C([0,1])\to C_1$ and
there is a finite dimensional $C^*$-subalgebra $F=\oplus_i F_i,$
where each $F_i$ is  simple  and ${\rm dim} F_i\ge N,$ with
$1_{F}=q$ such that
$$
\|\psi(f)-\phi(f)\|<\eta/2\,\,\,\,\, {\rm for\,\,\, all}\,\,\,
f\in {\cal G}_2 \andeqn
$$
$$
\|[\psi(g), b]\|=0
$$
for all $g\in C([0,1])$ and $b\in F.$ Set $F'=diag(F,F,...,F)$ in
$F\otimes M_l,$ $\psi'=\psi\otimes {\rm id}_{M_l},$
$\phi'=\phi\otimes {\rm id}_{M_l}$ and $P=diag(q,q,...,q)\in
M_l(C_1).$ With sufficiently small $\eta$ and large ${\cal G}_2,$
we have
$$
\|\psi'(g)-\phi'(g)\|<\ep/2\,\,\,{\rm for}\,\,\, g\in {\cal G}.
$$
We also have
$$
\|[\psi'(f),c]\|=0\,\,\,\,\,\,\,{\rm for}\,\,\,f\in C\andeqn c\in F'.
$$
These imply that
$$
\|[PxP, c]\|<\ep\,\,\,\,\,{\rm for\,\,\, all}\,\,\, x\in {\cal F}
\andeqn c\in F'.
$$
Note that $1_{F'}=P.$
We also have
$$
\tau(1-P)\le \sigma/2+l\sigma/2l=\sigma.
$$
By 4.7,
we conclude that $A$ is tracially approximately divisible.
When $C=\oplus C_i,$ it is clear that we can do exactly the same as above
for each summand. Let $d_i=1_{C_i}.$
If we find a matrix algebra $F_i\in d_iAd_i$ with rank greater than $N$
which commutes with $C_i,$ then $\oplus F_i$ commutes with $C.$
\QED

{\bf 5.5. Lemma} {\it Let $A$ be a unital nuclear simple
\CA\, with $TR(A)\le 1.$ Then for any $\ep>0,$ any $\sigma>0,$ any
integer $n>0,$ and any finite subset ${\cal F}\subset A,$
there exist mutually orthogonal projections $q, p_1,p_2,...,p_n$
with $q\preceq p_1$ and  $[p_1]=[p_i]$ ($i=1,2,...,n$),
a $C^*$-subalgebra $C\in {\cal I}$ with $1_C=p_1$
and completely positive linear contractions
$L_1: A\to qAq$ and $L_2:A\to C$ such that
$$
\|x-(L_1(x)\oplus diag(L_2(x),L_2(x),...,L_2(x)))\|<\ep \andeqn
$$
$$
\|L_i(xy)-L_i(x)L_i(y)\|<\ep,
$$
where $L_2(x)$ is repeated $n$ times,
for all $x,\,y\in {\cal F}$ and
$\tau(q)<\sigma$ for all $\tau\in T(A).$}

{\it Proof:}
From the proof of 5.4, we have the following:
For any $\eta>0,$ any integer $K>0,$ any integer $N>4Kn^2$
and finite subset ${\cal G}\subset
A$ (containing $1_A$),
there exists a projection $P\in A$ and a finite dimensional
$C^*$-subalgebra $B$ with $1_B=P$ such that

(i) $\|[P,x]\|<\eta$ for all $x\in {\cal G};$

(ii) every simple summand of $B$ has rank at least $N;$

(iii) there is a $C^*$-subalgebra
$D\in {\cal I}$ with $1_D=P$
such that $[d,g]=0$ for all $d\in D,$ $g\in B$ and
$$
{\rm dist}(x,D)<\eta\,\,\, {\rm for}\,\,\,x\in {\cal G};
\andeqn
$$

(iv) $5N[(1-P)] <[P].$

Let ${\cal F}_1\subset A$ be a finite subset (containing $1_A$)
and $\sigma>0.$
Since $A$ is nuclear, with sufficiently large ${\cal G}$ and
sufficiently small $\eta,$ by 3.2 of \cite{LnC},
there are unital completely positive linear contractions
$L_1': A\to (1-P)A(1-P)$ and $L_2': A\to D$ such that
$L_1'(a)=(1-P)a(1-P),$
$$
\|x-L_1'(x)\oplus L_2'(x)\|<\sigma
\andeqn
\|L_2'(x)-PxP\|<\eta+\sigma
$$
for all $x\in {\cal F}_1.$
It follows that, with sufficiently small $\sigma$ and $\eta,$
$$
\|L_i'(xy)-L_i'(x)L_i'(y)\|<\ep
$$
for all $x,y\in {\cal F}_1.$ Write $B=\oplus_{i=1}^kB_i,$ where
$B_i\cong M_{l(i)}$ with $l(i)\ge N,$ and denote by $C$ the
$C^*$-subalgebra generated by $D$ and $B.$ Note that
$C\cong\oplus_{i=1}^kD_0\otimes B_i,$ where $D_0\cong D.$ Let
$\pi_i: C\to D_0\otimes B_i$ be the projection. Denote
$\phi_i=\pi_i\circ L_2'.$ By (iii), we see that we may write
$\phi_i=diag(\psi_i,\psi_i,...,\psi_i),$ where $\psi_i: A\to
e_i(D_0\otimes M_{l(i)})e_i$ and $e_i$ is a  minimal rank-one
projection of $M_{l(i)}.$ Write $l(i)=k_in+r_i,$ where $k_i\ge
n>r_i\ge 0$ are  integers. We may rewrite
$$
\phi_i=diag(\Phi'_i, ...,\Phi'_i)\oplus \Psi'_i,
$$
where $\Phi'_i=diag(\psi_i,...,\psi_i): A\to D_0\otimes M_{k_i}$
is repeated $n$ times and
$\Psi_i'=diag(\psi_i,...,\psi_i): A\to D_0\otimes M_{r_i}.$

Define
$L_2=\oplus_{i=1}^k\Phi_i'$ and $L_1=L_1'\oplus_{i=1}^k\Psi_i'.$
We estimate that
$$
\tau((1-P)+\oplus_{i=1}^k\Psi_i'(1_A))
< (1/5N)\tau(P)+(1/4nK)\tau(P)<(1/2n)\tau(P)\le
\min(\sigma,\tau([L_2(1_A)]),
$$
provided that $1/K<\sigma.$
By 4.7, the lemma follows.
\QED

\vspace{0.1in}




The following follows from Lemma 5.5 immediately.

{\bf 5.6. Corollary}
{\it Let $A$ be a unital separable simple \CA\,
$TR(A)\le 1.$
Then for any $\ep>0,$ any $\sigma>0,$ any
integer $n>0,$ and any finite subset ${\cal F}\subset A,$
there exists a \SCA\, $C\in {\cal I}$
such that
$$
\|x-(1-p)x(1-p)\oplus {\rm diag}(y,y,...,y)\|<\ep
$$
where $y\in C$ and ${\rm diag}(y,y,...,y)\in M_n(C)$
and $p=1_{M_n(C)}$
for all $x,\in {\cal F}$ and
$\tau((1-p))<\sigma$ for all $\tau\in T(A).$
Moreover, we may require that $\|(1-p)x(1-p)\|\ge (1-\ep)\|x\|$
for all $x\in {\cal F}.$}

{\it Proof}:
Perhaps the last part of the statement needs an explanation.
In the proof of 5.5, we know that we may require
that $\|y\|\ge (1-\ep/2)\|x\|$ for all $x\in {\cal F},$
Thus we may replace $(1-p)x(1-p)$ by $(1-p)x(1-p)\oplus y$ and
replace $(1-p)$ by $1-p\oplus {\rm diag}(1_C,0,...,0).$
\QED

{\bf 5.7. Lemma} {\it Let
$B=\oplus_{i=1}^kB_i$ be a unital \CA\, in ${\cal I}$ (
where $B_i$ is a single summand).
For any $\ep>0,$ any finite subset ${\cal F}\subset B$ and any integer $L>0,$
there exist  a finite subset ${\cal G}\subset B$ depending on
$\ep$ and ${\cal F}$ but not $L,$ and $\dt=1/4L$
such that the following  holds.
If $A$ is a unital separable nuclear simple \CA\, with $TR(A)\le 1$
and $\phi_i: B\to A$ are
two \hm s satisfying the following:

{\rm (i)} there are $a_{g,i}, b_{g,j}\in A, i,j\le L$ with
$$
\|\sum_i a_{g,i}^*\phi_1(g)a_{g,i}-1_A\|<1/16
\andeqn
\|\sum_jb_{g,j}^*\phi_2(g)b_{g,j}-1_A\|<1/16
$$
for all $g\in {\cal G};$

{\rm (ii)} $(\phi_1)_*=(\phi_2)_*$ on $K_0(B);$ and,

{\rm (iii)} if $\|\tau\circ\phi_1(g)-\tau\circ \phi_2(g)\|<\dt$
for all $g\in {\cal G},$
then there exists a unitary $u\in A$ such that
$$
\|\phi_1(f)-u^*\phi_2(f)u\|<\ep
\,\,\,\rforal \,\,\,f\in {\cal F}.
$$}

{\it Proof:}
It is clear that we can reduce the general case to
the case in which $B$ has only one summand. Since the case in which $B=M_{l(i)}$
is well known to hold, we may assume
that $B=M_l(C([0,1])).$
Fix any $d_0>0.$ Condition (i), with sufficiently large ${\cal G},$ implies
that $Sp(\phi_i)$ is $d_0$-dense in $[0,1]$ ($i=1,2$).
By the proof of 3.7, therefore, for any $d_1>d_0,$ we may
assume that
$$
\phi_i(f)=\phi_i'(f)\oplus \sum_{j=1}^{N_i}f(t_{(i,j)})q(i,j)\,\,\,\,\,i=1,2,
$$
where $\{t_{(i,1)},t_{(i,2)},...,t_{(i,N_i)}\}$ ($i=1,2$)
 is $d_1$-dense in $[0,1]$ and
$q(i,1),...,q(i,N_i)$ are mutually (non-zero) orthogonal projections in $A.$
It is clear that without loss of generality, we may assume that
$t_j=t_{i,j}$ and $N_1=N_2.$
Since $TR(A)\le 1,$ we can find nonzero projections
$q(i,j)'\le q(i,j)$ such that
$q(1,j)'$ and $q(2,j)$ are unitarily equivalent.
By replacing $\phi_1$ by ${\rm ad}\,z\circ \phi_1$ for
some unitary $z,$ we may assume that $q(1,j)'=q(2,j)'.$
Then, by replacing $\phi_i'$ by $\phi_i'\oplus \sum_{j=1}^{N_i}f(t_j)(q(i,j)-
q(i,j)'),$ we may assume that, with $q_j=q(i,j)'$ and $N=N_1,$
$$
\phi_i(f)=\phi_i'(f)\oplus \sum_{j=1}^{N}f(t_j)q_j \,\,\,\,(i=1,2).
$$

Now let $Q=1-\sum_{j=1}^N q_j.$
We will apply 5.14 in \cite{Li}.
To do this, we let $r>0$ be as in the statement of 5.14 in \cite{Li} (but
with respect to $\ep/4$ and
${\cal L}_r\subset B$ (see also 5.2 in \cite{Li})).
Let $d=1/r$ and $\dt=1/4L.$
Set ${\cal G}\subset {\cal L}_r$ such that  the functions
$f_1,f_2,...,f_m$ required in 3.1 are all in ${\cal G}.$

Fix an integer $n>1.$ Let $e'\le QAQ$ such that $n[e']\le [q_i],$
$i=1,2,...,N.$ This is possible since $A$ is simple and has (SP).
Let $\eta>0,$ $e\in QAQ$ be any nonzero projection in $A$
with $\tau(e)<1/2L,$ $[e]\le [e']$ and $K>0$ be an
integer.
Since $A$ is a unital simple \CA\, with $TR(A)\le 1,$ there exist a projection
$P\in A$  and a unital $C^*$-subalgebra $C\in {\cal I}$ with
$1_C=P$ such that

(i) $\|[\phi_i'(g), P]\|<\eta,$

(ii) $P\phi_i(g)P\in_{\eta} C$ for all $g\in {\cal G}'$ and $i=1,2,$ and

(iii) $K[Q-P]\le [P]$ and $[Q-P]\le [e],$

where
${\cal G}'\supset {\cal G}\bigcup \{a_{g,i}, a_{g,i}^*, b_{g,j}, b_{g,j}^*,
g\in {\cal G} \andeqn i,j\le L\}.$
For any $\sigma>0,$ with sufficiently small $\eta,$
there exists \hm\, $\psi_i: B\to C$ such that
$$
\|\psi_i(g)-P\phi_i'(g)P\|<\sigma,
$$
$$
\|\sum_{i=1}c_{g,i}^*\psi_1(g)c_{g,i}-P\|<1/8
\andeqn
\|\sum_{j=1}d_{g,j}^*\psi_2(g)d_{g,j}-P\|<1/8
$$
for all $g\in {\cal G},$ where $c_{g,i},\,d_{g,j}\in C.$
It follows from 3.1 that
$$
|sp((\psi_i)_t)\cap T|\ge 1/L|Sp((\psi_i)_t)|
$$
for all $t\in [0,1]$ (or, both $\psi_1$ and $\psi_2$ have the property
$sdp(r,1/L)$ as in  5.13 in \cite{Li}), where $T$ has length at least
$1/r.$
We also have, if $\eta$ and $\sigma$ are  sufficiently small,
$$
\|\tau\circ \psi_1(g)-\tau\circ \psi_2(g)\|<1/2L
$$
for all $g\in {\cal G}.$
It follows from 5.14 in \cite{Li} that there exists a unitary
$v\in C$ such that
$$
\|\psi_1(f)-v^*\psi_2(g)v\|<\ep/4
\,\,\,\rforal \,\,\,f\in {\cal F}.
$$

We also have
$$
\|\phi_i'(f)-\psi_i(f)-(Q-P)\phi_i'(f)(Q-P)\|<\ep/4
$$
for all $f\in{\cal F}.$
Hence,
$$
\|\phi_i(f)-(\sum_{j=1}^Nf(t_j)q_j\oplus\psi_i'(f)\oplus
(Q-P)\phi_i''(f)(Q-P)\|<\ep/2
$$
for all $f\in {\cal F}.$
Since
$$
n[Q-P]\le [q_i],
$$
by, for example, Lemma 8 (i) in \cite{LR} (this was known earlier),
there exists (provided that $n$ is sufficiently large and
$d_1$ is sufficiently small, and these two numbers do not
depend on $\phi_i$ or $A$) a unitary $w\in
(1-P)A(1-P)$ such that
$$
\|\sum_{j=1}^Nf(t_j)q_j\oplus (Q-P)\phi_1''(f)(Q-P)
-w^*(\sum_{j=1}^Nf(t_j)q_j\oplus (Q-P)\phi_2''(f))(Q-P)w\|<\ep/2
$$
for all $f\in {\cal F}.$
Thus,
we obtain a unitary $u\in A$ such that
$$
\|\phi_1(f)-u^*\phi_2(f)u\|<\ep
\,\,\,\rforal\,\,\, f\in {\cal F}.
$$ \QED

 \vspace{0.1in}

{\bf 5.8. Theorem} {\it Let $A$ be a unital simple \CA\, with $TR(A)\le 1$
and $C$ be a $C^*$-subalgebra of $A$ in ${\cal I}.$
Then for any finite subset ${\cal F}\subset C$
and $\ep>0,$ there exist $\dt>0, $ $\sigma>0$ and a finite subset
${\cal G}\subset A$ satisfying the following:
if  $L_1, L_2:  A\to B$ are  two unital ${\cal G}$-$\dt$-multiplicative
\morp s, where $B$ is a unital simple \CA\, with $TR(B)\le 1,$ with
$(L_1|_{C})_*=(L_2|_C)_*$ on $K_0(C)$ and
$$
|\tau(L_1(g))-\tau\circ L_2(g)|<\sigma
$$
for all $g\in {\cal G}$ and for all $\tau\in T(B),$ then there
is a unitary $u\in A$ such that
$$
\|L_1(f)-u^*L_2(f)u\|<\ep
\,\,\,\rforal\,\,\,f\in {\cal F}.
$$}

{\it Proof:} Fix $\ep>0$ and a finite subset
${\cal F}\subset A.$
Let ${\cal G}_1\subset C$ be the finite subset
required by 5.7 (for a given $\ep>0$ and a given finite subset ${\cal F}$).
Suppose that $a_{g,i}\in A$ such that
$$
\|\sum_{i=1}^{n(g)}a_{g,i}^*ga_{g,i}-1_A\|<1/64
$$
for all $g\in {\cal G}_1.$
Set $L=\max\{n(g): g\in {\cal G}\}.$
Then, with sufficiently small $\dt>0$ and
large ${\cal G}\supset {\cal G}_1\bigcup\{a_{g,i}: g, i\},$
we have $b_{g,i,j}\in A$ such that
$$
\|\sum_{i=1}^{n(g)}b_{g,i,j}^*L_j(g)b_{g,i,j}-1_B\|<1/32
$$
for all $g\in {\cal G}_1$ and $j=1,2.$
Furthermore, for any $\eta>0,$ with sufficiently small $\dt,$
there is a \hm\, $\phi_j: C\to B$ ($j=1,2$) such that
$$
\|\phi_j(g)-L_j(g)\|<\eta
\andeqn
\|\sum_{i=1}^{n(g)}b_{g,i,j}^*\phi_j(g)b_{g,i,j}-1_B\|<1/16.
$$
for $g\in {\cal G}_1.$
We also require that $\sigma<1/4L.$
Then we see the conclusions of the theorem follow from
5.7 (and its proof) immediately.
\QED
\vspace{0.1in}

{\bf 5.9. Theorem} {\it Let $A$ be a unital simple \CA\, with $TR(A)\le 1$
and $B\in {\cal I}.$
Let $\phi_i: B\to A$ be two monomorphisms such that
$$
(\phi_1)_*=(\phi_2)_*: K_0(B)\to K_0(A)
\,\,\,and\,\,\,
\tau\circ\phi_1=\tau\circ\phi_2
$$
for all $\tau\in T(A).$
Then there is a sequence of unitaries $u_n\in A$ such that
$$
\lim_{n\to\infty}u_n^*\phi_1(x)u_n=\phi_2(x)
\,\,\,for\,\,\,all\,\,\,
x\in B.
$$}

{\it Proof:} As before, we reduce the general case
to the case in which $B=C([0,1]).$
Let $\ep>0$ and ${\cal F}\subset B$ be a finite subset.
Let ${\cal G}\subset B$ be the finite subset in the statement
of 5.7 (it does not depend on $L$).
Since $A$ is simple, there exists an integer $L>0$ and
$a_{i,g}, b_{i,g}\in A,$ $i=1,2,...,L$ (some of them could be zero)
such that
$$
\|\sum_ia_{i,g}^*\phi_1(g)a_{i,g}-1\|<1/16
\andeqn
\|\sum_ib_{i,g}^*\phi_2(g)b_{i,g}-1\|<1/16
$$
for all $g\in {\cal G}.$
Therefore the theorem follows from 5.7.
\QED

\section{The unitary group of a simple \CA\, $A$  with $TR(A)\le 1$ }

We start with the following observation:

{\bf 6.1}
Let $A$ be a unital \CA \, and $p, a\in A.$
Suppose that $p$ is a projection, $\|a\|\le 1$ and
$$
\|a^*a-p\|<1/16 \andeqn \|aa^*-p\|<1/16.
$$
A standard computation shows
that
$$
\|pap-ap\|<3/16 \andeqn \|pa-pap\|<3/16.
$$
Also
$
\|pa-a\|<1/2.
$
Set $b=pap.$ Then
$$
\|b^*b-p\|\le \|pa^*ap-pa^*a\|+\|pa^*a-p\|<1/16+1/16=1/8.
$$
So
$$
\|(b^*b)^{-1}-p\|<{1/8\over{1-1/8}}=1/7 \andeqn \||b|^{-1}-p\|<2/7,
$$
where the inverse is taken in $pAp.$
Set $v=b|b|^{-1}.$
Then $v^*v=p=vv^*$ and
$$
\|v-b\|<2/7.
$$
We denote
$v$ by ${\tilde a}.$
Suppose that $L: A\to B$ is a ${\cal G}$-$\dt$-multiplicative
\morp, $u$ is a normal partial isometry and a projection
$p\in B$ is given so that
$$
\|L(u^*u)-p\|<1/32.
$$
Note if $v'$ is another unitary in $pAp$ with
$\|v'-b\|<1/3,$ then $[v']=[v]$ in $U(pAp)/U_0(pAp).$
We define  ${\tilde L}$ as follows.
Let $L(u)=a.$ With small $\dt$ and large ${\cal G}, $
we denote by ${\tilde L}(u)$ the normal partial isometry
(unitary in a corner)
$v$ defined above. This notation will be used later.
Note also, if $u\in U_0(A),$ then, with sufficiently large
${\cal G}$ and sufficiently small $\dt,$ we may assume
that ${\tilde L}(u)\in U_0(B).$

\vspace{0.1in}

{\bf 6.2 Definition}
Let $A$ be a unital \CA. Let $CU(A)$ be the {\it closure} of
the commutator subgroup of $U(A).$
Clearly that the commutator subgroup  forms a normal subgroup of $U(A).$
It follows that $CU(A)$ is a normal subgroup of $A.$
It should be noted that $U(A)/CU(A)$ is commutative.
It is an easy fact that if $A=M_r(C(X)),$ where $X$ is a finite
CW complex of dimension 1, then $CU(A)\subset U_0(A).$
If $K_1(A)=U(A)/U_0(A),$ it is known and easy to verify
that every commutator is in $U_0(A).$ Therefore $CU(A)\subset U_0(A).$
If $u\in U(A),$ we will use ${\bar u}$ for the image of $u$ in $U(A)/CU(A),$
and if $F\subset U(A)$ is a subgroup of $U(A),$ then
${\bar F}$ is the image of $F$ in $U(A)/CU(A).$

If ${\bar u}, {\bar v}\in U(A)/CU(A)$ define
$$
{\rm dist}({\bar u}, {\bar v})
=\inf\{\|x-y\|: x,y\in U(A)\,\,\, {\rm such\,\,\, that}\,\,\,
{\bar x}={\bar u}, \, {\bar y}={\bar v}\}.
$$
If $u, v\in U(A)$ then
${\rm dist}({\bar u}, {\bar v})
=\inf\{\|uv^*-x\|: x\in CU(A)\}.
$
Let $g=\prod_{i=1}^n a_ib_ia_i^{-1}b_i^{-1},$ where $a_i, b_i\in U(A).$
Let ${\cal G}$ be a finite subset of $A,$ $\dt>0$ and
$L: A\to B$ be a ${\cal G}$-$\dt$-multiplicative \morp, where $B$ is a
unital \CA. From 6.1, for $\ep>0,$ if ${\cal G}$ is sufficiently large
and $\dt$ is sufficiently small,
$$
\|L(g)-\prod_{i=1}^n a_i'b_i'(a_i')^{-1}(b_i')^{-1}\|<\ep/2,
$$
where  $a_i', b_i'\in U(B).$
Thus, for any $g\in CU(A),$ with sufficiently large ${\cal G}$ and
sufficiently small $\dt,$
$$
\|L(g)-u\|<\ep
$$
for some $u\in CU(B).$ Moreover, for any finite subset
${\cal U}\subset U(B)$ and subgroup
$ F\subset U(B)$ generated by ${\cal U},$
and $\ep>0,$ there exists a finite subset ${\cal G}$ and
$\dt>0$ such that, for any ${\cal G}$-$\dt$-multiplicative
\morp\, $L: A\to B,$
$L$ induces a \hm\, $L^{\ddag}: {\bar F}\to U(B)/CU(B)$ such that
${\rm dist}(\overline{{\tilde L}(u)}, L^{\ddag}({\bar u}))<\ep$
for all $u\in {\cal U}.$
Note we may also assume that ${\bar F}\cap U_0(A)/CU(A)\subset
U_0(B)/CU(B).$

If $\phi: A\to B$ is  a \hm\, then
$\phi^{\ddag}: U(A)/CU(A)\to U(B)/CU(B)$ is  the induced
\hm. It is continuous.

\vspace{0.1in}

{\bf 6.3. Lemma} (N. C. Phillips)
{\it Let $A$ be a unital \CA\, and $2> d>0.$
Let $u_0, u_1,...,u_n$ be $n$ unitaries in $A$
such that
$$
u_n=1_A\andeqn
\|u_i-u_{i+1}\|\le d,\,\,\,i=0,1,...,n-1.
$$
Then there exists a unitary $v\in M_{2n+1}(A)$
with exponential length no more than $2\pi$
such that
$$
\|(u_0\oplus 1_{M_{2n}}(A))-v\|\le  d.
$$
Moreover, $v$ can be chosen in $CU(M_{2n+1}(A)).$}

The following is another version of the above lemma.

\vspace{0.1in}

{\bf 6.4. Lemma}
{\it Let $A$ be a unital \CA\, and $u\in U_0(A).$
Then for each $L>0,$ there is an integer $N(L)>0$
such that if $u=v\oplus (1-p)$ and $v\in U_0(pAp)$
with $cel(v)\le L$ in $pAp$ and
there are $N$ mutually orthogonal and mutually equivalent
projections in $(1-p)A(1-p)$ each of which is equivalent
to $p$ then
$$
cel(u)<2\pi+\ep
$$
for any $\ep>0.$
To be more precise, if $n>2L,$ then
$
cel(u)\le 2\pi+(L/n)\pi.
$
Furthermore, there is a unitary $w\in CU(A)$ such that
$
cel(uw)<(L/n)\pi.
$
 }

(See the proof of Theorem 3.8 in \cite{Ph3} and also
Cor. 5 in \cite{Ph1}.
It should be noted that a unitary in $M_2(A)$ with the form
${\rm diag}(u,u^*)$ is in $CU(M_2(A)).$
)

{\bf 6.5 Theorem }
{\it Let $A$ be a unital simple \CA\, with $TR(A)\le 1.$
Let $u\in U_0(A).$
Then, for any $\ep>0,$
there are unitaries $u_1, u_2\in A$ such that
$u_1$ has exponential length no more than $2\pi,$
$u_2$ is an exponential and
$$
\|u-u_1u_2\|<\ep.
$$
Moreover, $cer(A)\le 3+\ep.$
}

{\it Proof:} Let $\ep$ be a positive number.
Let $v_0, v_1,...,v_n\in U_0(A)$ such that
$$
v_0=u, v_n=1\andeqn \|v_i-v_{i+1}\|<\ep/16,\,i=0,1,...,n-1.
$$
Let $\dt>0.$
Since $TR(A)\le 1,$ there exists a projection
$p\in A$ and a unital $C^*$-subalgebra $B\subset A$ with $B\in {\cal I}$ and
with
$1_B=p$ such that

(1) $\|[v_i, p]\|<\dt,$ $i=0,1,...,n,$

(2) $pv_ip\in_{\dt} B,$ $0,1,...,n,$ and

(3) $2(n+1)[1-p]\le p.$

There are unitaries $w_i\in (1-p)A(1-p)$ with $w_n=(1-p)$ such that
$$
\|w_i-(1-p)v_i(1-p)\|<\ep/16,\,\,\,\,i=0,1,...,n
$$
for any given $\ep>0,$ provided $\dt$ is sufficiently small.
Furthermore, there is a unitary $z\in B$ such that
$$
\|z-pup\|<\ep/16.
$$
Therefore  (with $\dt<\ep/32$)
$$
\|u-w_1\oplus z\|<\ep/8.
$$

Write $z_1=w_1\oplus p.$ Since $2(n+1)[1-p]\le p,$
by 6.3, there is a unitary $u_1$ with exponential length
no more than $2\pi$ such that
$$
\|z_1-u_1\|<\ep/4.
$$
Now since $z\in B$ and it is well known that $B$ has exponential
rank $1+\ep,$ there is an exponential $u_2\in A$ such that
$$
\|u_2-(1-p)-z\|<\ep/3.
$$
Therefore
$$
\|u-u_1u_2\|<\ep.
$$
Since $cel(u_1)\le 2\pi,$ it follows from \cite{Rn} that
$cer(u_1)\le 2+\ep.$ Therefore $cer(u)\le 3+\ep.$ So $cer(A)\le
3+\ep.$ \QED

\vspace{0.1in}

{\bf 6.6. Lemma}
{\it Let $A$ be a unital \CA.

{\rm (1)} $U_0(A)/CU(A)$ is divisible.

{\rm (2)} If $u\in U(A)$ such that $u^k\in U_0(A).$
Then there is $v\in U_0(A)$ such that
${\bar v^k}={\bar u^k}$ in $U(A)/CU(A).$

{\rm (3)} Suppose that $K_1(A)=U(A)/U_0(A)$ and $G\subset
U(A)/CU(A)$ is finitely generated subgroup. Then one has $ G=G\cap
(U_0(A)/CU(A))\oplus \kappa(G), $ where $\kappa: U(A)/CU(A)\to
U(A)/U_0(A)$ is the quotient map. }

{\it Proof}:
Let $u\in U_0(A).$ Then there are $a_1,a_2,...,a_n\in A_{s.a.}$ such that
$
u=\prod_{j=1}^n exp(i a_j).
$
For any integer $k>0,$ let $v=\prod_{j=1}^n exp(i a_j/k).$
Then ${\bar v}^k={\bar u}.$
This proves (1).

To see (2), put $u^k=\prod_{j=1}^nexp(ia_j),$ where $a_j\in A_{s.a}.$
Let $v=\prod_{j=1}^nexp(ia_j/k).$
Thus $\overline{(uv^*)^k}={\bar 1}.$  So ${\bar v}^k={\bar u}^k.$
To see (3), we note that (1) implies
$0\to U_0(A)/CU(A)\to G+U_0(A)/CU(A)\to \kappa(G)\to 0$ splits.
\QED


\vspace{0.1in}

{\bf 6.7. Theorem} {\it Let $A$ be a unital simple  \CA\, with
$TR(A)\le 1$ and
$e\in A$ be a projection.
Let $\gamma: U(eAe)/CU(eAe)\to
U(A)/CU(A)$ be defined by $\gamma(\bar u)=\overline{u\oplus
(1-e)}.$ Then $\gamma$ is  a surjective (contractive)  \hm.
}

{\it Proof}:
It is clear that $\gamma$ is a \hm\, and is contractible.
We will show that $\gamma$ is also surjective.
Fix $u\in U(A).$
Let $N>0$ be an integer such that $N[e]\ge 1$ in $A.$
Fix $1/2>\ep>0$ and $0<\eta<\ep/8(N+1).$
It follows from 5.6 that there is a unitary
$z_1=s_0\oplus s_1\oplus s_1\cdots \oplus s_1,$
where $s_0\in U((e_0Ae_0))$ and $s_1\in U(C)$ which repeats $n+1$ times
($n\ge 3$),  where $C\in {\cal I}$ and $e_0\oplus 1_C$
is equivalent to a subprojection of $e$
 such that
$$
\|u-z_1\|<\eta/4.
$$
Note that $M_{n+1}(C)$ is a \SCA\, of $A.$
By replacing $s_0$ by $s_0\oplus s_1\oplus \cdots \oplus s_1,$ where
$s_1$ repeats several times, we may assume that $3\le n\le 4N+1.$
Without loss of generality, we may also assume that
$e_0\oplus 1_C\le e.$
Let $w=e_0\oplus s_1^n\oplus s_1^*\oplus s_1^*\oplus\cdots \oplus s_1^*,$
where $s_1^*$ repeats $n$ times.
Then $z_1w=s_0\oplus s_1^{n+1}\oplus 1_C\oplus 1_C\cdots \oplus 1_C.$
Put $v_1=s_0\oplus s_1^{n+1}\oplus [e-(e_0\oplus 1_C)]$ and put
$y=s_1^n\oplus s_1^*\oplus s_1^*\oplus\cdots \oplus s_1^*$
($s_1^*$ repeats $n$ times).
Then  $\det(y)=1$ (in $M_{n+1}(C)$).
Since $U_0(M_{n+1}(C))=U(M_{n+1}(C)),$ it follows from 2.4 in \cite{Tm1} that
$y\in CU(M_{n+1}(C)).$ Hence $w\in CU(A).$
Therefore $\overline{ v_1\oplus (1-e)}={\bar z_1}.$

Put $y_1=z_1^*u.$ Then $\|y_1-1_A\|<\eta/4.$
We now repeat the same argument. We obtain
$z_2=s_0\oplus s_1'\oplus\cdots s_1'\in U_0(A),$
where $s_0\in U_0(e_0'Ae_0')$ and where $s_1'$ repeats $n+1$ times,
 $s_1'\in U_0(C_1),$ $C_1\in {\cal I}$ and $e_0'\oplus 1_{C_1}$ is equivalent to
a subprojection of $e$ such that
$$
\|z_2-y_1\|<\eta/16.
$$
Without loss of generality, we may further assume that $ e_0'\oplus 1_{C_1}\le e.$
From the  fact that $\|y_1-1_A\|<\eta/4,$ we may assume
that $\|s_0'-e_0'\|<\eta/2$
and $\|s_1'-1_{C_1}\|<\eta/2.$ Put
$v_2=s_0'\oplus (s_1')^{n+1}\oplus (e-(e_0'\oplus 1_{C_1})).$
Then (since $n<4N+1$)
$$
\|v_2-e\|<\eta/2.
$$
As we have shown, we have $\overline{ v_2\oplus (1-e)}={\bar z_2}.$
Note that
$
\overline{v_1v_2\oplus (1-e)}=\overline{z_1z_2}$ and
$$
 \|z_1z_2-u\|\le \|z_1y_1-u\|+\|z_1y_1-z_1z_2\|=\|z_1y_1-z_1z_2\|<\eta/16.
$$
Also
$$
\|v_1v_2-v_1\|<\eta/2.
$$

Let $y_2=(z_1z_2)^*u.$ Then $\|y_2-1_A\|<\eta/16.$
We can continue the above argument. Consequently,
we obtain a sequence of unitaries $z_n\in U(A)$ and a sequence
of unitaries $v_n\in U(eAe)$ such that
$\overline{v_1v_2\cdots v_n\oplus (1-e)}=\overline{z_1z_2\cdots z_n},$
$$
\|z_1z_2\cdots z_n-u\|\to 0\andeqn \|v_1v_2\cdots v_n-v_1v_2\cdots v_m\|\to 0
$$
as $n, m\to\infty.$
Therefore  we obtain a unitary $v\in U(eAe)$ such that
$$
\hspace{2.4in}\overline{v\oplus (1-e)}={\bar u}.
\,\,\,\,\,\,\,\,\,\,\hspace{2.7in}\QED
$$

\vspace{0.1in}

{\bf 6.8. Lemma}
{\it Let $A$ be a unital \CA\, and ${\cal U}\subset U_0(A)$ be a finite
subset. Then, for any $\ep>0,$ there is a finite subset ${\cal G}\subset A$ and
$\dt>0$ satisfying the following:
for any ${\cal G}$-$\dt$-multiplicative contractive linear map $L: A\to B$
(for any unital \CA\, $B$),
there are unitaries $v\in B$ such that
$$
\|L(u)-v\|<\ep/2\andeqn cel(v)<cel(u)+\ep/2
$$
for all $u\in {\cal U}.$
}

{\it Proof}:
Suppose that
$
z_0(u)=u,z_j(u)\in U_0(A),$ $j=1,2,...,n(u)$ such that ${cel(u)\over{n(u)}}
\le 1/4$ and
$$
cel(z_j(u)(z_{j-1}(u))^*)< {cel(u)\over{n(u)}},\,\,\,j=1,2,...,n(u)
$$
for all $u\in {\cal U}.$
Let $N=\max \{n(u): u\in {\cal U}\}.$
It follows that (for sufficiently large ${\cal G}$ and sufficiently small
$\dt$) there are unitaries $w_j(u)\in U(B)$ such that
$$
\|L(z_j(u))-w_j(u)\|<\ep/8N\pi
$$
for all $j$ and $u\in {\cal U}.$
Thus for all $u\in {\cal U},$
$$
\|L(u)-w_0(u)\|<\ep/2\pi \andeqn cel(w_0(u))<n(u)
[{cel(u)\over{n(u)}}+(\ep/8N)2\pi]
<cel(u)+\ep/2.
$$
\QED

{\bf 6.9. Lemma}
{\it
Let $A$ be a unital simple \CA\, with $TR(A)\le 1$ and let $u\in CU(A).$
Then $u\in U_0(A)$ and for any $\ep>0,$ $cel(u)\le 8\pi+\ep.$
}

{\it Proof}:
We may assume that $u$ is actually in the commutator group.
Write $u=v_1v_2\cdots v_k,$ where each $v_i$ is a commutator.
We write $v_i=a_ib_ia_i^*b_i^*,$ where $a_i$ and $b_i$ are
in $U(A).$
Fix integers $N>0$ and $K>0.$
 Since $TR(A)\le 1,$  by Cor. 3.3, there is a projection
$p\in A$ and a \SCA\, $B\in {\cal I}$ with $1_p=B$
and $B=\oplus_{i=1}^lM_{m_i}(C([0,1]))\oplus_{j=1}^LM_{n_i},$
where $m_i, n_i\ge K$ such that
$$
\|a_i-(a_i'\oplus a_i')\|<\ep/4k,\,\,\,
\|b_i-(b_i'\oplus b_i'')\|<\ep/4k\,\,\,i=1,2,...,k
$$
$$
\|u-(\prod_{i=1}^k a_i'b_i'(a_i')^*(b_i')^*\oplus
a_i''b_i''(a_i'')^*(b_i'')^*)\|<\ep/8,
$$
$a'_i, b_i'\in U((1-p)A(1-p)),$  $a_i'', b_i''\in U_0(B)$
and $N[1-p]\le [p].$
Put $w=\prod_{i=1}^k a_i'b_i'(a_i')^*(b_i')^*$
and $z=\prod_{i=1}^k a_i''b_i''(a_i'')^*(b_i'')^*.$
Then $Det(z)=1.$
It follows from 3.4 in \cite{Ph2} (by choosing $K$ large)
we conclude that $cel(z)\le 6\pi$ in $pAp.$
It is standard to show that $a_i'b_i'(a_i')^*(b_i')^*\oplus (1-p)\oplus (1-p)$
is in $U_0(M_4((1-p)A(1-p)))$ and
it has exponential length no more that $4(2\pi)+\ep/8k.$
This implies that (in $U((1-p)A(1-p))$)
$cel(w\oplus (1-p))\le 8k\pi+\ep/2.$
Note the length only depends on $k.$
We can then choose $N=N(8k\pi+\ep)$ as in 6.4.
In this way, $cel(w\oplus p)\le 2\pi+\ep/2.$
It follows that
$$
cel((w\oplus p)((1-p)\oplus z))\le 8\pi+\ep/2.
$$
The fact that
$$
\|u-(w\oplus p)((1-p)\oplus z)\|<\ep/8
$$
implies that $cel(u)\le 8\pi+\ep.$
\QED

\vspace{0.1in}

{\bf 6.10. Theorem}
{\it Let $A$ be a unital
simple \CA\, with $TR(A)\le 1.$  Let $u, v\in U(A)$ such that
$[u]=[v]$ in $K_1(A)$ and
$$
u^k,v^k\in U_0(A) \andeqn cel((u^k)^*v^k)< L.
$$
Then for any $\ep>0,$
$$
cel(u^*v)\le 8\pi+L/k+\ep.
$$

Moreover, there is $y\in U_0(A)$ with $cel(y)<L/k+\ep$ such that
$\overline{u^*v}={\bar y}$ in $U(A)/CU(A).$
}

{\it Proof:}
Suppose that
$$
u^*v=\prod_jexp(ia_j) \andeqn (u^k)^*v^k=\prod_mexp(ib_m),
$$
where $a_j,b_m\in A_{sa}.$ Since $cel((u^k)^*v^k)< L,$ we may assume that
$\sum\|b_m\|<L$ (see \cite{Rn}).
Let $M=\sum_j\|a_j\|.$ (So $cel(u^*v)\le M$.)
Since $TR(A)\le 1$ , for any $\dt>0$ with
$\dt/(1-\dt)<\ep/2(M+L+1)$ and sufficiently small $\eta>0$ and
with a sufficiently large finite subset
${\cal G}$ (which contains
$a_j, b_m$),
 there exist a projection $p\in A$ and a unital
$C^*$-subalgebra $F\subset A$ with $1_F=p$ and $F\in {\cal I}$ such that

(1) $pxp\in_{\eta} F$ for all $x\in {\cal G},$

(2) $\|u-u_0\oplus u_1\|<\eta,$ $\|v-v_0\oplus v_1\|<\eta,$
and

(3) $u_0$ and $v_0$ are unitaries in $(1-p)A(1-p)$ and
$u_1$ and $v_1$ are unitaries in $pAp,$

(4) $cel(u_0^*v_0)\le M+1$ in $(1-p)A(1-p)$ and
$cel((u_1^k)^*v_1^k)<L$ in $F,$

(5) $\tau(1-p)<\dt$ for all $\tau\in T(A).$

(Note (4) follows from 6.8).

Write $F=\oplus_s^N F_s,$ where each $F_s=M_{n(s)}(C([0,1]))$ or
$F_s=M_{n(s)}.$ By corollary 3.3, we may assume that each
$n(s)>max(2\pi^2/\ep, K(1)),$ where $K(1)$ is the number described
in Lemma 3.4 in \cite{Ph2} (with $d=1$).

First consider the case in which $N=1$ and $F=M_K(C([0,1]))$
(so $K>max(2\pi^2/\ep, K(1))$).
Note that $cel((u_1^*)^kv_1^k)<L$ in $F.$
Therefore, by Lemma 3.3 (1) of \cite{Ph2},
there exists $a\in F_{sa}$ with $\|a\|<L$ such that
$$
det(exp(ia)(u_1^*)^kv_1^k)=1 \,\,\,\,\,\,({\rm for\,\,\, every\,\,\,}
t\in [0,1]).
$$
This implies that
$$
det((exp(ia/k)u_1^*v_1)^k)=1\,\,\,\,\,\,({\rm for\,\,\, every\,\,\,}
t\in [0,1]).
$$
Therefore
$$
det(exp(ia/k)u_1^*v_1)=exp(i2l\pi/k) \,\,\,\,\,\,({\rm for\,\,\,
every }\,\,\,t\in [0,1]
$$
for some $l=0,1,..,k-1.$
Note since determinant is a continuous function on $[0,1],$
the above function has to be constant (only one value of $l$
occurs). Set $f(t)=-2l\pi/k.$ Then $f\in C([0,1])_{sa}$ and
$\|f\|\le 2\pi.$ Note that $exp(if/K)\cdot 1_F$ commutes with
$exp(ia/k)$ and
$(exp(if/K)\cdot 1_F)exp(ia/k)=exp(i(f/K +a/k)).$
We have
$$
det((exp(if/K)\cdot 1_F) exp(ia/k)u_1^*v_1)=1.
$$
So, by 3.4 (and 3.1) in \cite{Ph2},
$$
cel(u_1^*v_1)\le 2\pi/K+L/k +6\pi\,\,\,{\rm (in\,\,\, }F{\rm )}.
$$
Moreover, by 2.4 in \cite{Tm1},
$z_1=(exp(if/K)\cdot 1_F) exp(ia/k)u_1^*v_1\in CU(F).$
Note the above also holds when $F=M_K.$
By considering each summand, the above also holds
for the case in which $N>1.$

Moreover, by Lemma 6.4, there is $y'\in CU(A)$ and $y''\in U_0(A)$ such that
$(u_0\oplus p)^*(v_0\oplus p)=y'y''$ and $cel(y'')<\ep/2.$
Note that $((1-p)\oplus z_1)y'\in CU(A).$
Therefore
$$
\overline{u^*v}=\overline{exp(if/K)\cdot 1\cdot  exp(ia/k) w}
$$ for
some $w\in U_0(A)$ with $cel(w)<\ep/2$ if $\eta$ is sufficiently
small. Therefore, by 6.9, with sufficiently small $\eta,$
$$
cel(u^*v)\le 2\pi/K+L/k +8\pi+\ep/2<8\pi+L/k+\ep.
$$
\QED

{\bf 6.11. Theorem}
{\it Let $A$ be a unital separable simple \CA\, with $TR(A)\le 1$
and $u\in U_0(A).$ Suppose that $u^k\in CU(A)$ for some
integer $k>0,$ then $u\in CU(A).$ In particular, $U_0(A)/CU(A)$
is torsion free.}

{\it Proof:}
The proof is essentially the same as that of
6.10.
Let $\ep>0$ and let
$v=\prod_{j=1}^ra_ib_i(a_i')^{-1}(b_i')^{-1}$ be such that
$$
\|u^k-v\|<\ep/64.
$$
Put $l=cel(u^k).$
Let $\dt>0$ be such that $2(l+\ep)\dt<\ep/64\pi.$
Fix  a finite subset ${\cal G}\subset A$
which  contains $u,u^k, v, a_i, b_i, a_i^{-1}, b_i^{-1}$
among other elements.

Since $TR(A)\le 1,$ there is a projection $p\in A$ and a unital
\SCA\, $F\in {\cal I}$ with $1_F=p$ such that

(1) $pxp\in_{\ep/64} F$ for all $x\in {\cal G},$

(2) $\|v-v_0\oplus v_1\|<\ep/32,$ $\|u-u_0\oplus u_1\|<\ep/32,$
and $\|u^k-u_0^k\oplus u_1^k\|<\ep/32,$

(3) $cel(u_0^k)<l+\ep/32$ in $U((1-p)A(1-p))$ and

(4) $\tau(1-p)<\dt$ for all $\tau\in T(A).$

(note (3) follows from 6.8 with large ${\cal G}$).

Here $u_0, v_0\in U((1-p)A(1-p))$ and $v_1, u_1\in U(F).$
Moreover, we may assume that there are $a_i', b_i'\in U(F)$ such that
$$
\|u_1^k-\prod_{j=1}^r a_i'b_i'(a_i')^{-1}(b_i')^{-1}\|<\ep/32.
$$
Put $w=\prod_{j=1}^r a_i'b_i'(a_i')^{-1}(b_i')^{-1}.$
Since  $U(F)=U_0(F),$ we may write  $w=\prod_{m=1}^Lexp(i d_m)$ for
some $d_m\in F_{s.a}.$ Put $w_k=\prod_{m=1}^Lexp(i d_m/k).$
Then $w_k^k=w.$
So
$$
cel((u_1)^k(w_k^*)^k)<{\ep\pi\over{32}}
$$
Write $F=\oplus_{s=1}^NF_s,$ where each $F_s=M_{r(s)}C([0,1])$ or
$F_s=M_{r(s)}.$ By 3.3, we may assume that each $n(s)>\max (16\pi^2/\ep,K(1)),$
where $K(1)$ is the number described in Lemma 3.4 in [Ph2] (with $d=1$).

As in the proof of 6.10,
$$
{\rm det}(exp(if/K)exp(ia/k)u_1w_k^*)=1
$$
for some $f\in F_{s.a}$
with $\|f\|\le 2\pi$ and
$a\in F_{s.a}$ with $\|a\|<\ep\pi/32$ (with $K>\max(16\pi^2/\ep, K(1))$).
By 2.4 in \cite{Tm1}, $exp(if/K)exp(ia/k)u_1w_k^*\in CU(F).$
We also have
$$
\|exp(if/K)exp(ia/k)-1_F\|<\ep/8+\ep/32.
$$
Thus
$$
{\rm dist}(\overline{u_1w_k^*},{\bar 1})<\ep/8+\ep/32.
$$
Since $det(w)=1,$ as in the proof of 6.10, we also have
$$
{\rm det}(exp(ig/K)w_k^*)=1
$$
for some $g\in F_{s.a}$ with $\|g\|\le 2\pi.$
Again, $exp(ig/K)w_k^*\in CU(F).$
But
$$
\|(exp(ig/K)w_k)w_k^*-1\|\le \|exp(ig/K)-1\|<\ep/4.
$$
So
$$
{\rm dist}(\overline{w_k},{\bar 1})<\ep/4.
$$
Therefore
$$
{\rm dist}({\bar u_1},{\bar 1})\le {\rm dist}({\bar u},\overline{ w_k})+
{\rm dist}(\overline{w_k},{\bar 1})<\ep/8+\ep/32+\ep/4<\ep/2
$$
in $U(F)/CU(F).$
On the other hand, by 6.4 and the choice of $\dt,$
$$
cel((u_0\oplus p)z)<\ep/(8\pi)
$$
for some $z\in CU(A).$
Thus
$$
\inf\{\|u_0\oplus u_1-x\|: x\in CU(A)\}<\ep/8+\ep/8+\ep/32+\ep/4<3\ep/4.
$$
This implies that
$$
\inf\{\|u-x\|: x\in CU(A)\}<\ep.
$$
Therefore $u\in CU(A).$ Consequently $U_0(A)/CU(A)$ is torsion free.
\QED

\vspace{0.1in}

{\bf 6.12. Corollary} {\it Let $B_n$ be a sequence of unital
simple \CA\, with $TR(B_n)\le 1.$
Let $\prod^b_nK_1(B_n)$ be the set of sequences
$z=\{z_n\},$ where $z_n\in K_1(B_n)$ and
$z_n$ can be represented by unitaries in $M_{K(z)}(B_n)$ for
some integer $K(z)>0$.}
{\it Then the kernel of the map
$$
K_1(\prod_n B_n)\to \prod^b_nK_1(B_n)\to 0
$$
is a divisible and torsion free subgroup of $K_1(\prod_n B_n).$}

{\it Proof:} By 6.5, the exponential rank of each $B_n$
is bounded by 4. Therefore that the kernel is divisible follows
 from the fact that each $B_n$ has stable
rank one (and has exponential rank bounded by $4$) (see
\cite{GL2}).
Suppose that $\{u_n\}\in U(M_K(\prod_n B_n))$ such that
$[\{u_n\}]$ is in the kernel
and $k[\{u_n\}]=0.$
By changing notation (with different $\{u_n\}$ and
larger $K$), we may assume that
$\{u_n^k\}\in U_0(M_K(\prod_n B_n)).$
Also each $u_n\in U_0(B_n).$
This implies that there is $L>0$ such that
$cel(u_n^k)\le L$
for all $n.$
It follows from 6.10 that
$$
cel(u_n)\le 8\pi+L+L/k +\pi/4\,\,\,\,\,\,\,{\rm for\,\,\, all}\,\,\, n.
$$
This implies (see for example \cite{GL2})
$\{u_n\}\in U_0(M_L(\prod_n B_n)).$
Therefore $[\{u_n\}]=0$ in $K_1(\prod_n B_n).$
So the kernel is torsion free.
\QED

\vspace{0.1in}

{\bf 6.13. Remark}
When ${\rm dim} X\le 1,$ we believe that the conclusion
of 3.4 in \cite{Ph2}  can be improved and $cel(u)$ should be
no more than $D(u)+2\pi.$ This could be achieved by
a modification of Phillips's argument as we were informed by
N. C. Phillips.  Consequently, in 6.10,  $8\pi$ could be
replaced by $4\pi$ and in 6.9, $8\pi$ could be replaced by
$4\pi.$ However, we will not use these better estimates.

\vspace{0.1in}

\section{ Homomorphisms from $U(C)/CU(C)$ to $U(B)/CU(B)$}

{\bf 7.1. Definition}
Let $Y$ be a connected finite CW complex with dimension
no more than three with torsion $K_1(C(Y))$ and $C'=PM_r(C(X))P,$ where
$X=S^1\vee S^1\vee\cdots S^1\vee Y$
and $P\in M_r(C(X))$ is a projection and $P$ has rank $R\ge 6.$
We assume that $S^1$ is repeated $s$ ($\ge 0$) times.
Note that the above includes the case that $X=Y=[0,1].$
Then $K_1(C')=tor(K_1(C'))\oplus G_1,$ where $G_1$ is $s$ copies of
${\bf Z}.$
Denote by $\xi$ the point in $X$ where each $S^1$ and $Y$ meet.
Rename each $S^1$ by $\Omega_i,$ $i=1,2,...,s.$
Denote by $z_i'$ the identity map on $i$-th  $S^1$ ($\Omega_i$).
Define $z_i''(\zeta)=\zeta$ on $\Omega_i$ (which is identified with
$S^1$) and $z_i''(\zeta)=\xi$ for all $\zeta\not\in \Omega_i.$
There is an obvious  \hm\,
$\Pi: PM_r(C(X))P\to D'=\oplus_{i=1}^k E_i,$ where $E_i\cong M_R(C(S^1)).$
 Note that if $s\ge 2,$ then $\Pi$ is not surjective.
We have that $G_1=K_1(D').$
We also use $\Pi_i: PM_r(C(X))P\to E_i$ which is the composition
of $\Pi$ with the projection from $D'$ to $E_i.$
Let $z$ be the identity map on $S^1.$

We may write
$$
P(t)=\begin{pmatrix}
                  P_1 & 0\\
                     0 & I\end{pmatrix}
$$
where $P_1$ is a projection with
rank 3 and $I={\rm diag}(1,1,...,1)$ with  $1$ repeating
${\rm rank}(P)-3$ times,

Note that, since ${\rm rank} (P)\ge 6,$  ${\rm tsr}(C(X))=2$ and
 ${\rm csr}(C(X))\le 2+1$
(by 4.10 of \cite{Rf1}). It follows that ${\rm csr}M_3(C')\le 2$
(by 4.7 in \cite{Rf2}).  It then follows from 5.3 of \cite{Rf2}
that $U(C')/U_0(C')=K_1(C').$ In particular, $CU(C')\subset U_0(C').$

Denote by $u_i={\rm diag}(z_i'',1,...,1),$ where $1$ is repeated
$r-4$ times.
If we write $z_i\in U(C),$ we mean the unitary
$$
z_i(t)=\begin{pmatrix} P_1 & 0\\
                       0 & u_i
                       .\end{pmatrix}
$$
If we write $z_i\in E_i,$ we mean $\Pi_i(z_i'').$
Note that in this case, $z_i$ has the form :
${\rm diag}(1,...,1,z,1,...,1),$ where $z$ is in the 4-th position and
there are $R-1$ many 1's.

Now let $C=\oplus_{j=1}^{l+l_1} C^{(j)},$ where $C^{(j)}$ is either of the form
$P_jM_{r(j)}(C(X_j))P_j$ for $j\le l,$ where $X_j$ is of the form $X$
described
above, $C^{(j)}=M_{r(j)}$, or $C^{(j)}=P_jM_{r(j)}(C(Y_j))P_j,$
where $Y_j$ is a finite CW complex with dimension no more than
$3,$ rank of $P_j$ is $R(j)\ge 6$ and $K_1(Y_j)$ is finite for
$l+1\le j\le l+l_1.$
  Let $D^{(j)}$ be as $D'$ above for each $j\le l.$
Let $\Pi^{(j)}$ be as $\Pi$ above for $C^{(j)}=P_jM_{r(j)}(C(X_j))P_j.$
Put $D=\oplus_{j=1}^lD^{(j)}$ and
 $\Pi=\oplus_{j=1}^{l} \Pi^{(j)}.$
Since $K_1(C)$ is finitely generated
and $U_0(C)/CU(C)$ is divisible (see 6.6), we may write
$$
U(C)/CU(C)=U_0(C)/CU(C)\oplus K_1(D)\oplus tor(K_1(C)).
$$
Let $\pi_1: U(C)/CU(C)\to K_1(D),$ $\pi_0: U(C)/CU(C)\to
U_0(C)/CU(C)$ and $\pi_{2}: U(C)/CU(C)\to tor(K_1(C))$
be fixed projection maps associated with the above
decomposition. To avoid possible confusion, by $\pi_i(U(C)/CU(C)),$
we mean a subgroup of $U(C)/CU(C),$ $i=0,1,2.$
We also assume that $\pi_1({\bar z_i})={\bar z_i}$ (in $U(C)/CU(C)$).

It is worth to point out that one could have $X=Y=[0,1].$

We will keep these notation in the rest of this section.

\vspace{0.1in}

{\bf 7.2. Lemma} {\it Let $C=\oplus_{i=1}^{l+l_1}
C_i$ be as above and ${\cal U}\subset U(C)$ be a finite subset and
$F$ be the group generated by ${\cal U}.$
Suppose that $G$ is a subgroup of
$U(C )/CU(C )$ which contains ${\bar F},$ $\pi_2(U(C)/CU(C )$
and $\pi_1(U( C)/CU (C )).$
Suppose that the composition map $\gamma: {\bar F}\to
U(D)/CU(D)\to U(D)/U_0(D)$ is injective  and $\gamma({\bar F})$ is free.
Let $B$ be a unital \CA\, and $\Lambda: G\to U(B)/CU(B)$
be a \hm\, such that $\Lambda(G\cap U_0(C)/CU(C))\subset U_0(B)/CU(B).$
Then there is a \hm \, $\beta: U(D)/CU(D)\to U(B)/CU(B)$
with $\beta(U_0(D)/CU(D))\subset U_0(B)/CU(B)$ and a \hm\,
$\theta: \pi_2(U(C)/CU(C))\to U(B)/CU(B)$ such that
$$\beta\circ
\Pi^{\ddag}\circ \pi_1({\bar w})= \Lambda({\bar w})(\theta\circ \pi_2({\bar
w}))
$$
for all $w\in F$ and such that $\theta(g)=\Lambda|_{\pi_2(U(C
)/CU( C))}(g^{-1})$ for $g\in \pi_2(U(C)/CU(C) ).$ Moreover,\\
$\beta(U_0(D)/CU(D))\subset U_0(B)/CU(B).$

If furthermore $A$ is a simple \CA\, with
$TR(B)\le 1$ and $\Lambda(U(C)/CU(C))\subset U_0(B)/CU(B),$
then $\beta\circ \Pi^{\ddag}\circ (\pi_1)|_{\bar F}=\Lambda|_{\bar F}.$

}

{\it Proof:} Let $\kappa_1: U(D)/CU(D)\to K_1(D)$ be the quotient
map. Let $\eta: \pi_1(U(C)/CU(C))\to K_1(D)$ be defined by
$\eta=\kappa_1\circ \Pi^{\ddag}|_{\pi_1(U(C)/CU(C))}.$
Note that $\eta$ is an isomorphism.
Since
$\gamma$ is injective and $\gamma({\overline{F}})$ is free, we conclude
that
$\kappa_1\circ \Pi^{\ddag}\circ \pi_1$ is  also injective on ${\bar F}.$
From this fact and the fact that
$U_0(C)/CU(C)$ is divisible (6.6), we obtain  a \hm\, $\lambda: K_1(D)\to
U_0(C)/CU(C)$ such that
$$
\lambda|_{\kappa_1\circ \Pi^{\ddag}\pi_1({\bar F})}=\pi_0\circ
((\kappa_1\circ \Pi^{\ddag}\circ \pi_1)|_{\bar F})^{-1}.
$$
Now define $\beta=\Lambda((\eta^{-1}\circ \kappa_1)
\oplus(\lambda\circ \kappa_1)).$
Then for any ${\bar w}\in {\bar F},$
$$
\beta(\Pi^{\ddag}\circ\pi_1({\bar w}))=
\Lambda[\eta^{-1}(\kappa_1\circ \Pi^{\ddag}(\pi_1({\bar w})))
\oplus  \lambda\circ \kappa_1(\Pi^{\ddag}(\pi_1({\bar w})))]
=
\Lambda(\pi_1({\bar w})\oplus
\pi_0({\bar w})).
$$

Now define $\theta: \pi_2(U(C)/CU(C))\to U(B)/CU(B)$
by $\theta(x)=\Lambda(x^{-1})$ for $x\in \pi_2(U(C)/CU(C)).$
Then
$$
\beta(\Pi^{\ddag}(\pi_1({\bar w})))=
\Lambda({\bar w})\theta(\pi_2({\bar w}))\,\,\,{\rm for}\,\,\, w\in F.
$$
To see the last statement, we assume that
$\Lambda(U(C)/CU(C))\subset U_0(B)/CU(B).$
Then $\Lambda(\pi_2(U(C)/CU(C)))$ is a torsion subgroup
of $U_0(B)/CU(B).$ By 6.11, $U_0(B)/CU(B)$ is torsion free.
Therefore $\theta=0.$
\QED

\vspace{0.1in}

{\bf 7.3. Lemma} {\it Let $A$ be  a unital separable simple \CA\,
with $TR(A)\le 1$ and $C$ be as described in 7.1. Let  ${\cal
U}\subset U(A)$ be a finite subset and $F$ be the subgroup
generated by ${\cal U}$ such that $(\kappa_1)|_{\bar F}$ is
injective and $\kappa_1({\bar F})$ is free, where  $\kappa_1:
U(A)/CU(A)\to K_1(A)$ is the quotient map. Suppose that $\alpha:
K_1 (C )\to K_1(A)$ is an injective \hm\, and $L: {\bar F}\to
U(C)/CU(C )$ is an injective \hm\, with $L({\bar F}\cap U_0(A
)/CU( A))\subset U_0(C)/CU(C)$ such that $\pi_1\circ L$ is
injective (see 7.1 for $\pi_1$) and
$$
\alpha\circ \kappa_1'\circ L(g)=\kappa_1(g)\,\,\,{\rm for\,\,\,
all}\,\,\, g\in {\bar F},
$$
where $\kappa_1': U( C) /CU( C)\to K_1(C )$ is the quotient map.
 Then
there exists a \hm\, $\beta: U( C)/CU( C)\to U(A)/CU(A)$
with $\beta(U_0(C)/CU(C ))\subset U_0(A)/CU(A)$ such that
$$
\beta\circ L({\bar w})={\bar w}\,\,\,\,{\rm for}\,\,\, w\in F.
$$
}

{\it Proof}:
Let $G$ be the preimage
 of $\alpha\circ \kappa_1'(U(C)/CU(C))$ under $\kappa_1.$
So we have the following short exact sequence:
$$
0\to U_0(A)/CU(A)\to G\to \alpha\circ \kappa_1'(U(C)/CU(C))\to 0.
$$
Since $U_0(A)/CU(A)$ is divisible, there exists
an injective \hm\, $\gamma: \alpha\circ \kappa_1'(U(C)/CU(C))
\to G$ such that $\kappa_1\circ \gamma(g)=g$ for
$g\in \alpha\circ \kappa_1'(U(C)/CU(C)).$
Since $\alpha\circ \kappa_1'\circ L(f)=\kappa_1(f)$ for all
$f\in {\bar F},$ we have ${\bar F}\subset G.$
Moreover,
$$
(\gamma\circ \alpha\circ \kappa_1'\circ L(f))^{-1}f\in U_0(A)/CU(A)
$$
for all $f\in {\bar F}.$
Define $\psi: L({\bar F})\to U_0(A)/CU(A)$ by
$$
\psi(x)=\gamma\circ \alpha\circ \kappa_1'\circ L([(L)^{-1}(x)]^{-1})L^{-1}(x)
$$
for $x\in L({\bar F}).$
Since $U_0(A)/CU(A)$ is divisible,
there is \hm\, ${\tilde \psi}: U(C)/CU(C)\to U_0(A)/CU(A)$
such that ${\tilde \psi}|_{L({\bar F})}=\psi.$
Now define
$$
\beta(x)=\gamma\circ \alpha\circ \kappa_1'(x){\tilde \psi}(x).
$$
Hence $\beta(L(f))=f$ for $f\in {\bar F}.$
\QED

\vspace{0.1in}

{\bf 7.4. Lemma}
{\it Let $B$ be a unital separable simple \CA\, with $TR(B)\le 1$
and $C$ be as in 7.1. Let $F$ be a group generated by a finite
subset ${\cal U}\subset U( C)$ such that $(\pi_1)|_{\bar F}$ is
injective. Let $G$ be a subgroup
containing ${\bar F},$ $\pi_0({\bar F}),$
$\pi_1(U( C)/CU( C))$ and $\pi_2(U(C )/CU( C)).$
Suppose that $\alpha: U(C)/CU(C)\to U(B)/CU(B)$ is a \hm\,
with $\alpha(U_0(C)/CU(C))\subset U_0(B)/CU(B))$. Then
for any
$\ep>0$ there is $\dt>0$ satisfying the following:
if
$\phi=\phi_0\oplus\phi_1: C\to B$ is a ${\cal G}$-$\eta$
-multiplicative \morp\, such that

  {\rm  (1)} both $\phi_0$ and $\phi_1$ are ${\cal G}$-$\eta$
-multiplicative and $\phi_0$ maps the identity of each summand
of $C$ into a projection,

{\rm (2)} ${\cal G}$ is sufficiently large and $\eta$ is sufficiently
small which depend only on $C$ and $F$ (so that $\phi^{\ddag}$ is
well defined on a $G$)

{\rm (3)} $\phi_0$ is homotopically trivial (see (vi) in section
1), $(\phi_0)_{*0}$ is a well-defined \hm\, and
$[\phi]|_{\kappa_1({\bar F})}=\alpha_*|_{\kappa_1({\bar F})},$
where $\alpha_*: K_1(C)\to K_1(B)$ induced by $\alpha$ and
$\kappa_1: U(C)/CU(C)\to K_1(C)$ is the quotient map,

{\rm (4)} $\tau(\phi_0(1_C))<\dt$ for all $\tau\in T(B),$

 then
there is a \hm\, $\Phi: C\to e_0Be_0$ {\rm (}$e_0=\phi_0(1_C)$ {\rm )}
such that

{\rm (i)} $\Phi$ is homotopically trivial and $\Phi_{*0}=(\phi_0)_{*0}$ and

{\rm (ii)}
$$
\alpha({\bar w})^{-1}(\Phi\oplus \phi_1)^{\ddag}({\bar w})=\overline{g_w},
$$
where $g_w\in U_0(B)$ and $cel(g_w)<\ep$ for
all $w\in {\cal U}.$
}

{\it Proof}:
It follows from 7.2 that there are
\hm s $\beta_1,\, \beta_2: U(D)/CU(D)\to U(B)/CU(B)$
with $\beta_i(U_0(D)/CU(D))\subset U_0(B)/CU(B)$ ($i=1,2$)
and \hm s $\theta_1,\theta_2: \pi_2(U(C)/CU(C))\to U(B)/CU(B)$
such that
$$
\beta_1\circ \Pi^{\ddag}(\pi_1({\bar w}))=
\alpha({\bar w})\theta_1(\pi_2({\bar w}))\andeqn
\beta_2\circ \Pi^{\ddag}(\pi_1({\bar w}))=\phi_1^{\ddag}({\bar w}^*)
\theta_2(\pi_2({\bar w}))
$$
for all ${\bar w}\in {\bar F}.$
Moreover $\theta_1(g)=\alpha(g^{-1})$ and $\theta_2(g)=
\phi_1^{\ddag}(g)$ if $g\in \pi_2({\bar F}).$
Since $\phi_0$ is homotopically trivial,
$$
\theta_1(g)\theta_2(g)\in U_0(B)/CU(B)\,\,\,{\rm for\,\,\, all}\,\,\,
g\in \pi_2({\bar F}).
$$
Since $\pi_2(U( C)/CU(C )$ is torsion and $U_0(B)/CU(B)$ is torsion
free, we conclude that
$$
\theta_1(g)\theta_2(g)={\bar 1}\,\,\,{\rm for\,\,\,all}\,\,\,g\in
\pi_2({\bar F}).
$$
To simplify notation,
without loss of generality, we may
assume that $C=\oplus_{j=1}^{1+l_1}C^{(j)}$ (with $l=1$) such that
$C^{(1)}=PM_r(C(X))P$ as described in 7.1 and $C^{(j)}$ is also
as described in 7.1 for $2=l+1\le j\le l_1+1.$
Let $D$ be as described in 7.1.

For each $w\in U(C),$ we may write $w=(w_1,w_2,...,w_{1+l_1})$
according to the direct sum $C=\oplus_{j=1}^{1+l_1}C^{(j)}.$ Note
that $\pi_1(w)=\pi_1(w_1).$ Let $\pi_1({\bar w})=({\bar
z_1}^{k(1,w)},{\bar z}_2^{k(2,w)}, ...,{\bar z}_s^{k(s,w)}),$
where $k(i,w)$ is an integer (here $z_i$ as described in 7.1).
 Then
$\Pi_i^{\ddag}(\pi_1({\bar w}_1))={\bar z_i}^{k(i,w)}.$
On the other hand, we may also write
$\Pi_i^{\ddag}(\bar w_1)=\overline{z_i^{k(i,w)}g_{i
,w}}$ for some
$g_{i,w}\in U_0(C(S^1, M_R)).$

Let $l=\max\{cel(g_{i,w}): w\in {\cal U}, 1\le i\le s\}.$ Choose
$\dt$ so that $(2+l)\dt<\ep/4\pi.$ Let $e_0=\phi_0(1_C)$ and
$e_1=\phi_1(1_C).$ Write $e_0=E_1\oplus E_2\oplus\cdots\oplus
E_{1+l_1},$ where $E_j=\phi_0(1_{C^{(j)}}),$ $j=1,2,...,1+l_1.$
Recall that $P$ has rank $R.$ Since $\phi_0$ is homtopically
trivial (see (vi) in Section 1), we may also write $E_1=
e_{01}\oplus\cdots \oplus e_{0R},$ where $\{e_{0i}: 1\le i\le R\}$
is a set of  mutually orthogonal and mutually equivalent
projections. Since $e_0Be_0$ is simple and has the property (SP),
$e_{01}$ can be written as a sum of $s$ mutually  orthogonal
projections. Thus $E_1=p_1\oplus p_2\oplus \cdots \oplus p_s,$
where each $p_i$ can be written as a direct sum of $R$ mutually
orthogonal and mutually equivalent projections
$\{q_{i,1},...,q_{i,R}\}.$ For each $q_{i1},$ we write
$q_{i1}=q_{i1,1}\oplus q_{i1,2},$ where both $q_{i1,1}$ and
$q_{i1,2}$ are not zero. Let $q=\sum_{i=1}^s q_{i1}.$ We may view
$E_1BE_1=\oplus_{i=1}^sM_R(q_{i1}Bq_{i1})=M_R(qBq).$

Let $z_i$ be as in 7.1. Put $x_i'\in U(q_{i1,1}Bq_{i1,1})$ such
that $ \overline{x_i'}=\beta_1({\bar z_i}),$ and $y_i'\in
U(q_{i1,2}Bq_{i1,2})$ such that $\overline{y_i'}=
\beta_2(\overline{z_i}),$ $i=1,2,...,s.$ This is possible because
of 6.7. Put $x_i=x_i'\oplus q_{i1,2},$ $y_i=y_i'\oplus q_{i1,1},$
 $i=1,2,...,s.$
Note that $x_iy_i=y_ix_i.$
Define $\Phi_1: D\to M_R(qBq)=\oplus_{i=1}^sM_R(q_{i1}Bq_{i1})$ by
$\Phi_1(f)=\sum_{i=1}^s f_i(x_iy_i),$ where $f=(f_1,f_2,...,f_s)$ and
$f_i\in C(S^1, M_R).$
Define $h(g)=\Phi_1(\Pi(g))\oplus\phi_1(g)$ for $g\in C.$
We compute that
$$
\overline{h(w)}=
\prod_{i=1}^s\overline{x_i^{k(i,w)}g_{i,w}(x_iy_i)y_i^{k(i,w)}}
 \phi_1^{\ddag}({\bar w})=
\beta_1(\Pi^{\ddag}(\pi_1({\bar w})))
\Phi_1^{\ddag}(\overline{\oplus_{i=1}^sg_{i,w}})
\beta_2(\Pi^{\ddag}(\pi_1({\bar
w})))\phi_1^{\ddag}(\overline{w})
$$
$$
=\alpha({\bar w})
\theta_1(\pi_2({\bar w}))\theta_2(\pi_2({\bar
w}))\Phi_1^{\ddag}(\overline{\oplus_{i=1}^sg_{i,w}})
=\alpha({\bar w})\Phi_1^{\ddag}(\overline{\oplus_{i=1}^sg_{i,w}})
$$
for all $w\in {\cal U}.$
Put $g_w'=\Phi_1(\oplus_{i=1}^sg_{i,w})\oplus(1-\phi_0(1_C)).$
Since $\tau(\phi_0(1_C))<\dt,$ by the choice
of $\dt,$ we conclude from Lemma 6.4 that there exists $w'\in CU(B)$ such that
$$
cel(w'g_w')<\ep/2\,\,\,\,{\rm for\,\,\,all}\,\,\, w\in {\cal U}.
$$
Note that $\Phi_1\circ \Pi$ factors through $D$ and
$(\Phi_1)_{*1}=0.$ In particular, $\Phi_1\circ \Pi$ is
homotopically trivial. Since $\phi_0$ is homotopically trivial, it
is easy to see that there  is a point-evaluation map $\Phi_2:
\oplus_{j=2}^{1+l_1}C^{(j)} \to
(\oplus_{j=2}^{1+l_1}E_j)B(\oplus_{j=2}^{1+l_1}E_j).$ Now define
$\Phi=\Phi_1\circ \Pi\oplus \Phi_2.$ We see that we can make (by a
right choice of $\Phi_2$) $\Phi|_{*0}=(\phi_0)_{*0}.$ It is clear
that $\Phi$ is homotopically trivial. Let $g_w''=\Phi_2(w)\oplus
(1-(e_0-E_1)).$ Since $\Phi_2(\sum_{j=2}^{1+l_1}C^{(j)})$ is
finite dimensional, $cel(\Phi_2(w))\le 2\pi$ (in
$U_0((e_0-E_1)B(e_0-E_1))$ for all $w\in {\cal U}).$ By the choice
of $\dt,$ we conclude that there is $w''\in CU(B)$ such that
$cel(w''g_w'')<\ep/2$ (see 6.4).  Put $g_w=w'g_w'w''g_w''.$ We
have, for all $w\in {\cal U},$
$$
\alpha({\bar w})^{-1}(\Phi\oplus \phi_1)^{\ddag}({\bar w})=
\overline{g_w}\,\,\,{\rm with}\,\,\, g_w\in U_0(B)
\andeqn cel(g_w)<\ep.
$$
\QED

{\bf 7.5. Lemma} {\it Let $B$ be a unital separable simple
\CA\, with $TR(B)\le 1$ and $C$ be as described in 7.1.
Let  ${\cal U}\subset U(B)$ be a finite subset and $F$ be
the subgroup generated by ${\cal U}$ such that $\kappa_1({\bar F})$ is
free, where  $\kappa_1: U(B)/CU(B)\to K_1(B)$ is the quotient map.
Let $\phi: C\to B$ be a \hm\, such that $\phi_{*1}$ is injective.
Suppose that
$j,L: {\bar F}\to U(C)/CU(C)$ are two injective \hm s
with $j({\bar F}\cap U_0(B)/CU(B)), L({\bar F}\cap U_0(B)/CU(B))
\subset U_0( C)/CU(C )$ such that
$\kappa_1\circ \phi^{\ddag}\circ  L=\kappa_1\circ \phi^{\ddag}\circ j= \kappa_1|_{\bar F}$
and all three are injective.

Then, for any $\ep>0,$ there exists $\dt>0$ such that
if $\phi=\phi_0\oplus \phi_1: C\to B,$ where $\phi_0$ and
$\phi_1$ are \hm s
satisfying the following:

{\rm (1)} $\tau(\phi_0(1_C))<\dt$ for all $\tau\in T(B)$ and

{\rm (2)} $\phi_0$ is homotopically trivial,

then there is a \hm\, $\psi: C\to e_0Be_0$ {\rm (}$e_0=\phi_0(1_C))$ {\rm )}
such that

{\rm (i)} $[\psi]=[\phi_0]$ in $KL(C, B)$ and

{\rm (ii)} $$ (\phi^{\ddag}\circ j({\bar w}))^{-1}
 (\psi\oplus \phi_1)^{\ddag}(L({\bar w}))=g_w,$$
where $g_w\in U_0(B)$ and $cel(g_w)<\ep$
for all $w\in {\cal U}.$
}

{\it Proof}:
The first part of the proof is essentially the same as
that of 7.3.
Let $\kappa_1': U(C)/CU(C)\to K_1(C)$ be the quotient map and let
$G$ be the preimage of
 $\phi_{*1}\circ \kappa_1'(U(C)/CU(C))$ under $\kappa_1.$
Since $U_0(B)/CU(B)$ is divisible, there exists
an injective \hm\, $\gamma: \phi_{*1}\circ \kappa_1'(U(C)/CU(C))
\to G$ such that $\kappa_1\circ \gamma(g)=g$ for
$g\in \phi_{*1}\circ \kappa_1'(U(C)/CU(C)).$
Since
$
\phi_{*1}\circ \kappa_1'\circ L(f)=\kappa_1(f)=\kappa_1(\phi^{\ddag}\circ
j(f))
$ for all
$f\in {\bar F},$ we have ${\bar F}\subset G.$
Moreover,
$$
[\gamma\circ \phi_{*1}\circ \kappa_1'\circ L(f)]^{-1} \phi^{\ddag}\circ j(f)
\in U_0(B)/CU(B)
$$
for all $f\in {\cal F}.$
Define $\psi:L({\bar F})\to U_0(B)/CU(B)$ by
$$
\psi(x)=[\gamma\circ \phi_{*1}\circ \kappa_1' (x)]^{-1}
[\phi^{\ddag}\circ j\circ L^{-1}(x)]
$$
for all $x\in L({\bar F}).$
Since $U_0(B)/CU(B)$ is divisible, there is a \hm\,
${\tilde \psi}: U(C)/CU(C)\to U_0(B)/CU(B)$ such that
${\tilde \psi}|_{L(\bar F)}=\psi.$
Define
$\alpha: U(C)/CU(C)\to U(B)/CU(B)$ by
$$
\alpha(x)=\gamma\circ \phi_{*1}\circ \kappa_1'(x){\tilde \psi}(x)\,\,\,
{\rm for\,\,\,all}\,\,\, x\in U(C)/CU(B).
$$
Note that
$$
\alpha(L(f))=\phi^{\ddag}\circ j(f)\,\,\,{\rm for\,\,\, all}\,\,\,
f\in {\bar F}.
$$

Now the lemma follows from 7.4.
\QED

\section{A uniqueness theorem and automorphisms
on simple \CA s with $TR(A)\le 1$}

{\bf 8.1. Definition} Let $A$ and $B$ be \CA s.
Two \hm s $\phi, \psi: A\to B$ are said to be {\it stably unitarily
equivalent} if for any monomorphism $h: A\to B,$ $\ep>0$ and
finite subset ${\cal F}\subset A,$ there exists an integer $n>0$
and a unitary $U\in M_{n+1}({\tilde B})$ (or in
$M_{n+1}(B),$ if $B$ is unital) such that
$$
\|U^*diag(\phi(a),h(a),h(a),...,h(a))U-diag(\psi(a),h(a),h(a), ...,h(a))\|<\ep
$$
for all $a\in {\cal F},$ where $h(a)$ is repeated $n$ times
on both diagonals.

Let $A$ and $B$ be \CA s and $\phi,\,\,\psi: A\to B$ be (linear) maps.
Let ${\cal F}\subset A$ and $\ep>0.$ We write
\vspace{-0.05in}
$$
\phi\sim_{\ep}\psi\,\,\,\,\,\, {\rm on}\,\,\, {\cal F},
$$
if there exists a unitary $u\in B$ such that
$$
\|ad(u)\circ \phi(a)-\psi(a)\|<\ep
\,\,\,\,{\rm for\,\,\, all\,\,\,} a\in {\cal F}.
$$ We write
\vspace{-0.05in }
$$
\phi\approx_{\ep} \psi\,\,\,\,\, {\rm on}\,\,\, {\cal F},
\,\,\,{\rm if}\,\,\,\,\,
\|\phi(a)-\psi(a)\|<\ep
\,\,\,{\rm for\,\,\, all\,\,\,} a\in {\cal F}.
$$

\vspace{0.1in}

{\bf 8.2. Definition} Let $A$ be a \CA.

(i) Denote by ${\bf P}(A)$ the set of all projections and
unitaries in $M_{\infty}({\widetilde{A\otimes C_n}}),$ $n=1,2,...,$
where $C_n$ is an abelian \CA\, so that
$$
K_i(A\otimes C_n)=K_*(A;{\bf Z}/n{\bf Z}).
$$
One also has the following exact sequence
$$
\begin{array}{ccccccc}
K_0(A) &\to & K_0(A, {\bf Z}/k{\bf Z}) & \to & K_1(A)\\
\uparrow_{\bf k} & & & & \downarrow_{\bf k}\\
K_0(A) & \leftarrow & K_1(A, {\bf Z}/k{\bf Z})
&\leftarrow & K_1(A)\\
\end{array}
$$
(see \cite{Sc}).
As in \cite{DL2}, we use the notation
$$
{\underline K}(A)=\oplus_{i=0,1, n\in {\bf Z}_+} K_i(A;{\bf Z}/n{\bf Z}).
$$
By
$Hom_{\Lambda}({\underline K}(A),{\underline K}(B))$
we mean all \hm s from ${\underline K}(A)$ to ${\underline K}(B)$
which respect the direct sum
decomposition and the so-called Bockstein operations (see \cite{DL2}).
Denote by $Hom_{\Lambda}({\underline K}(A),{\underline K}(B))^{++}$
those $\alpha\in Hom_{\Lambda}({\underline K}(A),{\underline K}(B))$
with the property that
$\alpha(K_0(A)_+\setminus\{0\})\subset K_0(B)_+\setminus \{0\}.$
It follows from \cite{DL2} that if $A$ satisfies the Universal
Coefficient Theorem, then $Hom_{\Lambda}({\underline K}(A),{\underline K}(B))
\cong KL(A,B).$
Moreover, one has the following short exact sequence,
$$
0\to Pext(K_*(A), K_*(B))\to KK(A,B)\to KL(A,B)\to 0
$$

A separable \CA\, $A$ is said to satisfies Approximate
Universal Coefficient Theorem (AUCT) if
$$
KL(A,B)\cong Hom_{\Lambda}({\underline K}(A),{\underline K}(B))
$$
for any $\sigma$-unital \CA\, $B$ (see \cite{Lnuct}).
A separable \CA\, $A$ which satisfies the UCT must satisfy
the AUCT. If $A$ satisfies the AUCT, for convenience, we will use
$KL(A,B)^{++}$ for $Hom_{\Lambda}({\underline K}(A),{\underline K}(B))^{++}.$

(ii) Let $L: A\to B,$ be a \morp. We also use $L$ for the
extension from $A\otimes {\cal K}\to B\otimes {\cal K}$ as well as
maps from ${\widetilde{A\otimes C_n}}\to {\widetilde{B\otimes
C_n}}$ for all $n.$ Given a projection $p\in {\bf P}(A),$ if $L:
A\to B$ is an ${\cal F}$-$\dt$-multiplicative \morp\, with
sufficiently large ${\cal F}$ and sufficiently small $\dt,$
$\|L(p)-q'\|<1/4$ for some projection $q'.$ Define $[L](p)=[q']$ in
${\underline K}(B).$ It is easy to see this is well defined (see
\cite{Ln2}). Suppose that $q$ is also in ${\bf P}(A)$ with
$[q]=k[p]$ for some integer $k.$ By adding sufficiently many
elements (partial isometries) in ${\cal F},$ we can assume that
$[L](q)=k[L](p).$ Similarly as in 6.1, one can do the same for
unitaries. Let ${\cal P}\subset {\bf P}(A)$ be a finite subset. We
say $[L] |_{\cal P}$ is well defined if $[L](p)$ is well defined
for every $p\in {\cal P}$ and if $[p']=[p]$ and $p'\in {\cal P},$
then $[L](p')=[L](p).$ This always occurs if ${\cal F}$ is
sufficiently large and $\dt$  is sufficiently small. In what
follows we write $[L]|_{\cal P}$ when $[L]$ is well defined on
${\cal P}.$

(iii)
Let $A=\oplus_{i=1}^nA_i,$ where each $A_i$ is a unital
\CA. Suppose that $L:A\to B$ is a ${\cal G}$-$\ep$-multiplicative
\morp. For any $\eta>0,$ if ${\cal G}$ is large enough and
$\ep$ is small enough, we may assume that
$$
\|L(1_{A_i})-p_i\|<\eta,\,\, \|p_jL(1_{A_i})\|<\eta \andeqn \|L(1_{A_i})p_j\|<\eta
$$
for some projection $p_i\in B$ and $i\not=j.$
Let $b=p_1L(1_{A_1})p_1.$ Then, with sufficiently small
$\eta,$ we may assume that $b$ is invertible in $p_1Bp_1.$
Define $L_1(a)=b^{-1/2}p_1L(a)p_1b^{-1/2}.$
Then $L_1(1_{A_1})=p_1.$ Consider
$(1-p_1)L(1-p_1).$ It is clear that, for any $\dt>0,$ by induction and choosing
a sufficiently large ${\cal G}$ and sufficiently small
$\eta$ and $\ep,$
$$
\|L-\Psi\|<\dt,
$$
where $\Psi(a)=\oplus_{i=1}^n L_i(1_{A_i}a1_{A_i})$ for $a\in A.$
So, to save notation in what follows, we may assume that
$L=\oplus_{i=1}^nL_i,$ where each $L_i: A_i\to B$ is a completely
positive contraction which maps $1_{A_i}$ to a projection in $B$
and $\{L_1(1_{A_1}), L_2(1_{A_2}),...,L_n(1_{A_n})\}$ are mutually
orthogonal.

\vspace{0.1in}

{\it Throughout the rest of this section, ${\bf A}$ denotes
the class of separable nuclear \CA s satisfying the
Approximate Universal Coefficient
Theorem.}

\vspace{0.1in}

{\bf 8.3. Lemma} (Theorem 4.4 in \cite{LnA}) {\it
Let $B$ be a unital \CA\,
and let $A$ be a unital
\CA\,in ${\bf A}$
which is a unital $C^*$-subalgebra
of $B.$  Let
$\alpha: A\to B$ and $\beta: A\to B$ be two  \hm s.
Then $\alpha$ and $\beta$ are stably approximately unitarily
equivalent if  $[\alpha]=[\beta]$ in $KK(A,B)$
and if $A$ is simple or
$B$ is simple.
}

\vspace{0.1in}

The following is a modification of Theorem 4.8 in \cite{LnA}. A
proof was given in the earlier version of this paper. Since then a
more general version of the following appeared in \cite{Lnuct}
(Theorem 3.9). We will omit the original proof and view the
following as a special case of Theorem 3.9 in \cite{Lnuct}

\vspace{0.1in}

{\bf 8.4 .  Theorem} (cf. Theorem 4.8 in \cite{LnA})
{\it Let $B$ be a  \CA\, with stable rank one
and
$cel(M_m(B))\le k$ for some $k\ge \pi$ and for all $m,$
and let $A$ be a unital simple
\CA\,in ${\bf A}$
which is a  $C^*$-subalgebra
of $B.$   Let
$\alpha: A\to B$ and $\beta: A\to B$ be two  \hm s.
Then $\alpha$ and $\beta$ are stably approximately unitarily
equivalent if  $[\alpha]=[\beta]$ in $KL(A,B).$}

\vspace{0.1in} The following uniqueness theorem is a  modification
of Theorem 5.3 in \cite{LnA}.

\vspace{0.1in}

{\bf 8.5. Theorem} (5.3 in \cite{LnA}) {\it
 Let $A$ be a unital simple \CA\, in ${\bf A}$  and
 ${\bf L}: U(M_{\infty}(A))\to {\bf R}_+$ be a map.
 For any $\ep>0$ and any finite subset ${\cal F}\subset A$ there
 exist a positive number $\dt>0,$ a finite subset ${\cal G}\subset A,$
 a finite subset ${\cal P}\subset {\bf P}(A)$
 and an integer $n>0$ satisfying the following:
 for any unital simple  \CA\, $B$ with $TR(B)\le 1,$
 if $\phi,\, \psi,\, \sigma: A\to C$
 are three ${\cal G}$-$\dt$-multiplicative contractive completely
 positive linear
 maps with
 $$
 [\phi]|_{\cal P}=[\psi]|_{\cal P},
 $$
$$
cel({\tilde \phi}(u)^*{\tilde \psi}(u))\le {\bf L}(u)
$$
for all $u\in U(A)\cap {\cal P}$
and $\sigma$ is unital,
 then there is a unitary $u\in M_{n+1}(B)$ such that
 $$
 \|u^*diag(\phi(a),\sigma(a),\cdots,\sigma(a))u
 -diag(\psi(a),\sigma(a),\cdots,\sigma(a))\|<\ep
 $$
 for all $a\in {\cal F},$ where $\sigma(a)$ is repeated $n$ times.
}

{\it Proof:}
Suppose that the theorem is false.
Then there are $\ep_0>0$ and a finite subset ${\cal F}\subset A$
such that there are a sequence of positive numbers $\{\dt_n\}$ with
$\dt_n\downarrow 0,$ an increasing sequence of finite subsets
$\{{\cal G}_n\}$ whose union is dense in the unit ball of $A,$ a
sequence of finite subsets $\{{\cal P}_n\}$ of
${\bf P}(A)$
with $\overline{\bigcup_{n=1}^{\infty}{\cal P}_n}={\bf P}(A)$
and with $U_n=U(A)\cap{\cal P}_n,$
a sequence $\{k(n)\}$ of integers ($k(n)\nearrow \infty$)
and sequences $\{\phi_n\},$
$\{\psi_n\}$ and $\{\sigma_n\}$ of ${\cal G}_n$-$\dt_n$-multiplicative
positive linear maps from $A$ to $B_n$ with
$[\phi_n]|_{{\cal P}_n}=[\psi]|_{{\cal P}_n}$
and
$$
cel({\tilde \phi_n}(u)^*{\tilde \psi_n}(u))\le {\bf L}(u)
$$
for all $u\in U_n$
satisfying the following:
$$
\inf\{\sup\{\|u^*diag(\phi_n(a),\sigma_n(a),\cdots,\sigma_n(a))u
-diag(\psi_n(a), \sigma_n(a),\cdots,\sigma_n(a))\|: a\in {\cal F}\}\}\ge \ep_0
$$
where $\sigma_n(a)$ is repeated $k(n)$ times and
the infimum is taken over all unitaries in $M_{k(n)+1}(B_n).$

Set $D_0=\oplus_{n=1}^{\infty}B_n$ and $D=\prod_{n=1}^{\infty}B_n.$
Define $\Phi,\, \Psi,\,\Sigma: A\to D$ by
$\Phi(a)=\{\phi_n(a)\},$ $\Psi(a)=\{\psi_n(a)\}$ and
$\Sigma(a)=\{\sigma_n(a)\}$ for $a\in A.$ Let $\pi: D\to D/D_0$ be
the quotient map and
set ${\bar \Phi}=\pi\circ \Phi,$ ${\bar \Psi}=\pi\circ\Psi$
and ${\bar \Sigma}=\pi\circ \Sigma.$
Note that ${\bar\Phi}, {\bar \Psi}$ and ${\bar \Sigma}$ are
monomorphisms.
For any $u\in U_k,$
$$
cel({\tilde \phi_n}(u)^*{\tilde \psi_n}(u))\le {\bf L}(u)
$$
for all sufficiently large $n$ ($>k$).
This implies that there is an equi-continuous path $\{v_n(t)\}$
($t\in [0,1]$) such that
$$
v_n(0)={\tilde \phi_n}(u) \andeqn v_n(1)={\tilde \psi_n}(u)
$$
(see, for example, Theorem 1.1 of \cite{GL2}). Therefore, we
conclude that
$$
[{\bar \Phi}]|_{K_1(A)}=[{\bar \Psi}]|_{K_1(A)}.
$$

Given an element $p\in {\cal P}_k\setminus U_k$ (for some $k$),
we claim that
$$
[{\bar\Phi}(p)]=[{\bar \Psi}(p)].
$$
We have (see Proposition 2.1 in \cite{GL2})
$$
K_0(\prod B_n)=\prod^b K_0(B_n)
\andeqn
K_0(D/D_0)=\prod^b K_0(B_n)/\oplus K_0(B_n),
$$
where $\prod^bK_0(B_n)$ is the sequences of elements
$\{[p_n]-[q_n]\},$ where $p_n$ and $q_n$ can be represented by
projections in $M_L(B_n)$ for some integer $L.$ Since each
$TR(B_n)\le 1,$ $B_n$ has stable rank one and $K_0(B_n)$ is weakly
unperforated. By Proposition 2.2 in \cite{GL2}, each $B_n$ has
$K_i$ -divisible rank $T$ with $T(n,k)=1.$ By 6.5,
$cer(M_k(B_n))\le 4$ for all $k$ and $n,$ and the kernel of the
map from $K_1(\prod_n B_n)$ to $\prod^b_n K_1(B_n)$ is divisible
and torsion free (see 6.12). By the proof of part (2) of Theorem 2.1 in
\cite{GL2}, we also have
$$
K_i(\prod_n B_n,{\bf Z}/m{\bf Z})\subset \prod_n K_i(B_n, {\bf Z}/m{\bf Z}),
\,\,\,m=2,3,....
$$
(In fact, by 6.10, each $B_n$ has exponential length divisible
rank $E$ with $E(L,k)=8/\pi+L/k+1$ so that Theorem 2.1 (2) in
\cite{GL2} can be applied directly. See also part (2) of Corollary
2.1 in \cite{GL2}.))

 Since
$[\phi_n(p)]=[\psi_n(p)]$ in $K_0(B_n)$ or in $K_i(B_n, {\bf
Z}/m{\bf Z})$ ($i=0,1,$ $m=2,3,...$) for large $n,$
$$
[{\bar\Phi}(p)]=[{\bar \Psi}(p)].
$$
Then ${\bar\Phi}_*={\bar \Psi}_*.$
Therefore $[\bar\Phi]=[\bar \Psi]$ in $KL(A, \prod_n B_n/\oplus_n B_n).$

By applying 8.4, we obtain an integer $N$ and a unitary $u\in M_{N+1}(D/D_0)$
such that
$$
\|u^*diag({\bar \Phi}(a),{\bar \Sigma}(a),\cdots, {\bar \Sigma}(a))u
-diag({\bar \Psi}(a),{\bar \Sigma}(a).\cdots, {\bar \Sigma}(a))\|<\ep_0/3
$$
for all $a\in {\cal F},$ where ${\bar \Sigma}(a)$ is repeated $N$ times.
It is easy to see (see 1.3 in \cite{Ln2} for example)
 there is a unitary $U\in M_{N+1}(D)$ such that
$\pi(U)=u$ and for each $a\in {\cal F}$ there exists $c_a\in M_{N+1}(D_0)$ such
 that
 $$
 \|U^*diag(\Phi(a),\Sigma(a).\cdots,\Sigma(a))U-
 diag(\Psi(a),\Sigma(a),\cdots, \Sigma(a))+c_a\|<\ep_0/3
 $$
 where $\Sigma(a)$ is repeated $N$ times.
 Write $U=\{u_n\},$ where $u_n\in M_{N+1}(B_n)$
 are unitaries.
Since $c_a\in M_{N+1}(D_0)$ and ${\cal F}$ is finite, there is $N_0>0$
such that for $n\ge N_0$
$$
\|u_n^*diag(\phi_n(a),\sigma_n(a),\cdots,\sigma_n(a))u_n-
diag(\psi_n(a),\sigma_n(a),\cdots,\sigma_n(a))\|<\ep_0/2
$$
for all $a\in {\cal F},$ where $\sigma_n$ is repeated $N$ times.
This contradicts the assumption that the theorem is false.
\QED

\vspace{0.1in}

{\bf 8.6. Theorem} {\it Let $A$ be a separable unital nuclear simple
\CA\, with $TR(A)\le 1$ satisfying the AUCT and let
${\bf L}: U(A)\to {\bf R}_+.$
Then for any $\ep>0$ and any finite subset ${\cal F}\subset A,$
there exist $\dt_1>0,$ an integer $n>0,$
a finite subset ${\cal P}\subset
{\bf P}(A),$ a finite subset $S\subset A$
satisfying the following:

{\rm (i)} there exist  mutually orthogonal projections
$q,p_1,...,p_n$  with $q\preceq p_1$ and $p_1,...,p_n$  mutually
unitarily equivalent, and there exists a $C^*$-subalgebra $C\in
{\cal I}$ with $1_C=p_1$ and unital $S$-$\dt_1/2$-multiplicative
\morp s $\phi_0: A\to qAq$ and $\phi_1: A\to C$ such that
$$
\|x-(\phi_0(x)\oplus(\phi_1(x),\phi_1(x),...,\phi_1(x)))\|<\dt_1/2
$$
for all $x\in S,$ where $\phi_1(x)$ is repeated $n$ times,

moreover, there exist a finite subset ${\cal G}_0\subset A,$ a
finite subset ${\cal P}_0$ of  projections in $M_{\infty}(C),$ a
finite subset ${\cal H}\subset A_{s.a.},$ $\dt_0>0$ and $\sigma>0$
(which depend on the choices of $C$),

for any unital simple \CA \, $B$ with $TR(B)\le 1$ and  any two
${\cal S}\cup {\cal G}_0$-$\dt$-multiplicative completely positive
linear contractions $L_1, L_2: A\to B$ for which the following
hold (with $\dt=\min\{\dt_1, \dt_0\}$):

\vspace{0.1in}

{\rm (ii)}\,\,\,\,\,\,\,\,\,\,\,\,\,\,\, $
[L_1]|_{{\cal P}\cup{\cal P}_0}=[L_2]|_{{\cal P}\cup {\cal P}_0},
$

\vspace{0.1in}

{\rm (iii)} \,\,\,\,\,\,\,\,\,\,\,\,\,\,\,$ |\tau\circ
L_1(g)-\tau\circ L_2(g)|<\sigma $ \,\,\,for all $g\in {\cal H}$
and $\tau\in T(A),$

\vspace{0.1in}

{\rm (iv)}\,\,\,\,\,\,\,\,\,\,\,\,\,\,\,
$e=L_1\circ\phi_0(1_A)=L_2\circ \phi_0(1_A)$
is a projection,

\vspace{0.1in}

{\rm (v)} \,\,\,\,\,\,\,\,\,\,\,\,\,\,\,
$
cel({\tilde L}_1(\phi_0(u))^*{\tilde L}_2((\phi_0(u))))\le {\bf L}(u)
$
(in $U(eBe)$) for all $u\in U(A)\cap {\cal P},$

\vspace{0.1in}

then there exists a unitary $U\in B$ such that
$$
ad(U)\circ L_1\approx_{\ep} L_2
\,\,\,{\rm on}\,\,\,{\cal F}.
$$
}

Note that (i) holds as long as $TR(A)\le 1$ and does not depend
on ${\bf L},$ $\ep$ and ${\cal F}.$

{\it Proof:} Since $TR(A)\le 1,$
 8.5 applies. Fix a finite subset ${\cal F}\subset A,$ $\ep>0$
and ${\bf L}.$
Let $\dt_1>0,$ ${\cal G}_1\subset A,$ ${\cal P}\subset {\bf P}(A)$ and
$n$ be as required by  8.5 corresponding
to ${\bf L},$ $\ep/4$ and ${\cal F}.$

Let $S={\cal G}_1,$ $\eta=\min(\dt_1/2, \ep/4).$
 Let  $q, p_1,...,p_n,$ $\phi_0$ and
$\phi_1$ satisfy (i) given by 5.5.
Let $\dt_0'>0,$ $\sigma_1>0,$ and a finite subset ${\cal G}_0'\subset
A$ (${\cal G}_2\supset \phi_1({\cal G}_1)$) be as required by 5.8
corresponding to $C,$ the finite subset $\phi_1(S)$ and $\eta.$
Let ${\cal  P}_0$ contain
a finite set of projections in $C$ which generates
$K_0(C).$

Now choose a finite subset ${\cal G}_0''\subset C$ and $\dt_2>0$ so that
any ${\cal G}'$-$\dt_2$-multiplicative \morp\,
$L$ from $C$ to any \CA\, gives a well defined
map from $K_0(C).$

Let ${\cal G}_0={\cal G}_0'\bigcup {\cal G}_0'',$
$\dt_0=\min(\dt_2,  \dt_0)$ and $\sigma=\sigma_1/2(n+1),$
and let ${\cal H}=\{ {a^*+a\over{2}}, {a-a^*\over{2i}}: a\in {\cal G}_0'\}.$
Let $L_i: A\to B$ be two ${\cal S}\cup {\cal G}_0$-$\dt$-multiplicative
\morp s ($i=1,2$) which satisfy (ii)-(v).
Note that
$[L_1]|_{K_0(C)}=[L_2]|_{K_0(C)}.$
Let $e_1=\phi_1(1_A).$ By assumption, $e_1$ is a projection in $B.$
By the choice of $\dt_2,$ $\sigma_1$ and ${\cal G}_2,$ applying 5.8,
we obtain a unitary $v\in e_1Be_1$ such that
$$
\|L_1(x)-v^*L_2(x)v\|<\eta
$$
for all $x\in \phi_1(S).$ Therefore,
$$
\|L_1\phi_1(a)-ad(v)\circ L_2\circ \phi_1(a)\|<\eta
$$
for all $a\in {\cal G}_1.$
To simplify notation, without loss of generality, we may assume
that $L_2\circ \phi_1=ad(v)\circ L_2\circ \phi_1.$

Now, by applying 8.5,
we have
$$
L_1\circ\phi_0\oplus (L_1\circ \phi_1,
L_1\circ\phi_1,...,L_1\circ\phi_1) \sim_{\ep/4}
L_2\circ\phi_0\oplus (L_1\circ \phi_1,
L_1\circ\phi_1,...,L_1\circ\phi_1)\,\,\,{\rm  on}\,\,\,{\cal F}.
$$
From the above ($\eta<\ep/4$), we obtain
$$
L_2\circ\phi_0 \oplus (L_1\circ \phi_1,
L_1\circ\phi_1,...,L_1\circ\phi_1) \sim_{\ep/2}
L_2\circ\phi_0\oplus (L_2\circ\phi_1,..., L_2\circ\phi_1)
\,\,\,{\rm on}\,\,\,{\cal F}.
$$ Therefore
$$
\hspace{1in}L_1\sim_{\ep} L_2\,\,\,\,\,\,\,{\rm  on}\,\,\,{\cal F}.
\,\,\,\,\,\,\,\,\,\,\,\,\,\,\,\,\,\,\,\,\,\,\,\,\,\,\,\,\,
\,\,\,\,\,\,\,\,\,\,\,\,\,\,\,\,\,\,\,\,\,\,\,\,\,\,\,\,\,\,\,\,\,\,\,\,\,\,\,\,\,\,\,\,\,\,\,\,\,\,\,\,\,\,\, \QED
$$
It turns out that when $K_1(A)$ is torsion the ``uniqueness
theorem" can be stated in a much easy way.

\vspace{0.1in}

 {\bf 8.7. Theorem} {\it Let $A$ be a
unital separable simple \CA\, with $TR(A)\le 1$ and with torsion
$K_1(A).$ For any $\ep>0$ and any finite subset ${\cal F}\subset
A$ there exist $\dt>0,$ $\sigma>0,$ a finite subset ${\cal
P}\subset {\bf P}(A)$ and a finite subset ${\cal G}\subset A$
satisfying the following: for any unital simple \CA\, $B$ with
$TR(B)\le 1,$ any two ${\cal G}$-$\dt$-multiplicative completely
positive linear contractions $L_1, L_2: A\to B$ with
$$
[L_1]|_{\cal P}=[L_2]|_{\cal P}
$$
and
$$
\sup_{\tau\in T(B)}\{|\tau\circ L_1(g)-\tau\circ L_2(g)|\}<\sigma
$$
for all $g\in {\cal G},$
there exists a unitary $U\in B$ such that
$$
ad(U)\circ L_1\approx_{\ep} L_2
\,\,\,{\rm on}\,\,\,{\cal F}.
$$
}

{\it Proof:}
Define ${\bf L}: U(A)\to {\bf R}_+$
as follows.
If $u\in U_0(A),$ we let ${\bf L}(u)=2(cel(u)+\pi/16);$
if $u\in U(A)\setminus U_0(A),$ but
$[u]$ has order $k>1$ in $K_1(A)$ (therefore $u^k\in U_0(A)$
since $A$ has stable rank one),
let
$
{\bf L}(u)=8\pi+{2cel(u^k)\over{k}}+\pi/16.
$

Let ${\cal G}_1={\cal G}({\cal F}, \ep/2),$
$\dt_1=\dt({\cal F},\ep/2),$ $n=n({\cal F}, \ep/2),$
$S_1=S({\cal F},\ep/2)$ and ${\cal P}_1={\cal P}({\cal F},\ep/2)$ be
as required
in Theorem 8.6  (corresponding to ${\cal F},$ $\ep/2$ and ${\bf L}).$
Choose a finite subset ${\cal G}_2$ and a positive number $\eta<\min(\ep/2,\dt_1/2)$
satisfying the following:
if $H_i: A\to B$ are both ${\cal G}_2$-$\eta$-multiplicative with
$$
H_1\approx_{\eta} H_2\,\,\,\,{\rm on}\,\,\, {\cal G}_2
$$
then $[H_1]|_{{\cal P}_1}=[H_2]|_{{\cal P}_1}.$

If $u\in U(A),$ denote by $k(u)$ the smallest integer for which
$[u^{k(u)}]\in U_0(A).$ Let $K=\max\{k(u): u\in {\cal P}_1\cap
U(A)\setminus U_0(A)\}.$ Let ${\cal P}_2={\cal P}_1\bigcup
\{u^{k(u)}: u\in U(A)\cap {\cal P}_1\}.$ Set $V=\{
v^1,v^2,...,v^{K},v\in {\cal P}_1\cap (U(A)\setminus U_0(A))\}.$
We may assume that ${\cal G}_2\supset V\bigcup {\cal G}_1\bigcup
S_1.$ Set $0<\dt_2<\min(\dt_1/2, \eta/4, \ep/4, 1/K^2 512).$

Since $A$ is a nuclear simple \CA\, with $TR(A)\le 1,$ from
Lemma 5.5,
 there exist a $C^*$-subalgebra
  $F\in {\cal I}$ of $A$  with $p_1=1_F$ and ${\cal G}_2$-$\dt_2$-multiplicative
  completely positive linear contraction $\phi_1: A\to F$ such that
  $$
   id_A(x) \approx_{\dt_2} qxq\oplus diag(\phi_1(x),\phi_1(x),...,\phi_1(x))
           \,\,\,{\rm for\,\,\,all }\,\,\, x\in {\cal G}_2,
                 $$
 where $\phi_1$ is repeated $n$ times and $(1-q)A(1-q)=M_n(p_1Ap_1).$
Set $\phi_0(x)=qxq$ for $x\in A.$
Now let ${\cal P}_0$, ${\cal G}_0,$ $\dt_0>0,$ $\sigma>0$
and finite subset ${\cal H}
\subset A_{s.a.}$ be as required in (i) of 8.6.

By 6.8 (with perhaps larger ${\cal G}_2$), we also  have
(in $U(qAq)$)
$$
cel({\tilde \phi_0}(u))\le cel(u)+\pi/128
$$
for $u\in U_0(A)\cap {\cal P}_1;$  and
$$
cel({\tilde \phi_0}(u^k))\le cel(u^k)+\pi/128 \andeqn
$$
$$
\|{\tilde \phi_0}(u^k)-{\tilde \phi_0}(u)^k\|<1/128
$$
for $u\in (U(A)\setminus U_0(A))\cap {\cal P}_1$ and
$[u]$ has order $k$ in $K_1(A).$

Let ${\cal G}\supset {\cal G}_2\bigcup \phi_0({\cal G}_2)\cup {\cal G}_0$
be a finite
subset  and $0<\dt<\dt_2.$
Let ${\cal P}={\cal P}_1\bigcup\{q, p_1\}.$
Suppose that $L_i$ are two unital ${\cal G}$-$\dt$-multiplicative
\morp s that satisfy
$$
[L_1]|_{\cal P}=[L_2]|_{\cal P}
\,\,\,\,{\rm and}\,\,\,
\sup_{\tau\in T(B)}\{|\tau\circ L_1(g)-\tau\circ L_2(g)|\}<\sigma/2
$$
for all $g\in {\cal H}.$

For any (finite subset) ${\cal G}_3\subset  {\cal G}_2\bigcup \phi_0({\cal G}_2)$
and $0<\dt_3<\dt_2/2,$
by considering $(L_i)|_{qBq\oplus (1-q)B(1-q)},$
from (iii) of 8.2, with possibly larger ${\cal G}$ and smaller $\dt,$
there are ${\cal G}_3$-$\dt_3$-multiplicative \morp s
$L_i': A\to B$ ($i=1,2$) such that
$$
\|L_i(g)-L_i'(g)\|<\dt_3\,\,\,\,{\rm for\,\,\, all}\,\,\, g\in {\cal G}_3
$$
and $L_i'(q)$ is a projection ($i=1,2$).
Furthermore, we may assume
that
$$
[L_1']|_{\cal P}=[L_2']|_{\cal P}\andeqn
|\tau\circ L_1'(g)-\tau\circ L_2'(g)|<\sigma
$$
for all $g\in {\cal H}$ and $\tau\in T(B).$
Replacing $L'_2$ by $ad(W)\circ L'_2$ for a suitable unitary
$W\in B,$ we may assume
that $e=L_1'(q)=L_1'(\phi_0(1_A))=L_2'(q).$

Let $\Lambda_i=L_i'\circ \phi_0.$
With sufficiently large ${\cal G}_3$
and sufficiently small $\dt_3$ (applying 6.8), we may assume that
$$
cel({\widetilde \Lambda_i}(u))\le cel(u)+\pi/64,
$$
$i=1,2,$ $u\in U_0(A)\cap {\cal P}_1;$
$$
\|{\widetilde \Lambda_i}(u^k)-{\tilde \Lambda_i}(u)^k\|<1/64
\andeqn cel({\widetilde \Lambda_i}(u^k))\le cel(u^k)+\pi/64
$$
for $u\in (U(A)\setminus U_0(A))\cap {\cal P}_1$ but $[u]$ has order $k$
in $K_1(A).$
Then, for $u\in U_0(A)\cap {\cal P}_1$
$$
cel({\tilde L_1'}(\phi_0(u))^*{\tilde L_2'}(\phi_0((u)))\le
2(cel(u)+\pi/64)\le {\bf L}(u);
$$
and, for $u\in (U(A)\setminus U_0(A))\cap {\cal P}_1$ with order $k,$
since
$$
cel({\tilde L_1'}(\phi_0(u^k))^*{\tilde L'}_2(\phi_0(u^k)))
<2[cel(u^k)+\pi/64],
$$
by 6.10,
$$
cel({\tilde L_1'}(\phi_0((u))^*{\tilde L'}_2(\phi_0((u)))\le
8\pi+{2[cel(u^k)+\pi/64]\over{k}}+\pi/64<{\bf L}(u).
$$

Now Theorem 8.6 provides a unitary $U\in B$ such that
 $$
 ad(U)\circ L_2'\approx_{\ep/2} L_1'
  \,\,\,\,{\rm on}\,\,\, {\cal F}.
$$
This implies that
$$
L_2\sim_{\ep} L_1\,\,\,\,\,\,\,{\rm on}\,\,\,  {\cal F}.
$$
 \QED

\vspace{0.1in}

The following is a characterization of approximately inner
automorphisms (for the case in which $K_1(A)$ is torsion).

\vspace{0.1in}

{\bf 8.8. Theorem} {\it Let $A$ be a unital nuclear simple
\CA \, with $TR(A)\le 1$ and  with torsion $K_1(A)$  which satisfies the AUCT.
Then an automorphism $\alpha: A\to A$ is approximately inner
if and only if $[\alpha]=[{\rm id}_A]$ in $KL(A,A)$
and $\tau\circ\alpha(x)=\tau(x)$ for
all $x\in A$ and $\tau\in T(A).$}

{\it Proof:} If $\alpha$ is approximately inner, then
it is clear that
$$
\tau\circ\alpha(x)=\tau(x)
$$
for all $x\in A$ and $\tau\in T(A).$ The  ``only if" part
follows from 4.5 in \cite{LnA}.
It is also clear that the ``if part" follows from
8.7.

\QED

\section{The existence theorems}

{\bf 9.1. Definition}
Let $A$ and $B$ be two unital stably finite
\CA s and let $\alpha: K_0(A)\to K_0(B)$ be a positive \hm\, and
$\Lambda: T(B)\to T(A)$ be a continuous affine map.
We say $\Lambda$ is compatible to $\alpha$ if
$\Lambda(\tau)(x)=\tau(\alpha(x))$ for all $x\in K_0(A),$
where we view $\tau$ as a state on $K_0(A).$
Let $S$ be a compact convex set. Denote by $Aff(S)$ the set
of all (real) continuous affine functions on $S.$
Let $\Lambda: S\to T$ be a continuous affine
map from $S$ to another  compact convex set  $T.$
We denote by $\Lambda_{\natural}: Aff(T)\to Aff(S)$
the unital positive linear continuous map defined by $\Lambda_{\natural}(f)(s)
=f(\Lambda(s))$ for $f\in Aff(T).$
A positive linear map $\xi: AffT(A)\to AffT(B)$ is said to be
compatible to $\alpha$ if
$\xi({\hat p})(\tau)=\tau(\alpha(p))$ for all $\tau\in T(B)$ and any projection
$p\in M_{\infty}(A).$
Let $A$ be a unital \CA\, (with at least one normalized trace).
Define $Q: A_{sa}\to AffT(A)$ by $Q(a)(\tau)=\tau(a)$ for $a\in A.$
Then $Q$ is a unital positive linear map.

   A \CA\, $A$ is said to be $KK$-{\it attendable}
   (for simple  \CA s with tracial rank no more than 1 ) if
      for any unital separable simple \CA \, $B$ with $TR(B)\le 1,$
  any $\alpha\in Hom_{\Lambda}
  (\underline{K}(A),\underline{K}(B))^{++}$ (see 8.2)
  and any finite subset ${\cal P}\subset {\bf P}(A)$
 with $[1_A]\subset {\cal P},$
there exists a sequence of
completely positive linear
 contractions $L_n: A\to B\otimes {\cal K}$
  such that
  $$
    \|L_n(a)L_n(b)-L_n(ab)\|\to 0 \andeqn  [L_n]|_{\cal P}=\alpha|_{\cal
    P}\,\,\,\,
     \rforal \,\,\,a,\,b\in A.
     $$

As in \cite{LnC},
       if for any $\ep>0$ and any finite subset ${\cal F}\subset
          A,$ there exists a $C^*$-subalgebra $A_1$ of $A$ which is
             $KK$-attendable such that
                $
                   {\cal F}\subset_{\ep} A_1,
                      $
then $A$ is $KK$-attendable.

A unital nuclear separable simple \CA\, $A$  with $TR(A)\le 1$ is
said to be  {\it pre-classifiable} if it satisfies the Universal
Coefficient Theorem and if $A$ is $KK$-attendable, and, in
addition to the above, for any  unital separable nuclear simple
\CA\, with $TR(A)\le 1$ and any continuous affine map
$\Lambda:T(B)\to T(A)$ compatible to $\alpha,$ \vspace{-0.05in}
$$
\sup_{\tau\in T(B)}\{|\Lambda(\tau)(a)-\tau\circ L_n(a)|\}\to 0
\,\,\,\rforal \,\,a\in A.
$$
Or, equivalently, for any contractive
positive linear map $\xi: AffT(A)\to AffT(B)$ compatible to
$\alpha,$ \vspace{-0.05in}
$$
\sup_{\tau\in T(B)}\{|\xi(Q(a))(\tau)-\tau\circ L_n(a)|\}\to
0\,\,\,\rforal \,\,\,a\in A_{sa}.
$$

{\bf 9.2 } If $h: A\to B$ is a unital \hm, then $h$ induces a
unital positive affine  map $h_{\sharp}: AffT(A)\to AffT(B).$ The
map $h_{\sharp}$ is contractive. Suppose that $Y$ is a compact
metric space and $P\in M_l(C(Y))$ is a non-zero projection with
constant rank. It is known and easy to see that
$$
AffT(PM_l(C(Y))P)=AffT(M_l(C(Y))=C_{\bf R}(Y).
$$

{\bf 9.3. Theorem} {\it Let $A$ be a simple unital \CA\, with at
least one tracial state. Then for any  affine function $f\in
Aff(T(A))$ with $\|f\|\le 1$ and any $\ep>0,$ there exists an
element $a\in A_{sa}$ with $\|a\|<\|f\|+\ep$ such that
$\tau(a)=f(\tau)$ for all $\tau\in T(A).$ Furthermore, if $f\ge
0,$ we can choose $a\ge 0.$}

{\it Proof:} We prove this using the results in \cite{CP}.
By 2.7 in \cite{CP}, we may identify
$T(A)$ with the real part of the
unit sphere of $(A^q)^*$ (see \cite{CP} for the notation).
By 2.8 in \cite{CP}, it suffices to consider those
$f\in Aff(T(A))_+.$ There is an element $b\in (A^q)^{**}$
such that $b(\tau)=f(\tau)$ for all $\tau\in T(A).$ Since
$f$ is (weak-*) continuous, $b\in A^q.$
Since $b(\tau)\ge 0$ for all $\tau\in T(A),$ by 6.4 in \cite{CP},
there is $c\in A_+$ and $z\in A_{sa}$ with
$\tau(z)=0$ for all $\tau\in T(A)$ (i.e., $z\in A_0$ using
the notation in \cite{CP})  such that $b=c+z.$
Now the theorem follows from
2.9 in \cite{CP}.
\QED

{\bf 9.4. Lemma} {\it Let $A$ be a separable unital \CA.
Let $\ep>0$ and ${\cal F}\subset A$ be a finite subset.
Then there exists $\dt>0$ and a finite subset ${\cal G}\subset A$
satisfying the following:
For any  unital separable \CA\, $C$ with at least one tracial state
and any
unital contractive positive linear maps $L: A\to C$
which is ${\cal G}$-$\dt$-multiplicative,
then, for any $t\in T(C)$
there is a trace $\tau\in T(A)$ such that
$$
|\tau(a)-t(L(a))|<\ep\,\,\, for\,\,\,all\,\,\, a\in {\cal
F}.
$$
}

{\it Proof:} Otherwise, there would be an $\ep_0>0$ and a finite
subset ${\cal F}\subset A,$ a sequence of unital separable
\CA\, $C_n,$ a sequence of unital contractive positive
linear map $L_n: A\to C_n$
such that
$$
\lim_{n\to\infty}\|L_n(a)L_n(b)-L_n(ab)\|=0\rforal a, b\in A,
$$
and a sequence
$t_{n}\in T(C_{n})$ such that
$$
\inf\{ \max \{|t(a)-t_{n}(L_{n}(a))|: a\in {\cal
F}\}: t\in T(A)\}\ge \ep_0
$$
for all $n.$
Let $s_{n}$ be a state of $A$ which extends $t_{n}\circ
L_{n}.$ Let $\tau$ be a weak limit of $\{s_{n}\}.$
So there is a subsequence $\{n_k\}$
such that $\tau(a)=\lim_{m\to\infty}s_{n_k}(a)$ for all $a\in A.$
 It is
routine to check that $\tau\in T(A).$ Therefore, there exists
$K>0,$ such that
$$
|\tau(a)-t_{n_k}(L_{n_k}(a))|<\ep_0/2
$$
for all $k\ge K.$ We obtain a contradiction. \QED

{\bf 9.5. Lemma} {\it Let $A=C(X),$ where $X$ is a
path connected finite CW-complex. Let $B$ be
a unital separable nuclear non-elementary simple \CA\,
with $TR(B)\le 1$ and
$\Lambda: T(B)\to T(A)$ be a  continuous affine map.
Then, for any $\sigma>0$ and any finite subset ${\cal H}\subset AffT(A),$
there exists a unital monomorphism $h: A\to B$ such that
the image of $h$ is in a $C^*$-subalgebra $B_0\in {\cal I}$ and
$$
\|h_{\sharp}(f)-\Lambda_{\natural}(f)\|<\ep
$$
for all $f\in {\cal H},$ where $h_{\sharp}, \,  \Lambda_{\natural}
: AffT(A)\to AffT(B)$ are the maps induced by
$h$ and $\Lambda,$ respectively.

Moreover, if
there is positive \hm\, $\alpha: K_0(A)\to K_0(B)$ with
$\alpha([1_A])=[1_B]$ and $\Lambda_{\natural}$ is compatible to $\alpha,$
then the above is also true for $A=PM_l(C(X))P,$ where $P\in M_l(C(X))$ is a projection  in $M_l(C(X)).$
Furthermore, if $X$ is contractible, we can also require that $h_{*0}=\alpha.$}

{\it Proof:} We will apply 2.5 of \cite{Li2}. We may assume that
the identity of $A$ is contained in ${\cal H}.$ Fix $\ep>0.$ By 9.3, for
each $f\in \Lambda_{\natural}({\cal H}),$ there is $a_f\in B_{sa}$
with $\|a_f\|\le \|f\|$ such that
$$
\sup_{\tau\in T(B)}\{|\tau(a_f)-\tau(f)|\}<\ep/32.
$$
Let ${\cal G}=\{a_f: f\in  \Lambda_{\natural}({\cal H})\}.$
Let $N(2)$ be the number described in 2.5 of \cite{Li2}
corresponding to $\ep/4$ and ${\cal H}.$

Since $TR(B)\le 1,$  by 4.10, for any $\dt>0$ and any finite
subset ${\cal G}_1\subset A,$ (we assume that ${\cal G}\subset
{\cal G}_1$), there exists a nonzero projection $p\in B$ and a
unital $C^*$-subalgebra $C\in {\cal I}$ with $1_C=p$ such that

(1) $\|[x,p]\|<\dt$ for all $x\in {\cal G}_1,$

(2) $pxp\in_{\dt} C$ for all $x\in {\cal G}_1,$

(3) $\tau(1-p)<\ep/64,$ and

(4) $C=\oplus_{i=1}^n C_i$ with $C_i=M_{m(i)}(C([0,1]))$ and
$m(i)\ge N(2).$

We have
$$
|\tau(pxp)-\tau(x)|<\ep/32\,\,\,{\rm for\,\,\,all}\,\,\, \tau\in
T(B)\andeqn\rforal\,\,x\in {\cal G}.
$$ Since $B$ is nuclear,
by 3.2 in \cite{LnC}, there exists a unital completely positive
contraction $L': pBp\to C$ such that
$$
\|L'(pxp)-pxp\|<\ep/32 \,\,\,\rforal x\in {\cal G}.
$$
Define $L:
B\to C$ by $L(b)=L'(pbp)$ for $b\in B.$ With sufficiently small
$\dt,$ we have
$$
\sup_{\tau\in T(B)}\{|t\circ L(x)-\tau(x)|:\tau\in T(B)\andeqn
t=\tau/\tau(p) \}<\ep/16
$$
for all $x\in {\cal G}.$
Choose an integer $N>0$ so that
$$
|g(t')-g(t')|<\ep/32\,\,\,{\rm if}\,\,\, |t'-t''|\le 1/N\,\,\,
{\rm for\,\,\, all}\,\,\, g\in \Lambda_{\natural}({\cal H}).
$$
Here we view $Aff(C)=\oplus_{i=1}^n C_{\bf R}([0,1])$ and $t'$ and
$t'$ are points in the same $[0,1].$ Let $0=t_{i,0}<t_{i,1}<\cdots
<t_{i,N}=1$ such that $t_{i,k}-t_{i,k+1}=1/N$ and $i$ indicates
the $i$th interval. We also view them as tracial states of $C,$
i.e., $t_{i,k}(f)=tr(f(t_{i,k}))$ for $f\in \oplus_{i=1}^nC([0,1],
M_{m(i)}),$ where $tr$ is the normalized trace on $M_{m(i)}.$

By applying Lemma 9.4, with sufficiently large ${\cal G}_1$ and
small $\dt,$ we may assume that there are $\tau_{i,k}\in T(B)$
such that
$$
|\tau_{i,k}(a)-t_{i,k}(L(a))|<\ep/32
$$
for all $a\in {\cal G}.$ Let $\Delta$ be the convex hull of
$t_{i,k}$ in $T(C).$ Define $\gamma_1:
\Delta\to T(B)$ by $\gamma_1(\sum_{i,k}\alpha_{i,k}t_{i,k})
=\sum_{i,k}\alpha_{i,k}\tau_{i,k}$ for
$\alpha_{i,k}\ge 0$ and  $\sum_{i,k}\alpha_{i,k}=1.$
Therefore
$$
|\gamma_1(t)(a)-t(L(a))|<\ep/32\,\,\,\rforal a\in {\cal G}
$$
and $t\in \Delta.$
For each
$\tau\in T(C),$ one has, for $f\in \oplus_{i=1}^n C([0,1],
M_{m(i)}),$ that
$$
\tau(f)=\sum_{i=1}^n\int tr(f(t))d\mu_i,
$$
where $\mu_i$ is a Borel measure on $[0,1].$ Define $\gamma_2:
T(C)\to \Delta$ by
$$
\gamma_2(\tau)(f)=\sum_{i=1}^n
\sum_{k=0}^N\alpha_{i,k}tr(f(t_{i,k}))
$$
for $f\in C,$ where $\alpha_{i,k}=\mu_i(A_k)$ and $A_k=[t_{i,k},
t_{i,k+1}),$ $k=0,1,2,...,N.$
Put $\gamma=\gamma_1\circ \gamma_2.$
It is clear that $\gamma$ is an affine continuous map from $T(C)$
into $T(B).$ Moreover,
$$
|\gamma(t)(a)-t(L(a))|<\ep/16
$$
for all $f\in {\cal G}$ and $t\in T(C).$


By 2.5 in \cite{Li2}, there is a unital \hm\, $h_1: A\to C$ such
that
$$
\|(h_1)_{\sharp}(g)-\gamma_{\natural}\circ\Lambda_{\natural}(g)\|<\ep/16
$$
for all $g\in {\cal H}.$
In the following estimation, if $\tau\in T(B),$ we denote ${\bar
\tau}=(1/\tau(p))\tau|_{C},$ where $C$ is regarded as a subalgebra
of $A.$ We also note that
$\Lambda_{\natural}(g)(\gamma({\bar \tau}))=\gamma({\bar\tau})(a_f)$
for $f=\Lambda_{\natural}(g)$ and $g\in Aff(T(A)).$
We estimate, that, for $g\in {\cal H},$ $f=\Lambda_{\natural}(g),$
and for any $\tau\in T(B),$
$$
|(h_1)_{\sharp}(g)({\bar \tau})- \Lambda_{\natural}(g)(\tau)|
$$
$$
\le |(h_1)_{\sharp}(g)({\bar \tau})-\gamma_{\natural}\circ
\Lambda_{\natural}(g)({\bar \tau})|+|\gamma({\bar
\tau})(a_f)-{\bar\tau}(L(a_f))|+|{\bar
\tau}(L(a_f)-\Lambda_{\natural}(g)(\tau)|
$$
for $g\in {\cal H}.$ Each of the first two terms of the right side
of the inequality is no more than $\ep/16.$ For the last term, we
note that
$$
\Lambda_{\natural}(g)(\tau)=\tau(a_f) \andeqn |{\bar
\tau}(L(a_f))-\tau(a_f)|<\ep/16
$$
for all $g\in {\cal H}.$ Thus we have
$$
|(h_1)_{\sharp}(g)({\bar \tau})- \Lambda_{\natural}(g)(\tau)|
<3\ep/16\,\,\,{\rm for\,\,\,all}\,\,\, g\in {\cal H}.
$$

Since $(1-p)B(1-p)$ is non-elementary and simple, by p.61 in
\cite{AS}, there exists a positive element $b\in (1-p)B(1-p)$ with
$sp(b)=[0,1].$ From this we know that there is a unital
$C^*$-subalgebra $C_0$ of $(1-p)B(1-p)$ such that $C_0\cong
C([0,1]).$ It is well known that there is a unital monomorphism
$h_2: A\to C_0.$ Finally, we let $h=h_1+h_2.$ It is clear that $h$
meets the requirements of the conclusion of the (first part of)
lemma.

For the second part, we note there is an integer $N>0$ and a
projection $e\in M_N(A)$ such that $eM_N(A)e\cong C(X).$ Let $d\in
B$ be a projection such that $\alpha([e])=[d].$ Since
$\Lambda_{\natural}$ is compatible to $\alpha,$
$\Lambda_{\natural}$ induces a unital positive map $\zeta:
AffT(eAe)\to AffT(dBd).$ From what we have proved, we obtain a
\hm\, $h_1: eAe\to dBd$ as required (with $\dt\cdot \ep,$ where
$\dt=\inf\{\tau(d): \tau\in  T(B)\}$). Set $h_2=h_1\otimes {\rm
id_{M_l}}: eAe\otimes M_l\to dAd\otimes M_l.$ Since $eAe\cong
C(X),$ we may assume that $P\in eAe\otimes M_l.$ Let $p\in
dAd\otimes M_l$ such that $p=h_2(P).$ Since $\Lambda_{\natural}$
is compatible to $\alpha,$ $[p]=[1_B].$ Since  $B$ has stable rank
one, $p$ is unitarily equivalent to $1_B.$ Therefore $p(dAd\otimes
M_l)p\cong A.$ The second part of the lemma follows.

For the last part of the lemma, we note that $PM_l(C(X))P=M_s(C(X))$ for
some integer $0<s\le l,$
since $X$ is contractible. Let $e_{11}$ be a minimal projection of $A.$
Choose a projection $d\in B.$ For any finite subset ${\cal H}_1\subset
(e_{11}Ae_{11})_{s.a},$ from what we have shown, we obtain
a \hm\, $h': e_{11}Ae_{11}\to dBd$ such that
$$
\|h'_{\sharp}(f)-\Lambda_{\natural}(f)\|<\ep/s.
$$
View $M_r(dBd)$ as a unital hereditary \SCA\, of $B.$
Put $h=h'\otimes {\rm id}_{M_r}.$ It is clear that $h$ meets
the requirements of the lemma.
\QED

\vspace{0.1in}

{\bf 9.6. Corollary} {\it Let $A\in {\cal I},$ $B$ be
a unital separable nuclear simple \CA\, with $TR(B)\le 1,$
$\gamma: K_0(A)\to K_0(B)$ be a positive \hm\, and
$\Lambda: T(B)\to T(A)$ be an affine continuous map which is compatible
to $\gamma.$ Then, for any $\sigma>0$ and any finite
subset ${\cal G}\subset A,$
there exists a unital monomorphism $\phi: A\to B$
such that
$$
\sup_{\tau\in T(B)}\{|\tau\circ\phi(g)-\Lambda(\tau)(g)|\}<\sigma
$$
for all $g\in {\cal G}$
and
$\phi_*=\gamma.$}

{\it Proof:} Note that 9.5 holds for $A=M_n.$
It is then clear that, by considering each summand of
$A,$ the corollary follows from 9.5.
\QED

\vspace{0.1in}

{\bf 9.7. Proposition} {\it Every $KK$-attendable, unital separable
nuclear simple  \CA\, $A$ with $TR(A)\le 1$ which satisfies the AUCT
is pre-classifiable.}

{\it Proof:} Let $A$ be a $KK$-attendable separable
nuclear simple \CA\, with $TR(A)\le 1$ satisfying  the AUCT and $B$ be
a unital nuclear separable simple \CA\, with $TR(B)\le 1.$

Let $\alpha\in Hom_{\Lambda}(\underline{K}(A), \underline{K}(B))^{++},$
${\cal P}\subset {\bf P}(A)$ be a finite subset containing
$[1_A],$ and  $\Lambda: T(B)\to T(A)$ be a continuous map
which is comparable to $\alpha|_{K_0(A)}.$
Suppose that $e\in B$ is a projection such that
$\alpha(1_A)=e.$
To save notation, without loss of generality, we may assume that
$B=e(B\otimes {\cal K})e.$
Let $\{\dt_n\}$ be a decreasing sequence of positive numbers
with $\lim_{n\to\infty}\dt_n=0.$ For each $n,$ since $A$ is a unital
simple \CA\, with $TR(A)\le 1,$
there are nonzero projections $p_n\in A$ and a $C^*$-subalgebra
$C_n\in {\cal I}$ with $1_{C_n}=p_n,$ and a sequence of unital completely
positive linear contractions $\Phi_n: A\to C_n$
such that

(1) $\|[x,p_n]\|<\dt_n,$

(2) $\|p_nxp_n-\Phi_n(x)\|<\dt_n$

(3) $\|x-(p_nxp_n\oplus \Phi_n(x))\|<\dt_n$
for all $x\in A$ with $\|x\|\le 1$ and

(4) $\tau(1-p_n)<1/2n$ for all $\tau\in T(A).$

Denote by $\Psi_n(x)=(1-p_n)x(1-p_n)+\Phi_n(x)$ (for $x\in A$).
Note that
$$
\|\Psi_n(ab)-\Psi_n(a)\Psi_n(b)\|\to 0
\andeqn \|\Phi_n(ab)-\Phi_n(a)\Phi_n(b)\|\to 0
$$
for all $a,\,b\in A$ as $n\to\infty.$

Since $A$ is $KK$-attendable, for each $n,$
there exists a sequence of completely positive linear contractions
$L_n: A\to B\otimes {\cal K}$
such that
$$
[\Psi_n]|_{\cal P}=[{\rm id}]|_{\cal P}, [L_n]|_{\cal
P}=\alpha|_{\cal P},\,\,\, [L_n\circ \Psi_n]|_{\cal
P}=\alpha|_{\cal P},
$$
$$
\|L_n\circ \Psi_n(ab)-L_n\circ \Psi_n(a)L_n\circ \Psi_n(b)\|\to 0
\andeqn\,\, \|L_n\circ \Phi_n(ab)-L_n\circ \Phi_n(a)L_n\circ
\Phi_n(b)\|\to 0
$$
as $n\to\infty$ for all $a,b\in A.$
Suppose that $C_n=\oplus_{i=1}^{t(n)}D_{n,i},$
where $D_{n,i}\cong M_{(n,i)}$ or $D_{n,i}\cong M_{(n,i)}(C([0,1])).$
Let $d_{n,i}={\rm id}_{D_{n,i}}.$ We may also assume
that, for each $n,$ $L_n(d_{n,i})$ is a projection (see 8.2 (iii)) and
$$
[L_n]([d_{n,i}])=\alpha([d_{n,i}])\,\,\,{\rm for \,\,\,all}\,\,\,
n, i.
$$

Let $\gamma_n: T(A)\to T(C_n)$ be defined by
$\gamma_n(\tau)=(1/\tau(p_n))\tau|_{C_n}.$
Let ${\cal G}_n$ be a finite subset (containing
generators) of
$C_n$ and let $\{d_n\}$ be a decreasing sequence
of positive numbers with $\lim_{n\to\infty}d_n=0.$ For large $n,$ by applying
9.6,
we obtain a \hm\, $h_n: C_n\to q_nBq_n,$ where
$[d_{n,i}]=\alpha([d_{n,i}]),$ such that
$$
|\tau\circ h_n(g)-\gamma_n\circ \Lambda(\tau)(g)|<d_n \,\,\,{\rm for\,\,\,all}
\,\,\, \tau\in T(B)\andeqn
$$
for all $g\in {\cal G}_n.$
Put
$\phi_n(x)=L_n\circ((1-p_n)x(1-p_n))+h_n\circ\Phi_n(x)$
for $x\in A.$
It is easy to see, by choosing a large $n,$
$\phi_n: A\to B$  meets the requirements of Definition 9.1.
\QED

\vspace{0.1in}

{\bf 9.8. Lemma} {\it  Let $A$ be a unital \CA,  $B$ be a unital
separable simple \CA\, with $TR(B)\le 1$ and $F\in {\cal I}$ be a
$C^*$-subalgebra of $B.$ Let $G$ be a subgroup generated by a
finite subset of ${\bf P}(A).$ Suppose that there is an ${\cal
F}$-$\dt$ -multiplicative \morp\, $\psi: A\to F\subset B$ such
that $[\psi]|_{G}$ is well defined. Then for any $\ep>0,$ there
exists  a finite dimensional $C^*$-subalgebra $C\subset B$ and an
${\cal F}$-$\dt$--multiplicative \morp\, $L: A\to C\subset B$ such
that
$$
[L]|_{G\cap K_0(A,{\bf Z}/k{\bf Z})}=[\psi]|_{G\cap K_0(A,{\bf Z}/k{\bf Z})},
\andeqn \tau(1_C)<\ep
$$
for all tracial states $\tau$ in $T(B)$
and for all $k\ge 1$ so that
$G\cap K_0(A,{\bf Z}/k{\bf Z})\not=\{0\},$
where $L$ and $\psi$ are viewed as maps to $B.$
Furthermore, if $[\psi]_{G\cap K_0(A)}$ is
positive, so is $[L]|_{G\cap K_0(A)}.$ }

{\it Proof:} This is a minor modification of the proof of 4.2 in \cite{LnC}.
Let $0<\ep<1.$ Without loss of generality, we may assume
that $F=C([0,1])\otimes M_n.$
Let $q_1\in F$ be a minimal projection.
Suppose that
$$
G\cap K_0(A, {\bf Z}/k{\bf Z})=\{0\}\,\,\,\,
{\rm for}\,\,\,k>K.
$$
By 5.5, with $m=2lK!+1$ and $1/l<\ep/n,$ we may write $q_1=q+\sum_{i=1}^mp_i,$
where $[q]\le [p_1],$ $q, p_1,...,p_m$ are mutually orthogonal
projections,  $[p_1]=[p_i],$ $i=1,2,...,m$ and $\tau(p_1)<1/2l<\ep/2n.$
Set $e_1=q+p_1$ and $q_0=\sum_{j=2}^{2l+1}p_j.$
Then $[e_1]+K![q_0]=[q_1]$ in $K_0(B)$ and
$\tau(e_1)<\ep/n$ for all tracial states $\tau$ on $B.$
From this we obtain a $C^*$-subalgebra $C$ of $B$
such that $C\cong M_n$ and its minimal projection
is equivalent to $e_1.$
In particular, $\tau(1_C)<\ep.$
Let $\phi: F\to M_n\to C$ be a unital \hm,
where the map $F\to M_n$ is a point-evaluation.
Let $L=\phi\circ \psi,$ $j_1: F\to B$ and
$j_2: C\to B$ be embeddings.
By the choice of $q_1,$ $[e_1]$ and $[q_1]$ have the same
image in $K_0(B)/kK_0(B)$ for $k=1,2,...,K.$
Therefore $(j_1)_*=(j_2\circ \phi)_*$ on $K_0(F, {\bf Z}/k{\bf Z})$
for all $k\le K.$
Since $K_1(F)=K_1(C)=0,$ by the six-term exact sequence
in 8.2 (see 1.6 of \cite{LnC}), both $[L]$ and $[\psi]$ map
$K_0(A, {\bf Z}/k{\bf Z})$ to $K_0(B)/kK_0(B)$ and factor through
$K_0(F,{\bf Z}/k{\bf Z}).$
Therefore
$$
[L]|_{G\cap K_0(A, {\bf Z}/k{\bf Z})}=
[\psi]|_{G\cap K_0(A, {\bf Z}/k{\bf Z})},\,\,\,k=1,2,...,K.
$$
The general case in which $F$ is a direct sum of $M_l(C([0,1])$ follows
immediately.
\QED

\QED

{\bf 9.9 Lemma} {\it Let $C=\oplus_{j=1}^nC_j,$ where each
$C_j=P_jM_{s(j)}(C(X_j))P_j,$ $P_j$ is a projection in
$M_{s(j)}(C(X_j))$ and $X_j$ is a path connected compact metric
space with finitely generated $K_1(C_j),$ $K_0(C(X_j))={\bf
Z}\oplus tor(K_0(C_j)),$ $K_1(C(X_j))$ and $K_0(C_j))\subset
\{(z,x): z\in {\bf N},\, {\rm or}\,\,\, (z,x)=(0,0)\}.$ Then $C$
is $KK$-attendable.}

{\it Proof:}
Clearly, by considering each summand separately, we may assume that
$C$ has only one summand.
It is also clear that one can reduce
the general case to the case in which  $C=M_k(C(X))$ and
$K_i(C)$ satisfies the condition described in the lemma.

Let $A$ be a unital simple \CA\, with $TR(A)\le 1$ and
$\alpha\in KL(C,A)^{++}.$
Let ${\cal P}\subset {\bf P}(A)$ be a finite subset
and $G$ be the subgroup generated by
${\cal P}.$
By \cite{DL1}, for any finite subset ${\cal G}\subset C$ and $\eta>0,$
there exists a  ${\cal G}$-$\eta/2$-multiplicative \morp\,
$\psi: C\to M_N(B)$ for some large integer $N$ such that
$$
[\psi]|_{G}=\alpha|_{G}+h|_{G}
$$
for some one-point evaluation (at $\zeta$) $h: C\to M_N(B).$
Since $TR(M_N(B))\le 1,$  for any $\sigma_1>0$ and
$\ep/N>\sigma_2>0,$ there exist a projection $p\in M_N(B)$
and a unital $C^*$-subalgebra $F\subset M_N(B)$ with
$F\in {\cal I}$ and
with $1_F=p$ such that

(1) there are ${\cal G}$-$\eta$-multiplicative \morp s
$L_1: C\to F$ and $L_2: C\to (1-p)M_N(B)(1-p)$
such that
$$
\|\psi(x)-L_1(x)\oplus L_2(x)\|<\sigma_1
$$
for all $x\in {\cal G},$ and

(2) $\tau(1-p)<\sigma_2$ for all $\tau\in T(M_N(B)).$

With sufficiently small $\sigma_1,$ we may assume
that
$$
[\psi]_{G}=[L_1]|_{G}+[L_2]|_{G}.
$$
Suppose that
$$
G\cap K_0(C, {\bf Z}/k{\bf Z})=\{0\}\,\,\,{\rm for}\,\,\,k>K.
$$

By Lemma 9.8, there exists  a projection $e\le p$ with
$\tau(e)<\sigma_2$ for all $\tau\in T(M_N(B))$ and a unital ${\cal
G}$-$\eta$-multiplicative \morp\, $L_1': C\to F_1,$ where $F_1$ is
a $C^*$-subalgebra of $pM_N(B)p$ with $1_{F_1}=e$ such that ${\rm
dim} F_1<\infty$ and
$$
[L_1']|_{G\cap K_0(C, {\bf Z}/k{\bf Z})}=[L_1]|_{G\cap  K_0(C, {\bf Z}/k{\bf Z})}
$$
for all $k$ so that
$G\cap K_0(C, {\bf Z}/k{\bf Z})\not=\{0\}.$
So, in particular,
$$
([L_1']+[L_2])|_{G\cap K_0(C, {\bf Z}/k{\bf Z})}=
\alpha|_{G\cap K_0(C, {\bf Z}/k{\bf Z})}+h|_{G\cap K_0(C, {\bf Z}/k{\bf Z})}.
$$

Without loss of generality, we may assume that $\alpha([1_C])=[1_B].$
Note that, with a small $\sigma_2,$
we have $[(1-p)+e] \le [1_B].$
Let $q\in pM_N(B)$ be a projection with $[q]=[1_B]-[(1-p)+e].$
By applying an inner automorphism on $M_N(B)$ if necessary,
without loss of generality, we may assume that
$p, e, q\le 1_B.$
Let $h_0: C\to qBq$ be the unital one-point evaluation (at $\zeta$)
(mapping $1_B$ to $q$).

We now define $\Psi=L_1'\oplus L_2\oplus h_0.$ Note if $D$ is a
\CA\, of finite dimension or $D\in {\cal I},$ then $K_1(D,{\bf
Z}/k{\bf Z})=0.$
 Since the images of $L_1'$ and $h_0$ are in some
$C^*$-algebras belong to ${\cal I},$
$$
([L_1']+[h_0])|_{G\cap K_1(C)}=0,\,\,\,
([L_1']+[h_0])|_{G\cap K_1(C,{\bf Z}/k{\bf Z})}=0
\andeqn
$$
$$([L_1']+[h_0])|_{G\cap tor(K_0(C))}=0.
$$
Hence we compute that
$$
[\Psi] |_{G\cap K_i(C)}=\alpha|_{G\cap K_i(C)},\,\,\,\,i=0,1,
\andeqn [\Psi]|_{G\cap K_1(C,{\bf Z}/k{\bf Z})}=\alpha_{G\cap K_1(C,{\bf Z}/k{\bf Z})}.
$$
The proof is complete if we can show, in addition,
that $[\Psi]|_{G\cap K_0(C,{\bf Z}/k{\bf Z})}=\alpha|_{G\cap K_0(C,{\bf Z}/k{\bf Z})}.$
We note that there is an (unnatural) splitting short exact sequence
$$
0\to K_0(C)/kK_0(C)\to K_0(C,{\bf Z}/k{\bf Z})\to Tor(K_1(C,{\bf Z}/k{\bf Z}))
\to 0.
$$
From $[\Psi]|_{G\cap K_0(C)}=\alpha|_{G\cap K_0(C)},$ we conclude that
$$
[\Psi]|_{G\cap K_0(C)/kK_0(C)}=\alpha|_{G\cap K_0(C)/kK_0(C)}.
$$
On the other hand, it is easy to compute that
$$
K_0(C,{\bf Z}/k{\bf Z})={\bf Z}/k{\bf Z}\oplus K_0(C_0(Y),{\bf Z}/k{\bf Z}),
$$
where $Y$ is the locally compact space formed by taking the point
$\zeta$ away from $X$ and the summand ${\bf Z}/k{\bf Z}\subset
K_0(C)/kK_0(C).$ Since both $h$ and $h_0$ are point-evaluation at
$\zeta,$  we have
$$
[h]|_{K_0(C_0(Y), {\bf Z}/k{\bf Z})}=[h_0]|_{K_0(C_0(Y),{\bf Z}/k{\bf Z})}=0.
$$
It follows that
$$
[\Psi]_{G\cap K_0(C, {\bf Z}/k{\bf Z})}
=\alpha|_{G\cap K_0(C, {\bf Z}/k{\bf Z})}.
$$
This proves the  lemma.

\QED

\vspace{0.1in}

Let $X$ be a connected compact metric space and let $A=C(X).$
Denote by $\rho_A; K_0(A)\to {\bf Z}$ the rank  map. In 9.9,
we assume that ${\rm ker}\rho_{C(X_j)}$ is a finite group. One
should note that the proof of 9.9 does not work when ${\rm
ker}\rho_{C(X_j)}$ contains an infinite cyclic subgroup. This
happens
because $L_1'$ would not kill ${\rm ker}\rho_{C(X_j)}.$ In 9.10
below, we apply a result of L. Li and a result of Villadsen to
avoid this problem.

\vspace{0.1in}

Recall a \CA\, $A$ is said to be locally AH if
for any $\ep>0$ and any finite subset ${\cal F}\subset A,$
there exists a  $C^*$-subalgebra $B\subset A$
with $B=PM_l(C(X))P$ for some
compact metric space $X$ and where $P$ is a projection in $M_l(C(X)),$
such that
$$
dist(x, B)<\ep\,\,\,\,\,\,{\rm for\,\,\, all}\,\,\,x\in {\cal F}.
$$

{\bf 9.10. Proposition} {\it  Let $A$ be a separable unital simple
\CA\, with $TR(A)\le 1.$ If $A$ is locally AH, then $A$ is
pre-classifiable.}

{\it Proof:} It follows from \cite{Lnuct} that $A$ satisfies the
AUCT. We may assume that $A={\overline{\cup_{n=1}^{\infty} A_n}},$
where each $A_n$ is a finite direct sums of
$P_{n,i}M_{r(n,i)}C(X_{n,i})P_{n,i}$ and $X_{n,i}$ is a path
connected finite CW complex. One may assume that $1_{A_n}=1_A.$
Put $j_n: A_n\to A$ the embedding. Consider $j_n\times \alpha.$ If
$A_n$ has only one summand, then $K_0(A_n)={\mathbb Z}\oplus {\rm
ker}\rho_{A_n}.$ Since $\alpha\in KL(A,B)^{++},$ $(j_n\times
\alpha) \in KL(A_n,B)^{++}.$ By considering each summand
separately, we may assume $A_n$ has only one summand. Since $A$ is
simple, by 9.7 and  9.1, it suffices to show that, $A=C(X)$ is
$KK$-attendable for every path connected finite CW complex $X.$

Let $\alpha\in KK(A,B)^{++}.$
Suppose that $\alpha(1_A)=[p]\,(\not=0),$ where $p\in M_l(B)$
is a projection.

Fix a unital nuclear simple \CA\, $B$ with $TR(B)\le 1.$ By
\cite{V}, there is a unital simple \CA\, $C$ which is direct limit
of \CA s in 9.9 such that
$$
(K_0(C), K_0(C)_+, [1_C], K_1(C))=
(K_0(B), K_0(B)_+, [1_A], K_1(B)).
$$
By \cite{RS}, there exists $\beta\in KK(C,B)$ which gives
the above isomorphism.

Let $\alpha\in KL(A,B)^{++}$ and $\gamma=\alpha\times \beta^{-1}
\in KL(A,C)^+.$ Since $K_i(C(X))$ is finitely generated,
$KL(A,C)=KK(A,C).$ In particular, $\gamma(K_0(A)_+\setminus
\{0\})\subset K_0(C)_+\setminus \{0\}.$ By \cite{Li3}, there is a
\hm\, $h: A\to pM_l(C)p$ such that $[h]=\gamma.$ But by 9.9, since
each \CA\, described in 9.9 is $KK$-attendable, $C$ is
$KK$-attendable (see 9.1). Let $\ep>0$ and fix  finite subsets
${\cal F}\subset A$ and ${\cal P}\subset {\bf P}(A).$ Let ${\cal
G}=h({\cal F})\subset C$ and ${\cal Q}=[h]({\cal P})\subset {\bf
P}(C).$ Let $\Lambda: C\to B$ be a ${\cal G}$-$\ep$-multiplicative
\morp\, such that
$$
[\Lambda]|_{\cal Q}=\beta|_{\cal Q}.
$$
Define $L=\Lambda\circ h.$
Then $L: A\to B$ is a ${\cal F}$-$\ep$--multiplicative
\morp\, such that
$$
[L]|_{\cal P}=\alpha|_{\cal P}.
$$
So $A$ is $KK$-attendable.
\QED

{\bf 9.11. Lemma} {\it Let $A$ be a unital separable \CA, $\{{\cal
F}_k\}$ be an increasing sequence of finite subsets of the unit
ball of $A$ such that $\cup_k{\cal F}_k$ is dense in the unit ball
of $A,$  and let $\Phi_n: A\to A$ be a sequence of unital \morp s
such that $\lim_{n\to\infty}\|\Phi_n(a)\|=\|a\|$ for all $a\in A$
and
$$
\sum_{k=n}^{\infty}\|\Phi_k(ab)-\Phi_k(a)\Phi_k(b)\|<\sum_{k=n}^{\infty}\dt_n,
$$
for all $a,b\in {\cal G}_n,$ and for any finite subset ${\cal P}\subset {\bf P}(A),$
$[\Phi_n]|_{\cal P}=[{\rm id}]|_{\cal P} $ for all sufficiently large $n,$
where ${\cal G}_1={\cal F}_1,$ ${\cal G}_{n+1}\supset
\cup_{k=1}^n\Phi_k({\cal F}_n)\cup{\cal F}_n\cup\Phi_n({\cal G}_{n}),$ $n=1,2,...,$ and
 where $\sum_{n=1}^{\infty}\dt_n<\infty.$
Let $B=\lim_{n\to \infty} (A, \Phi_n)$ be the generalized inductive limit
{\rm in the sense of \cite{BK1})}.
Then $\{\Phi_n\}$ induces an isomorphism:
$$
(K_0(B), K_0(B)_+, [1_B], K_1(B))=(K_0(A), K_0(A)_+, [1_A], K_1(A)).
$$
}

{\it Proof}:  The proof is standard. We sketch here. Write
$K_i(A)=\bigcup_{n=1}^{\infty} G_n^{(i)},$ where each $G_n^{(i)}$
is a finitely generated subgroup of $K_0(A).$ Let
$\Phi_{n,n+m}=\Phi_{n+m}\circ \Phi_{n+m-1}\circ\cdots \circ
\Phi_n$ and $\Psi_n: A\to B$ be the map induced by the inductive
system which maps the $n$-th $A$ to $B.$ For each $G_n^{(i)},$ we
may assume that $ [\Psi_m]|_{G_n^{(i)}} $ is well defined for all
$m\ge n.$ The assumption that $[\Phi_n]|_{\cal P}=[{\rm
id}_A]|_{\cal P}$ for all sufficiently large $n$ implies that
$[\Psi_m]|_{G_n^{(i)}}=[\Psi_{m'}]|_{G_n^{(i)}}$ for all $m, m'\ge
n.$ This gives a \hm\, $\beta_i: K_i(A)\to K_i(B)$ ($i=0,1$).

Suppose that $p_1, p_2,v\in M_l(B)$ such that $v^*v=p_1$ and
$vv^*=p_2.$ There is a sequence $\{\Psi_{n_k}(a_k)\},$ where
$a_k\in M_l(A),$ such that it converges to $v.$ Since $v^*v=p_1,$
we have $\Psi_{n_k}(a_k^*a_k)\to p_1$ and $\Psi_{n_k}(a_ka_k^*)\to
p_2.$ Therefore we may assume that
$$
\|\Psi_{n_k}(a_k^*a_k-(a_k^*a_k)^2)\|<1/2^{k+1}
\andeqn \|\Psi_{n_k}(a_k^*a_k)-p_1\|<1/2^{k+1}.
$$
Since $\|\Psi_m(x)\|=\limsup\|\Phi_{m,n}(x)\|$ for all $x\in A$ and
$m\ge 1,$
by passing to a subsequence and possibly replacing
$a_k$ by $\Phi_{n_k,m_k}(a_k),$ $\Psi_{n_k}$ by $\Psi_{m_k},$
if necessary, we may assume that
$$
\|a_k^*a_k-(a_k^*a_k)^2\|<1/2^k, \,\,\, k=1,2,....
$$
It is standard that there is a partial isometry $v_k$ and a
projection $q_k\in A$ such that
$$
v_k^*v_k=q_k\andeqn \|v_k-a_k\|<1/2^{k-1}
$$
for all large $k.$
Let $q_k'=v_kv_k^*.$
Note also, for any $\ep>0,$
we have
$$
\|\Psi_{n_k}(q_k)-p_1\|<\ep \andeqn \|\Psi_{n_k}(q_k')-p_2\|<\ep
$$
for all large $k.$ Hence $[\Psi_{n_k}](q_k)=[p_1]$ and
$[\Psi_{n_k}](q_k')=[p_2].$ This, in particular, implies that
$[p_1]$ is in the image of $\beta_0.$ It follows that  $\beta_0$
is surjective. Note also that $[q_k]=[q_k']$ in $K_0(A).$ It
follows that
 $\beta_0$ is also injective. It is also easy to check
 from the definition that $\beta_0$ preserves the order.

In the above, if we let $w^*w=p_1$ and $ww^*\le p_2,$ then exactly
the same argument shows that there  are partial isometries $v_k\in
A$ such that $v_kv_k^*=q_k,$ $v_kv_k^*\le q_k'$ and
$\Psi_{n_k}(v_k)\to v,$ $\Psi_{n_k}(q_k)\to p_1$ and
$\Psi_{n_k}(q_k')\to vv^*\le p_2.$ These imply that $\beta_0$ is
an
 order
isomorphism.

A similar argument shows that $\beta_1$ is an isomorphism and
$K_1(A)=K_1(B).$ \QED

\vspace{0.05in}

{\bf 9.12 Theorem} {\it Let $A$ be a unital separable nuclear
simple \CA\, with $TR(A)\le 1$ satisfying the AUCT. Then there
exists a unital separable nuclear simple \CA\, $B$ with $TR(B)=0$
satisfying AUCT and the following:

{\rm (1)} $(K_0(A), K_0(A)_+, [1_A], K_1(A))
=(K_0(B), K_0(B)_+, [1_B], K_1(B)),$

{\rm (2)} there exists a sequence of \morp s $\Phi_n: A\to B$
such that

{\rm (i)} $\lim_{n\to\infty}\|\Phi_n(ab)-\Phi_n(a)\Phi_n(b)\|=0$ for
$a, b\in A,$

{\rm (ii)} For each finite subset ${\cal P}\subset {\bf P}(A)$
there exists an integer $N>0$ such that
$$
[\Phi_n]|_{\cal P}=[\alpha]|_{\cal P}
$$
for all $n\ge N,$
where $\alpha\in KL(A,B)$ which gives
an identification in {\rm (1)} above.
}

{\it Proof:} Let $\{{\cal P}_n\}$ be an increasing sequence of
finite subsets of ${\bf P}(A)$ such that the union is dense in
${\bf P}(A).$ In particular the union of the subgroups generated
by the images of ${\cal P}_n$ in $\underline{K}(A)$ is
$\underline{K}(A).$ Let $\{{\cal F}_n\}$ be an increasing sequence
of finite subsets of the unit ball of $A$ whose union is dense in
the unit ball of $A.$  Let $\{\dt_n\}$ be a decreasing sequence of
positive numbers such that $\sum_{n=1}^{\infty}\dt_n<\infty.$
Without loss of generality, we may assume that any ${\cal
F}_n$-$\dt_n$-multiplicative \morp\, $L_n$ defined on $A$
well-defines $[L_n]|_{{\cal P}_n}.$ Furthermore we may assume that
$({\cal F}_n, \ep_n, {\cal G}_n, {\cal P}_n,\dt_n)$ (for some
finite subsets ${\cal G}_n$ and ${\cal P}_n$) forms a 5-tuple as
defined in 6.8 of \cite{Lnuct} for $l=1,$ $b\ge \pi$ and $M=1.$ We
may also assume that ${\cal F}_n\subset {\cal G}_n.$

Let ${\cal G}_1'$ be the union of ${\cal G}_1$ and $\{(a+a^*)_+,
(a+a^*)_-,(a-a^*)_+,(a-a^*)_-: a\in {\cal F}_1\}.$ Since $A$ is
simple, if $b\in ({\cal G}_1')_+$ is a nonzero element, there are
$x_i(b)\in A$ such that
$$
\sum_{i=1}^{n(b)}x_i(b)^*bx_i(b)=1_A.
$$
Let ${\cal G}_1''$ be the union of ${\cal G}_1'$ and $\{x_i(b),
x_i(b)^*: b\in ({\cal G}_1')_+\}.$ Since $TR(A)\le 1,$ there
exists a \SCA\, $C_1\in {\cal I}$ with $1_{C_1}=p_1$ such that

($i_1$) $\|p_1a-ap_1\|<\dt_1/4$ for all  $a\in {\cal G}_1'',$

($ii_1$) $p_1ap_1\in_{\dt_1/4} C_1$ and $\|(1-p_1)a(1-p_1)\|\ge
(1-\dt_1/4)\|a\|$ for all  $a\in {\cal G}_1''$ (see 5.6) and

($iii_1$) $\tau(1-p_1)<\dt_1/4$ for all $\tau\in T(A).$

We may also assume  that

($iv_1$)
$\|\sum_{i=1}^{n(b)}(1-p_1)x_i(b)^*(1-p_1)b(1-p_1)x_1(b)(1-p_1)-(1-p_1)\|
<1/8$ and there are $z_i(b), b'\in C_1$ such that
$\|b-b'\|<\dt_1/2$ and
$\|\sum_{i=1}^{n(b)}z_i(b)^*b'z_i(b)^*-p_1\|<1/8$ for all $b\in
({\cal G}_1')_+.$

Write $C_1=M_{d(1,1)}(C([0,1]))\oplus \cdots \oplus M_{d(s_1,1)}(C([0,1])
\oplus M_{d(s_1+1,1)}\oplus \cdots \oplus M_{d(t_1,1)}.$
Put $F_1=M_{d(1,1)}\oplus \dots \oplus M_{d(t_1,1)}.$
Note that ${\rm dim}F_1<\infty.$
Define $\pi_1: C_1\to F_1$ to be a surjective point-evaluation
map.
We identify $F_1$ with the unital \SCA\, of $C_1$
(scalar matrices). So $F_1\subset C_1\subset p_1Ap_1.$
Note $\underline{K}(C_1)=\underline{K}(F_1)$ and
$[\pi_1]$ gives the identification.
Since $A$ is nuclear, there is a \morp\, $L_1': p_1Ap_1\to C_1$
such that
$$
\|L_1'(p_1ap_1)-p_1ap_1\|<\dt_1/4\,\,\,{\rm for\,\,\, all}\,\,\,
a\in {\cal G}_1''.
$$
Define $L_1(a)=\pi_1\circ L_1'(p_1ap_1)$ for $a\in A.$ Then $L_1$
is ${\cal G}_1''$-$\dt_1/2$-multiplicative. Define $\phi_1: A\to
A$ by $\phi_1(a)=(1-p_1)a(1-p_1)\oplus L_1(a)$ for $a\in A.$ It is
clear that $\phi_1$ is ${\cal F}_1$-$\dt_1/2$-multiplicative.
Furthermore,
$$
[\phi_1]|_{{\cal P}_1}=[{\rm id}_A]|_{{\cal P}_1},
\|\phi_1(a)\|\ge (1-\dt_1/2)\|a\| \,\,\,{\rm for}\,\,\, a\in {\cal F}_1
\andeqn
\|\sum_{i=1}^{n(b)}\phi_1(x_i(b))^*\phi_1(b)\phi_1(x_i(b))-1_A\|<1/4
$$

Let ${\cal G}_2'$ be the union of ${\cal G}_2,$ $\phi_1({\cal
G}_1'')$ and $\{(a+a^*)_+,(a+a^*)_-, (a-a^*)_+, (a-a^*)_-: a\in
{\cal F}_2\bigcup \, \phi_1({\cal G}_1'')\}.$ Since $A$ is simple,
for each nonzero positive element $b\in {\cal G}_2',$ there are
$x_1(b),...,x_{n(b)}(b)\in A$ such that
$$
\sum_{i=1}^{n(b)}x_i(b)^*bx_i(b)=1_A.
$$

Now let ${\cal G}_2''$ be a finite subset of $A$ containing ${\cal
G}_2',$ $\{1-p_1, p_1\},$ a generating set of $F_1$ and $\{x_i(b),
x_i(b)^*: b\in ({\cal G}_2')_+\}.$ Let ${\cal Q}_2$ be the finite
subset which is the union of ${\cal P}_2,$ $1-p_1,$ $p_1$ and
contains at least one minimal projection of each summand of $F_1,$
By choosing a possibly larger ${\cal G}_2''$ and smaller $\dt_2$
we may assume that any ${\cal G}_2''$-$\dt_2$-multiplicative
\morp\, $L$ well-defines $[L]|_{{\cal Q}_2}.$ Since $TR(A)\le 1,$
there is a \SCA\, $C_2\in {\cal I}$ with $1_{C_2}=p_2$ such that

($i_2$) $\|p_2a-ap_2\|<\dt_2/4$ for all $a\in {\cal G}_2'',$

($ii_2$) $p_2ap_2\in_{\dt_2/4} C_2$ and $\|(1-p_2)a(1-p_2)\|\ge
(1-\dt_2/4)\|a\|$ for all $a\in {\cal G}_2''$ (see 5.6),

($iii_2$) $\tau(1-p_2)<\dt_2/4$ for all $\tau\in T(A),$ and

($iv_2$)
$\|\sum_{i=1}^{n(b)}(1-p_2)x_i(b)^*(1-p_2)b(1-p_2)x_i(b)(1-p_2)-(1-p_2)\|
<1/32$ and there are $z_i(b), b'\in C_2$ such that
$\|b-b'\|<\dt_1/2$ and
$\|\sum_{i=1}^{n(b)}z_i(b)^*b'z_i(b)^*-p_2\|<1/32$ for all $b\in
({\cal G}_2')_+.$

Write $C_2=M_{d(1,2)}(C([0,1]))\oplus \cdots M_{d(s_2,2)}(C([0,1]))
\oplus M_{d(s_2+1,2)}\oplus \cdots \oplus M_{d(t_2,2)}.$
Put $F_2=M_{d(1,2)}\oplus \dots \oplus M_{d(t_2,2)}.$
Note that ${\rm dim}F_2<\infty.$
Define  $\pi_2: C_2\to F_2$ to be a surjective point-evaluation
map.
We identify $F_2$ with the unital \SCA\, of $C_2$
(scalar matrices). So $F_2\subset C_2\subset p_2Ap_2.$
Note that $\underline{K}(C_2)=\underline{K}(F_2)$
 and $[\pi_2]$ gives the identification.
Since $A$ is nuclear, there is a \morp\, $L_2': p_2Ap_2\to C_2$
such that
$$
\|L_2'(p_2ap_2)-p_2ap_2\|<\dt_1/4\,\,\,{\rm for\,\,\, all}\,\,\,
a\in {\cal G}_2''.
$$
Define $L_2(a)=\pi_2\circ L_2'(p_2ap_2)$ for $a\in A.$
Then $L_2$ is ${\cal F}_2$-$\dt_2/2$-multiplicative.
Define $\phi_2: A\to A$ by $\phi_2(a)=(1-p_2)a(1-p_2)\oplus L_2(a)$
for $a\in A.$ It is clear that $\phi_2$ is
${\cal F}_2$-$\dt_2/2$-multiplicative.
Furthermore,
$$
[\phi_2]|_{{\cal Q}_2}=[{\rm id}_A]|_{{\cal Q}_2},
\|\phi_2(a)\|\ge (1-\dt_2/4)\|a\|\,\,\, {\rm (}a\in {\cal
G}_2''{\rm )} \andeqn
\|\sum_{j=1}^n(b)\phi_2(x_j(b))\phi_2(b)\phi_2(x_j(b))-1_A\|<1/16.
$$

Since ${\rm dim} F_1<\infty,$ we may also assume
that there exists an injective \hm\,
$h_2: F_1\to A$ such that
$$
\|h_2-(L_2)|_{F_1}\|<1/4
$$
(see for example 2.3 of \cite{LnC2}).
We continue the construction of $L_n, \pi_n$ and $\phi_n$
in this fashion.

Let ${\cal G}_{n+1}'$ be the union of ${\cal G}_{n+1},$  $
\phi_1({\cal G''}_n),...,\phi_n({\cal G''}_n)$ and
$\{(a+a^*)_+,(a+a^*)_-, (a-a^*)_+, (a-a^*)_-: a\in {\cal
F}_n\bigcup_{k=1}^n  \phi_k({\cal G}_n'')\}.$ Since $A$ is simple,
for each nonzero positive element $a\in {\cal G}_{n+1}',$ there
are $x_1(a),...,x_{n(a)}(a)\in A$ such that
$$
\sum_{i=1}^{n(a)}x_i(a)^*ax_i(a)=1_A.
$$

Now let ${\cal G}_{n+1}''$ be a finite subset of $A$ containing
${\cal G}_{n+1}',$ $\{1-p_n, p_n\},$ a generating set of $F_n$ and
$\{x_i(b), x_i(b)^*: b\in ({\cal G}_{n+1}')_+\}.$ Let ${\cal
Q}_{n+1}$ be the finite subset which is the union of ${\cal
P}_{n+1},$ $1-p_{n+1},$ $p_{n+1}$ and at least one minimal
projection of each summand of $F_n.$ By choosing a possibly larger
${\cal G}_{n+1}''$ and smaller $\dt_n$ we may assume that any
${\cal G}_{n+1}''$-$\dt_n$-multiplicative \morp\, $L$ well-defines
$[L]|_{{\cal Q}_n}.$

Since $TR(A)\le 1,$ there is a \SCA\, $C_{n+1}\in {\cal I}$ with
$1_{C_{n+1}}=p_{n+1}$ such that

($i_{n+1}$) $\|p_{n+1}a-ap_{n+1}\|<\dt_{n+1}/4$ for all $a\in
{\cal G}_n'',$

($ii_{n+1}$) $p_{n+1}ap_{n+1}\in_{\dt_{n+1}/4} C_{n+1}$ and
$\|(1-p_{n+1})a(1-p_{n+1})\|\ge (1-\dt_{n+1}/4)\|a\|$ for all
$a\in {\cal G}_n'',$

($iii_{n+1}$) $\tau(1-p_{n+1})<\dt_n/4$ for all $\tau\in T(A)$ and

($iv_{n+1}$)
$\|\sum_{i=1}^{n(b)}(1-p_{n+1})x_i(b)^*(1-p_{n+1})b(1-p_{n+1})x_i(b)(1-p_{n+1})-(1-p_{n+1})\|
<1/2^{n+3}$ and there are $z_i(b), b'\in C_2$ such that
$\|b-b'\|<\dt_{n+1}/2$ and
$\|\sum_{i=1}^{n(b)}z_i(b)^*b'z_i(b)^*-p_{n+1}\|<1/2^{n+3}$ for all
$b\in ({\cal G}_{n+1}')_+.$

Write $C_{n+1}=M_{d(1,n+1)}(C([0,1]))\oplus \cdots \oplus M_{d(s_{n+1},n+1)}(C([0,1])
\oplus M_{d(s_{n+1}+1,n+1)}\oplus \cdots \oplus M_{d(t_{n+1},n+1)}.$
Put $F_{n+1}=M_{d(1,n+1)}\oplus \dots \oplus M_{d(t_{n+1}, n+1)}.$
Note that ${\rm dim}F_{n+1}<\infty.$
Define $\pi_{n+1}: C_{n+1}\to F_{n+1}$ to be a surjective point-evaluation map.
We identify $F_{n+1}$ with the unital \SCA\, of $C_{n+1}$
(scalar matrices). So $F_{n+1}\subset C_{n+1}\subset p_{n+1}Ap_{n+1}.$
Note that $\underline{K}(C_{n+1})=\underline{K}(F_{n+1})$
 and $[\pi_{n+1}]$ gives the identification.
Since $A$ is nuclear, there is a \morp\, $L_{n+1}': p_{n+1}Ap_{n+1}\to C_{n+1}$
such that
$$
\|L_{n+1}'(p_{n+1}ap_{n+1})-p_{n+1}ap_{n+1}\|<\dt_n/4\,\,\,{\rm for\,\,\, all}\,\,\, a\in {\cal F}_{n+1}.
$$
Define $L_{n+1}(a)=\pi_{n+1}\circ L_{n+1}'(p_{n+1}ap_{n+1})$ for
$a\in A.$ Then $L_{n+1}$ is ${\cal
F}_{n+1}$-$\dt_{n+1}/2$-multiplicative. Define $\phi_{n+1}: A\to
A$ by $\phi_{n+1}(a)=(1-p_{n+1})a(1-p_{n+1})\oplus L_{n+1}(a)$
Then $\phi_{n+1}$ is ${\cal G}_{n+1}''$-$\dt_n$-multiplicative.
Furthermore,
$$
[\phi_{n+1}]|_{{\cal Q}_{n+1}}=[{\rm id}_A]|_{{\cal Q}_{n+1}},
\|\phi_{n+1}(a)\|\ge (1-\dt_{n+1}/2)\|a\|\,\,\,{\rm
for\,\,\,all}\,\,\, a\in {\cal G}_{n+1}''
$$
$$
\andeqn
\|\sum_{j=1}^{n(b)}\phi_{n+1}(x_j(b))\phi_{n+1}(b)
\phi_{n+1}(x_j(b))-1_A\|<1/2^{n+1}.
$$
Again, we may also assume that there exists an injective
\hm\, $h_{n+1}: F_n\to A$ such that
$$
\|h_{n+1}-(L_{n+1})|_{F_n}\|<1/2^{n+1}.
$$

We then define $B=\lim_n(A, \phi_n).$ This is a generalized
inductive limit in the sense of \cite{BK1}(but $\{\phi_n\}$ is in
fact asymptotically multiplicative). Note that $B$ is a unital
separable \CA. By 5.13 in \cite{BK1}, $B$ is nuclear.
From ($iv_{n+1}$) and the construction
above, it is easy to check that
$B$ is simple.
Since for each $n$
$$
[\phi_{n+1}]|_{{\cal Q}_{n+1}}=[{\rm id}_A]|_{{\cal Q}_{n+1}},
$$
by 9.11, we have that
$$
(K_0(B), K_0(B)_+, [1_B], K_1(B))
=(K_1(A), K_1(A)_+, [1_A], K_1(A)).
$$
From the construction it is also standard to show
that $TR(B)=0$ (see for example the proof of 4.3 in \cite{LnB}
and also \cite{LnC2} ).
Let $\phi_{n, \infty}: A\to B$ be the  map
from the $n$th $A$ to $B$ induced by the inductive limit system.
 Set also $\Phi_n=\phi_{n,\infty}.$
It is clear $\{\Phi_n\}$ satisfies (i) and (ii). To see that $B$
satisfies AUCT, we note that $B$ satisfies the property (P)
described in 6.8 of \cite{Lnuct}. It follows from 6.13 and 6.16 of
\cite{Lnuct} that $B$ also satisfies the AUCT.

\QED

\section{The classification theorem}

To establish the classification of separable nuclear simple \CA s
with tracial rank  no more than one,
we first present, for each given set of invariant,
an example of  a simple AH-algebra which has the invariant set.
Then we will show every separable simple nuclear
\CA\, $A$ with $TR(A)\le 1$ and with the same set
of invariant which satisfies the (AUCT) is
in fact isomorphic to one of the  examples of the above mentioned simple
AH-algebra.

The following theorem was
proved by J. Villadsen. The extra conditions (1) and (2) are not
new either. It has appeared implicitly in several places
including in Villadsen's proof.

If $X$ is a convex set, the extremal points of $X$ is denoted
by $\partial_e(X).$

{\bf 10.1 Theorem} (cf. \cite{V})
{\it
Suppose that $G$ is a countable, partially ordered abelian group
which is simple, weakly unperforated with the Riesz interpolation property,
that $G/{\rm tor}(G)$ is non-cyclic, $u\in G_+,$ $H$ is a countable abelian
group, $\Delta$ is a metrizable Choquet simplex and $\lambda:\Delta\to S(G,u)$
is a continuous affine map with $\lambda(\partial_e\Delta)=\partial_eS(G,u).$
Then there is a simple AH-algebra $A=\lim_{n\to\infty}(A_n, h_n)$
with $TR(A)\le 1$
and with  $A_n=C_1\oplus C_2\oplus \cdots \oplus C_{m(n)},$
where $C_1$ is of the form as described in 7.1 (a single summand)
 and $C_j$ is of the form $C([0,1])\otimes M_{m(j)}$
(for $j>1$), such that

{\rm (1)} $h_n=h_n^{(0)}\oplus h_n^{(1)}\oplus h_n^{(2)},$
where $h_n^{(0)},h_n^{(1)}$ factors through a \CA\, in ${\cal I},$
and $h_n$ is injective, in particular, $h_n^{(0)}$ is homotopically
trivial,

{\rm (2)} $\tau\circ h_{n+1,\infty}\circ h_n^{(0)}(1_{A_n})\to 0$ uniformly
on $T(A),$

{\rm (3)}
 $\tau\circ h_{n+1,\infty}\circ h_n^{(2)}(1_{A_n})\to 0$
uniformly on $T(A),$

{\rm (4)} $(h_n)_{*1}$ is injective and

{\rm (5)}
 $(K_0(A), K_0(A)_+, [1_A], K_1(A), T(A), r_A)=
(G,G_+,u, H, \Delta, \lambda).$

}

{\it Proof}:
The proof of this is a combination of Villadsen's proof of
the main theorem in \cite{V} and the proof of 1.5 in \cite{Lnan}.
Let $A=\lim_{n\to\infty}(A_n,h_n)$ be as in 1.5 in \cite{Lnan}.
This algebra $A$ satisfies  (3), (4)
and
$(K_0(A), K_0(A)_+, [1_A], K_1(A))=(G, G_+, u, H).$
Moreover each $A_n$ can be chosen so it has  the form as required.
Here one needs one modification and one explanation.
In the proof of 1.5 in \cite{Lnan} we use $C_j=M_{m(j)}$ for
$j>1.$ But we can  map $C([0,1], M_{m(j)})$ into
$M_{m(j)}$ by a one point-evaluation and then map $M_{m(j)}$
into $C([0,1], M_{m(j)})$ (as constant functions). So we can assume
that $A_n$ has the required form. Note also that the new $A_n$ has
the same $K$-theory as the old one.
If $K_1(A)=F=\lim_{n\to\infty}(F_n,\gamma_n),$
in the proof of 1.5 in \cite{Lnan}, $K_1(A_n)=F_n$ and the map
$\Phi_{n,n+1}$ in the proof of 1.5 has the property $(\Phi_{n,n+1})_{*1}
=\gamma_n.$ However, since $F$ is a countable abelian group, one
can always assume that $F_n$ is finitely generated and $\gamma_n$ is
injective (by choosing $F_n$ as subgroups and $\gamma_n$ as
embeddings) so that (4) holds.

We will revise the map $h_n$ to meet the other requirements.
Villadsen's proof  in \cite{V} is to replace
$h_n$ by $\phi_n$ without changing its $K$-theory in
such a way that one gets $\Delta$ as tracial space and $\lambda$
as paring.  We will follow his proof with a minor
modification.
Each $h_n$ may be written as $h_n'\oplus h_n'',$
where $h_n''$ is a point-evaluation, as in 1.5 of \cite{Lnan}.
Note that the Lemma in \cite{V} holds when $X_q^i$ is a compact
connected CW complex with dimension at least one but no more than
three.
Following Villadsen's proof,
by applying the Lemma in \cite{V} and its proof,
one can replace $h_n''$ to achieve exactly what the Lemma
in \cite{V} achieved.
It should be noted that
Villadsen's proof of the main theorem in \cite{V}
 works when $X_q^j$ has lower dimension (but
at least one), since the required maps $i_q^j: [0,1]\to X_q^j$ and
$k_q^j: X_q^j\to [0,1]$ still exist.
The new map obtained from Villadsen's proof has the form
${\tilde \psi}_n=h_n'\oplus
{\tilde h_n''},$  where ${\tilde h_n''}$ is homotopically trivial.
Furthermore, it can be chosen so that it factors through
a \CA\, in ${\cal I}.$
The construction of Villadsen then gives
a simple AH-algebra $B$ with $TR(B)\le 1$ and satisfies
(5). Moreover, one has $\tau\circ {\tilde \psi}_{n+1}\circ h_n(1_{A_n})\to 0$
uniformly on $T(B).$
 The construction does not change (3) and (4).
It is also easy to get (1) and (2).
For example, consider
$h_{n+1}'\circ h_n'\oplus h_{n+1}'\circ {\tilde h_n''}
\oplus {\tilde h_{n+1}''}\circ {\tilde \psi}_n.$
Note that
$h_{n+1}'\circ {\tilde h_n''}$ and ${\tilde h_{n+1}''}\circ {\tilde \psi_n}$ are
homotopically trivial and
$\tau\circ {\tilde \psi}_{n,\infty}\circ
h_{n+1}'\circ {\tilde h_n''}(1_{A_n})\to 0$ uniformly
on $T(B).$
\QED

 \vspace{0.1in}

{\bf 10.2 Definition} Let $C$ be a unital \CA. We denote by
$S_u(K_0(C))$ the set of states on $K_0(C),$ i.e., the set of
order and unit preserving \hm\, from $K_0(C)$ to (the additive
group) ${\mathbb R}.$ There is an affine map $\lambda: T(C)\to
S_u(K_0(C))$ such that $\lambda(t)([p])=t(p)$ for all projections
$p\in M_{\infty}(C)$ and $t\in T(C).$ Suppose that $C$ is stably
finite. It was proved in Theorem 6.1 of \cite{R2} (for simple case)
and Theorem 3.5 in \cite{BR} that each state in $S_u(K_0(C))$ is
induced by a quasitrace $t\in T(C).$ If $C$ is exact, or if
$TR(C)\le 1$ and simple, all quasitraces are traces (see (ix) in
4.9).

{\it  Let $A$ and $B$ be two unital \CA s. We say
$$\gamma:
(K_0(A), K_0(A)_+, [1_A], K_1(A), T(A))\to
(K_0(B), K_0(B)_+, [1_B], K_1(B), T(B))
$$
is an order isomorphism if there is an order isomorphism
$\gamma_0: (K_0(A), K_0(A)_+)\to (K_0(B), K_0(B)_+)$ which maps
$[1_A]$ to $[1_B],$ there is an isomorphism $\gamma_1: K_1(A)\to
K_1(B)$ and an affine  homeomorphism $\gamma_2: T(A)\to T(B)$ such that
$\gamma_2^{-1}(\tau)(x)=\tau(\gamma_0(x))$ for all $\tau\in T(B)$
and $x\in K_0(A),$ where we view $\tau$ as a state on $K_0(A).$
}
\vspace{0.1in}

{\bf 10.3 Theorem } {\it Let $A$ and $B$ be two unital separable nuclear
simple \CA s with $TR(A)\le 1$  and $TR(B)\le 1$
satisfying AUCT such that
$$
(K_0(B), K_0(B)_+, [1_B], K_1(B), T(B))
=(K_0(A), K_0(A), [1_A], K_1(A), T(A))
$$
in the sense of 10.2.
Then there is a sequence of \morp s
$\{\Psi_n\}$ from $A$ to $B$ such that

{\rm (i)}
 $\lim_{n\to\infty}\|\Psi_n(ab)-\Psi_n(a)\Psi_n(b)\|=0$ for all $a, b\in A,$

{\rm (ii)} for any finite subset set ${\cal P}\subset {\bf P}(A),$
       $$
        [\Psi_n]|_{\cal P}=\alpha|_{\cal P},
$$
for all sufficiently large $n,$ where $\alpha\in KL(A,B)^{++}$
(see 8.2) gives the above identification on
 $K$-theory and

{\rm (iii)} $$
\lim_{n\to\infty}\sup_{\tau\in T(B)}\{|\tau\circ \Psi_n(a)-\xi(Q(a))(\tau)|\}=0
$$
for all $a\in A_{s.a},$ where $\xi: AffT(A)\to AffT(B)$ is the
affine isometry given above.
}


 {\it Proof:} It follows from
Theorem 9.12 that there is a unital separable simple nuclear \CA\,
$C$ with $TR(C)=0$ which satisfying the AUCT such that
$$
(K_0(A), K_0(A), [1_A], K_1(A))=(K_0(C), K_0(C)_+, [1_C], K_1(C))
$$
and a sequence of \morp s $L_n: A\to C $ satisfying condition (2) in
9.12. In particular, $[L_n]|_{\cal P}=\beta|_{\cal P},$ for any
finite subset ${\cal P}$ and all sufficiently large $n,$ where
$\beta\in KL(A, C)^{++}$ which gives the above identification on
$K$-theory. It follows from \cite{Lnmsri} that there is a unital
separable simple AH-algebra $C_1$ such that $ C_1\cong C. $ To
simplify notation, we may assume that $C_1=C.$

It follows from 9.10 that there exists a sequence of \morp s
$\Phi_n': C\to B$ such that

(i') $\lim_{n\to\infty}\|\Phi_n'(ab)-\Phi_n'(a)\Phi_n'(b)\|=0$ for all $a, b\in C,$

(ii') for any finite subset ${\cal Q}\subset {\bf P}(C),$
$$
[\Phi_n']|_{\cal Q}=(\beta^{-1}\times \alpha)|_{\cal Q},
\,\,\,\,\,\,\,\,\,\rforal \,\,\, {\rm sufficiently\,\,\,
large}\,\,\, n.
$$

Thus by choosing a subsequence $\{k(n)\}$ and defining
$\Psi_n=\Phi_{k(n)}'\circ L_n: A\to B$ we see that  $\Psi_n$
satisfies (i) and (ii). (In fact one can show that  $A$ is $KK$-attendable.)
We then apply the proof of 9.7, to obtain a
(new) sequence $\{\Phi_n\}$ which also satisfies (iii).

\QED

Using the argument of \cite{Lnmsri}, Zhuang Niu
gives a different proof the above theorem.

\vspace{0.1in}

{\bf 10.4. Theorem} {\it Let $A$ and $B$ be two unital separable
nuclear simple \CA s with $TR(A)\le 1$ and $TR(B)\le 1$ satisfying
the AUCT. Suppose that
$\lambda(\partial_e(T(A)))=\partial_e(S_u(K_0(A)))$ and
$\lambda(\partial_e(T(B)))=\partial_e(S_u(K_0(B))).$ Then $A$ is
isomorphic to $B$ if and only if there exists an order isomorphism
$$\gamma=(\gamma_0, \gamma_1, \gamma_2):
(K_0(A), K_0(A)_+, [1_A], K_1(A), T(A))\to (K_0(B), K_0(B)_+,
[1_B], K_1(B), T(B)),$$ where
$\gamma_2^{-1}(\tau)(x)=\tau(\gamma_0(x))$ for all $\tau\in T(B)$
and $x\in K_0(A)$ {\rm (}see 10.2{\rm )}. }

{\it Proof:} Let $A$ be as in the theorem. Note that $T(A)$ is a
Choquet simplex (see 4.9 (ix) for example). By 10.1 and 10.3 (as
well as 4.8), there is a unital simple AH-algebra
$B=\lim_{n\to\infty}(B_n, \phi_{n, n+1})$ with
$$
(K_0(A), K_0(A)_+, [1_A], K_1(A), T(A))\cong
(K_0(B), K_0(B)_+, [1_B], K_1(B), T(B)).
$$
where $B_n$ is as described in 10.1. Here $\phi_{n,n+1}$ is a
\hm\, from $B_n$ to $B_{n+1}.$ Put $\phi_{n,m}=\phi_{m-1,m}\circ
\cdots \circ\phi_{n,n+1}.$ Denote by $\psi_n: B_n\to B$ the \hm\,
induced by the inductive system. As in 10.1, we will assume that
$\phi_{n,m}$ and $\psi_n$ are injective and $(\phi_{n,m})_{*1}$ is
injective for all $m>n.$ In what follows, when it is convenient,
we may identify $B_n$ with $\psi_n(B_n)$ without further warning.
We also assume that the inductive limit satisfies the condition
(1)-(5) in 10.1. To prove the theorem it suffices to prove that
$A$ is isomorphic to this specially constructed $B.$ In
what follows, $\kappa$ is the quotient map
from $U(D)/CU(D)$ to $K_1(D)$ for a \CA\, $D.$

Let $\xi: AffT(A)\to AffT(B)$ be the affine isometry induced by
the above identification. Since both $A$ and $B$ satisfy the AUCT,
there is  $\alpha\in KL(A,B)^{++}$ which gives the isomorphism
$(K_0(A), K_0(A)_+, [1_A], K_1(A))$ $\to (K_0(B), K_0(B)_+, [1_B],
K_1(B)).$ Let ${\bf L}: U(B)\to {\bf R}_+$ be defined as follows:
if $u\in U_0(B),$ ${\bf L}(u)=2cel(u)+8\pi+\pi/16;$ if $u\in
U(B)\setminus U_0(B)$ and $u^k\in U_0(B)$ (where $k$ is smallest
such positive integer) ${\bf L}(u)=16\pi +2cel(u^k)/k+\pi/16,$ and
if $u\in U(B)\setminus U_0(B)$ and $[u]$ is not of finite order in
$K_1(B),$ ${\bf L}(u)=16\pi+\pi/16.$ Fix any $\ep>0$
($\ep<\pi/128$) and finite subset ${\cal F}\subset B.$ Let
$\dt'>0,$ integer $n>0,$ finite subsets ${\cal P}\subset {\bf
P}(B),$  ${\cal S}\subset B$
be as required in Theorem 8.6 (corresponding to $B,$ ${\bf L},$
$\ep/2$ and ${\cal F}$). There are mutually orthogonal projections
$q, p_1,p_2,..,p_n$ with $q\lesssim p_i,$ and with $[p_i]=[p_1]$
($i=1,2,...,n$) and there is a \SCA\, $C_1'\in {\cal I}$ with
$1_{C_1'}=q$ and there are unital ${\cal S}$-$\dt'/4$-multiplicative
\morp s $h_0: B\to qBq$ and $h_1: B\to C_1'$ such that
$h_0(x)=qxq$ and
$$
\|x-(h_0(x)\oplus h_1(x)\oplus\cdots\oplus h_1(x))\|<\dt'/16
$$
for all $x\in S,$ where $h_1$ is repeated $n$ times. Put
$C=M_n(C_1')\subset (1-q)B(1-q).$ Let ${\cal P}_0,$ ${\cal G}_0,$
${\cal H},$ $\dt_0>0$ and $\sigma_1>0$ also be  as required by 8.6 and
let $\dt=\min\{\dt_0, \dt'\}.$

We may assume that ${\cal P}_0$ contains at least one minimal
projection of  each summand of $C$. Note that, without loss of
generality, we may also assume that $q$ commutes with each element
in ${\cal S}$ and ${\cal H}.$

Without loss of generality (by omitting a possible error of $\dt/16$),
 we may also assume
that  each $u\in U(B)\cap {\cal P}$ has the form
$quq\oplus (1-q)u(1-q),$ where $quq\in U(qBq)$ and $(1-q)u(1-q)
\in U(C).$
We may further assume that $q\in B_1$ and $quq\in U(B_1).$
Let ${\cal U}'=\{quq: u\in U(B)\cap {\cal P}\}$
and let  $F$ be the subgroup of $U(qBq)$ generated by ${\cal U}'.$
Let ${\bar F}$ be the image of $F$ in $U(qBq)/CU(qBq).$
By 6.6 (3), we write ${\bar F}={\bar F}\cap U_0(qBq)/CU(qBq)\oplus \bar{F_0}
\oplus \bar{F_1},$ where $\bar{F_0}$ is torsion and
$\bar{F_1}$ is free. Note $\kappa({\bar F}_1)\cong{\bar F}_1.$
We may assume that ${\cal U}'={\cal U}_0\bigcup{\cal U}_1$ such that
$\overline{{\cal U}_0}$ generates
${\bar F}\cap U_0(qBq)/CU(qBq)\oplus {\bar F}_0$ and
$\overline{{\cal U}_1}$ generates $\overline{ F_1}.$
Note that modulo unitaries
in $CU(qBq),$ we can always make this assumption and it may cause us
no more than $8\pi$ in the estimation of the exponential length (see 6.9).
Without loss of generality, we may also assume
that $q\in B_1$ and ${\cal U}_0,$ ${\cal U}_1\subset qB_1q.$
Note we also assume that $K_1(B_m)\to K_1(B_{m+1})\to K_1(B)$ is
injective.

Let ${\cal G}_1'$ be a finite subset which contains
${\cal S},$ ${\cal G}_0$ and ${\cal H}$ as well as $q, p_1,...,p_n$ and a finite
generating set of $C_1'.$ It also contains ${\cal U}'.$ Without
loss of generality, we may further assume that ${\cal G}_1'\subset B_1.$
Note we still have that $q$ commutes with all elements in ${\cal
G}_1'.$ It follows from 10.3 (or 9.10) that there is a ${\cal
G}_1'$-$\dt/4$-multiplicative \morp\, $L_1: B\to A$ such that
$$
[L_1]|_{{\cal P}\cup {\cal P}_0}=\alpha^{-1}|_{{\cal P}\cup{\cal P}_0}
\andeqn
$$
$$
\sup_{\tau\in T(B)}\{|\tau\circ L_1(a)-\xi^{-1}(Q(a))(\tau)|\}
<\sigma/2
\,\,\,{\rm for\,\,\, all\,\,\,} a\in {\cal H},
$$
where $Q: A_{s.a}\to Aff(T(A))$ is
the evaluation map.

We assume that $L_1^{\ddag}$ is well defined on  the $\overline{
F}$ ($\subset U(qBq)/CU(qBq)$) (see 6.2). Define ${\bf L}_1(A)\to
{\bf R}_+$ exactly the same way as ${\bf L}$ above. Let ${\cal
F}_1$ be a finite subset of $A.$ Let $\dt_1'>0,$ integer $n_1>0,$
finite subsets ${\cal P}_1\subset {\bf P}(A),$ ${\cal S}_1,
\subset A$
(for $A,$ ${\bf L}_1,$ $\ep/4$ and $ {\cal F}_1$) be as required in
Theorem 8.6.
There are mutually orthogonal projections $q',
p_1',p_2',..,p_{n_1}'$ with $q'\lesssim p_i',$ and with
$[p_i']=[p_1']$ ($i=1,2,...,n_1$) and there is a \SCA\, $C_2'\in
{\cal I}$ with $1_{C_2'}=q'$ and there are unital
$S_1$-$\dt_1'/4$-multiplicative \morp s $h_0': A\to q'Aq'$ and
$h_1: A\to C_2'$ such that $h_0(x)=q'xq'$ and
$$
\|x-(h_0'(x)\oplus h_1'(x)\oplus\cdots\oplus
h_1'(x))\|<\dt_1'/16\,\,\,\rforal\,\,\,x\in S_1,
$$
where $h_1'$ is repeated $n_1$ times. We may assume that
$L_1({\cal S})\subset {\cal S}_1.$ Put $C_1=M_{n_1}(C_2')\subset
(1-q)A(1-q).$ Let ${\cal P}_{01},$ ${\cal G}_{01},$ ${\cal H}_1,$
$\dt_{01},$ and $\sigma_{01}$ also be  as required by 8.6.
Let $\dt_1=\min\{\dt_1', \dt_{01}\}.$ We may assume that
$\dt_1<\dt/2,$ $\sigma_1<\sigma/4$ and ${\cal P}_{01}$
contains at least one minimal projection in each summand of
$C_1$. Furthermore, without loss of generality, we may also
assume that $q'$ commutes with each element in ${\cal H}_1$ and
${\cal S}_1$,
and that $[{\cal P}_1]\supset [L_1]({\cal P}\cup {\cal P}_{0}).$

Without loss of generality, we may assume that each
 $u\in U(A)\cap {\cal P}_1$ has
the form $q'uq'+(1-q')u(1-q'),$ where $q'uq'\in U(q'Aq')$ and
$(1-q')u(1-q')\in U((1-q')C_1(1-q')).$ Put ${\cal V}'=\{q'uq':
u\in U(A)\cap {\cal P}\}.$ Let $F'$ be the subgroup of $U(q'Aq')$
generated by ${\cal V}'.$ By 6.6 (3), we write
$\overline{F'}=\overline{F'}\cap U_0(q'Aq')/CU(q'Aq')\oplus
\overline{ F_0'}\oplus \overline{F_1'},$ where ${\overline{F_0'}}$
is torsion and $\overline{F_1'}$ is free. Note
$\kappa(\overline{F_1'})\cong \overline{F_1'}.$ We may also assume that
$\bar{F}'\supset L_1^{\ddagger}(\bar{F}).$ Since  $A$ is
simple and separable, $K_1(q'Aq')=K_1(A).$ Without loss of generality, we may
also assume that ${\cal V}'={\cal V}_0\bigcup{\cal V}_1,$ where
$\overline{{\cal V}_0}$ generates
$\overline{F'}\cap U_0(q'Aq')/CU(q'Aq')\oplus
\overline{F_0'}$ and $\overline{{\cal
V}_1}$ generates $\overline{F_1'}.$ Note again this assumption may cause us
no more than  $8\pi$ when we estimate the exponential length later
(see 6.9).

Let ${\cal G}_2'$ be a finite subset which contains ${\cal S}_1,$
${\cal G}_{01},$  ${\cal H}_1$ and
$L_1({\cal G}_1')$ as well as $q', p_1',..., p_{n_1}'$ and a
finite generating set of $C_2.$ It also contains ${\cal V}'.$ Note
$q'$ commutes with all elements in ${\cal G}_2'.$ It follows from
10.3 that there is a ${\cal G}_2'$-$\dt_1/2$-multiplicative
\morp\, $\Phi_1': A\to B$ such that
$$
[\Phi_1']|_{{\cal P}_1\cup{\cal P}_{01}}=\alpha|_{{\cal P}_1\cup{\cal P}_{01}}
\andeqn
$$
$$
\sup_{\tau\in T(A)}\{|\tau\circ
\Phi_1'(a)-\xi((Q(a))(\tau)|\}<\sigma_1/4 \,\,\,{\rm for\,\,\,
all\,\,\,} a\in {\cal H}_1\cup L_1({\cal H}).
$$
We may also assume that $(\Phi_1')^{\ddag}$ is well defined on
$\overline{F'}$ and $(\Phi_1'\circ L_1)^{\ddag}$ is well defined
on $\overline{F}.$ Without loss of generality (see 6.2 in
\cite{LnC}), we may assume that the image of $\Phi_1'$ is
contained in $B_n.$

Let $B_n'=qB_nq.$ Since $B$ is simple, it is known and easy to see
that, by choosing a possibly large $n,$ we may assume that the
rank of $q$ at each point is at least $6$ (in $B_n$). We have assumed that
$q\in B_n.$ So $B_n'$ is a corner of $B_n.$ By the construction,
we know that  $[\Phi_1'\circ L_1](q)$ is equivalent to $q.$ By
replacing $\Phi_1'$ by ${\rm ad} w\circ \Phi_1'$ for some unitary
$w$ if necessary, we may assume that
$$
\|\Phi_1'\circ L_1(q)-q\|<\dt/4.
$$
Define $\Lambda(b)=aq[\Phi_1\circ L_1(qbq)]qa$
(for $a=(q\Phi_1\circ L_1(q)q)^{-1/2}$ and ) for $b\in qBq.$
Note that
$$
\|\Lambda-\Phi_1'\circ L_1|_{qBq}\|<\dt/2.
$$

Write $B_n=\oplus_{j=1}^mB_n(j),$ where each $B_n(j)$
has the form $C^{(j)}$ as described in 7.1.
According to this direct sum decomposition, we may write
$q=q_1\oplus q_2\oplus \cdots \oplus q_l$ with $0\le l\le m$
and $q_j\not=0,$ $1\le j\le l.$
Choose an integer $N_1>0$ such that $N_1[q_j]\ge 3 [1_{B_n(j)}]$
for $j\le l.$ Note that we may assume
that $q_j$ has rank at least $6.$
By  applying an inner automorphism, we may assume that
$\oplus_{j=1}^lB_n(j)$ is a hereditary \SCA\, of $M_{N_1}(B_n').$
Since $F_1$ is finitely generated, with sufficiently large $n$,
we obtain (see 6.2) a \hm\,
$j: {\bar F}_1\to U(qB_n'q)/CU(qB_n'q)$
such that $\psi_n^{\ddag}\circ j={\rm id}_{{\bar F}_1}.$
Then (since $K_1(B_1)\to K_1(B_n)\to K_1(B)$ is injective),
$$
\kappa_1\circ \psi_n^{\ddag}\circ (\Phi_1'\circ
L_1)^{\ddag}|_{\overline{F_1}}= \kappa_1\circ
(\psi_n)^{\ddag}\circ j=(\kappa_1)_{\overline{F_1}},
$$
where $\kappa_1: U(qBq)/CU(qBq)\to K_1(qBq)$ is the quotient map.
Note that $K_1(qBq)=K_1(B).$ Let $\Delta_1$ be $\dt(\ep/4)$ as
described in 7.5. We may assume that $\Delta_1<\sigma_1/4.$ To
simplify notation, without loss of generality, we may assume that
$\psi_n(q)=q.$ By the assumption on $B,$  we may write that
$\psi_n|_{B_n'}=(\psi_n)_0\oplus (\psi_n)_1,$ where

(1) $\tau((\psi_n)_0(1_{B_n'}))<\Delta_1/2(N_1+1)^2$ for all $\tau\in T(B)$
and

(2) $(\psi_n)_0$ is homotopically trivial (but nonzero).

(see 10.1)\\
 It follows from 7.5 that there is a \hm\, $h: B_n'\to
e_0Be_0$ such that

(i) $[h]=[(\psi_n)_0]$ in $KL(B_n', B)$ and

(ii) $(\psi_n^{\ddag}\circ j({\bar w}))^{-1}
             (h\oplus (\psi_n)_1)^{\ddag}(\Lambda^{\ddag}({\bar w}))
=\overline{g_w},$ where $g_w\in U_0(qBq)$ and $cel(g_w)<\ep/4$ (in
$U(qBq)$) for all $w\in {\cal U}_1.$

Define (we have assumed that $B_n\subset M_{N_1}(B_n')$)
$$
h'=((h\oplus (\psi_n)_1)\otimes
 {\rm id}_{M_{N_1}})|_{\oplus_{j=1}^lB_n(j )}.
$$
and define $\Psi'=h'\oplus (\psi_n)|_{\oplus_{j=l+1}^m B_n(j)}.$
Let $\Phi_1=\Psi'\circ \Phi_1'.$
It is clear that (since $\Delta_1<\sigma_1/4$)
$$
[\Phi_1]|_{{\cal P}_1\cup {\cal P}_{01}}=[\Phi_1']|_{{\cal P}_1\cup {\cal P}_{01}}
\andeqn
|\tau\circ \Phi_1(a)-\tau\circ \Phi_1'(a)|<\sigma_1/2
$$
for all $a\in A_{s.a}$ and $\tau\in T(B).$ For all $w\in {\cal
U}_1,$ we have, by (ii) above,
$$
cel(w^*(\Phi_1\circ L_1(w)))<8\pi+\ep/4 \,\,\,\,\,\,{\rm in}\,\,\,
U(qBq).
$$
For $w\in {\cal U}_0,$ by 6.8, 6.10 and 6.9, we also
have
$$
cel(w^*(\Phi_1\circ L_1(w)))<2cel(w)+\pi/64\,\,\,\,\, {\rm
(or\,\,\,} <8\pi+2cel(w^k)/k+\pi/16{\rm )}\,\,\,{\rm in}\,\,\,
U(qBq)
$$
(depending if $[w]=0$ or $[w]$ has order $k$ in $K_1(B)$).
(Recall the definition of $h_0$ earlier in this proof
 in the next estimate.)
Therefore (even after we add $8\pi$ for the decomposition assumption
of $\overline{F}$)
$$
cel({\rm id}_B(h_0(u))^{-1}(\Phi_1\circ L_1(h_0(u)))<{\bf L}(u)\,\,\,
{\rm in}\,\,\,U(qBq)
$$
for all $u\in U(B)\cap {\cal P}_1.$
Since we also have
$$
[{\rm id}_B]|_{{\cal P}\cup {\cal P}_0}=[\Phi_1\circ L_1]|_{{\cal P}\cup {\cal P}_0}
\andeqn \sup_{\tau\in T(B)}\{|\tau(a)-\tau(\Phi_1\circ L_1(a))|\}<\sigma
$$
for all $a\in {\cal H},$
by 8.6, we obtain a unitary $W\in U(B)$ such that
$$
{\rm ad} W\circ \Phi_1\circ L_1\approx_{\ep/2} {\rm id}_B
\,\,\,\,{\rm on }\,\,\, {\cal F}.
$$
Replacing $\Phi_1$ by ${\rm ad}\,W\circ \Phi_1,$ we may assume
that
$$
\Phi_1\circ L_1\approx_{\ep/2} {\rm id}_B
\,\,\,\,{\rm on}\,\,\, {\cal F}.
$$

Now let ${\cal F}_2\subset B.$ We may assume that
${\cal F}_2\subset B_{m_1'}$ ($m_1'>n$). Let $\dt_2'>0,$ integer $n_2>0,$
finite subsets ${\cal P}_2\subset {\bf P}(B),$
${\cal S}_2, \subset B,$
be as required
by Theorem 8.6 (for $B,$ ${\bf L},$ $\ep/16$ and ${\cal F}_2$).
There are mutually orthogonal projections $q'',
p_1'',p_2'',..,p_n''$ with $q''\lesssim p_i'',$ and with
$[p_i'']=[p_1'']$ ($i=1,2,...,n$) and there is a \SCA\, $C_3'\in
{\cal I}$ with $1_{C_3'}=q'',$ and there are unital ${\cal
S}_2$-$\dt_2'/4$-multiplicative \morp s $h_0'': B\to q''Bq''$ and
$h_1'': B\to C_3'$ such that $h_0''(x)=q''xq''$ and
$$
\|x-(h_0''(x)\oplus h_1''(x)\oplus\cdots\oplus
h_1''(x))\|<\dt_2'/4
$$
for all $x\in {\cal S}_2,$ where $h_1''$ is repeated $n_2$ times.
We assume that ${\cal S}_2\supset \Phi_1({\cal S}_1).$ Put
$C_2=M_{n_2}(C_3')\subset (1-q'')B(1-q'').$ Let ${\cal P}_{02},$
${\cal G}_{02},$ ${\cal H}_2,$ $\dt_{02}>0$ and $\sigma_2>0$ also be as
required by 8.6. Let $\dt_2=\min\{\dt_2',\dt_{02}\}.$ We may
assume that $\sigma_2<\sigma_1/4$ and $\dt_2<\dt_1/4.$
We may also assume that ${\cal P}_{02}$ contains at least one
minimal projection of each summand of $C_2$ and $[{\cal
P}_2]\supset [\Phi_1({\cal P}_1\cup {\cal P}_{01})].$ Furthermore,
we may assume that each $u\in U(B)\cap{\cal P}_2$ has the form
$q''uq''\oplus (1-q'')u(1-q''),$ where $q''uq''\in U(q''Bq'')$ and
$(1-q'')u(1-q'')\in U(C_2).$ Put ${\cal W}=\{q''uq'': u\in
U(B)\cap{\cal P}_2\}.$ Let $F''$ be the subgroup generated by
${\cal W}.$ Write $\overline{F''}=\overline{F''}\cap
U_0(q''Bq'')/CU(q''Bq'') \oplus \overline{F_0''}\oplus
\overline{F_1''},$ where $\overline{F_0''}$ is torsion and
$\overline{F_1''}$ is free. We may also assume that
$\overline{F}''\supset \Phi_1^{\ddagger}(\overline{F}').$ We may further
assume, without loss of generality,
 that ${\cal W}={\cal W}_0\bigcup {\cal W}_1,$ where
${\cal W}_0$ generates $\overline{F''}\cap U_0(q''Bq'')/CU(q''Bq'')
\oplus \overline{F_0''}$
and ${\cal W}_1$ generates $\overline{F_1''}.$

Let ${\cal G}_3'$ be a finite subset which contains ${\cal S}_2,$
${\cal G}_{02},$
$q'', p_1'',...,p_{n_1}'',$ ${\cal H}_2,$ $\Phi_1({\cal G}_2''),$
a generating set of $C_2$ and
${\cal W}.$
Without loss of generality, we may assume (see 6.2 in \cite{LnB}) that
$\Phi_1(A)\subset B_{m}$ ($m>m_1>n$).
By choosing a larger $m,$ we may assume that there
is a \morp\, $J: B\to B_m$ such that
$$
\|J(a)-{\rm id}(a)\|<\dt_2/8\,\,\,{\rm for \,\,\, all}\,\,\,
a\in {\cal G}_3'.
$$
There is a projection ${\tilde q'}\in B_m$ such that
$$
\|\Phi_1(q')-{\tilde q'}\|<\dt_2/2.
$$
We may write $B_m=\oplus_{j=1}^s B_m(j).$ As above, by choosing a
possibly larger $m,$ we may also assume that ${\tilde q'}$ has rank at
least $6$ at each point. We write that ${\tilde q'}=q_1'\oplus
q_2'\oplus\cdots q_l'$ according to the direct sum decomposition
($q_j'\not=0$ for $1\le j\le l$ and $l\le s$). Suppose that
$N_2>0$ is an integer such that $N_2[q_j']> 3[1_{B_m(j)}]$ for
$1\le j\le l.$ Set $B_m'={\tilde q'}B_m{\tilde q'}.$
Note $\Phi_1^{\ddagger}$ is injective on $\bar{F}_1'.$
We
may further assume that ${\cal G}_3'$ contains $q_1', q_2',...,q_l'$
and  a generating set of $B_m'$ and $B_m.$

Now let $L_2': B\to A$ be a ${\cal
G}_3'$-$\dt_2/16(N_2+1)^2$-multiplicative \morp\, (10.3) such that
$$
[L_2']|_{{\cal P}_2}=\alpha^{-1}|_{{\cal P}_2\cup{\cal P}_{02}}\andeqn
$$
$$
\sup_{\tau\in T(A)}\{|\tau\circ L_2'(a)-\xi^{-1}(Q(a))(\tau)
|\}<\sigma_2/4\,\,\, {\rm for\,\,\, all}\,\,\, a\in {\cal H}_2\cup
\Phi_1({\cal H}_1) .
$$
We may assume that $(L_2')^{\ddagger}$ is well defined on
$\overline{F''}.$ Suppose that $e\in A$ is a projection such that
$$
\|L_2'\circ \Phi_1(q')-e\|<\dt_2/4.
$$
Since $q'\in {\cal P}_1,$ $[e]=[q']$ in $K_0(A).$
Thus, by replacing $L_2'$ by  ${\rm ad} u'\circ L_2'$ for some
unitary $u'\in A,$ we may assume that
$e=q'.$ Without loss of generality, to simplify notation,
we may assume $L_2'\circ \Phi_1(q')=q'.$
Note $\overline{F'_1}$ is free and $(\phi_{m,M})_{*1}$ ($M>m$) is injective.
We compute that
$$
\alpha^{-1}\circ (\psi_m)_{*1}\circ \kappa_1'\circ
(\Phi_1)^{\ddag}(g)=\kappa_1(g)
$$
for all $g\in \overline{F_1'},$ where $\kappa_1':
U(B_m')/CU(B_m')
\to K_1(B_m')$ and
$\kappa_1: U({\tilde q'}B{\tilde q'})/CU({\tilde q'}B{\tilde q'})\to
K_1({\tilde q'}B{\tilde q'})$ are
the quotient maps. Note that $K_1({\tilde q'}B{\tilde q'})=
K_1(B).$
With $\Phi_1^{\ddagger}$ playing the role of $L$ and
$\alpha^{-1}\circ (\psi_{m})_{*1}$
playing the role of $\alpha$ in 7.3, by applying 7.3, we obtain
 a \hm\, $\beta: U(B_m')/CU(B_m')\to U(q'Aq')/CU(q'Aq')$
with $\beta(U_0(B_m')/CU(B_m'))\subset U_0(q'Aq')/CU(q'Aq')$
such that
$$
\beta\circ (\Phi_1^{\ddagger})({\bar w})={\bar w}
$$
for all ${\bar w}\in \overline{F_1'}.$ Let $\Delta_2'=\dt(\ep/16)$
be as described in 7.4. It follows from the assumption on $B$ that
there is $M>m$ such that $\phi_{m,M}=\phi_{m,M}^{(0)}\oplus
\phi_{n,M}^{(1)}: B_m\to B_M$ such that $ \phi_{m, M}^{(0)}$ is
(nonzero) homotopically trivial and $\tau(\psi_{M}\circ
\phi_{n,M}^{(0)}(1_{B_m'}))<\Delta_2'/4(N_2+1)^2$ for all $\tau\in
T(B).$ To simplify notation, without loss of generality, we may
also assume that $e_0'=L_2'\circ \psi_M\circ
\phi_{m,M}^{(0)}(1_{B_m'})$ and $e_1'=L_2'\circ \psi_M\circ
\phi_{m,M}^{(0)}(1_{B_m'})$ are mutually orthogonal projections
(see 8.2 (ii) ). Note that $\psi_M\circ \phi_{m,M}=\psi_m.$ It
follows from 7.4 (by the choice of $\Delta_2'$ and with $\beta$
playing the role of $\alpha,$ $L_2'\circ \psi_M\circ
\phi_{m,M}^{(0)}$ playing the role of $\phi_0$ and $L_2'\circ
\psi_M\circ \phi_{m,M}^{(1)}$ playing the role of $\phi_1$ in
7.4), we obtain a \hm\, $\Phi': B_m'\to e_0'Ae_0'$ such that

(i') $\Phi'$ is homotopically trivial, $\Phi'_{*0}=
[L_2']\circ (\psi_M\circ \phi_{m,M}^{(0)})_{*0}|_{K_0(B_m')}$ and

(ii') $(\beta( \Phi_1^{\ddagger} ({\bar w}))^{-1})(\Phi' \oplus
(L_2'\circ \psi_M\circ \phi_{m,M}^{(1)}))^{\ddag}(
\Phi_1^{\ddagger}({\bar w}))=\overline{g_{ w}},$ where $g_w\in
U_0(q'Aq')$ and $cel(g_w)<\ep/4$ (in $U(q'Aq')$) for all $w\in
{\cal V}_1.$ As in the construction of $\Psi',$ we obtain a \hm\,
${\tilde \Phi'}: B_m \to q'Aq'$ such that ${\tilde \Phi'}$ is
homotopically trivial, ${\tilde \Phi'}_{*0}=[L_2']\circ
(\psi_M\circ \phi_{m,M})_{*0}$ and ${\tilde \Phi'}|_{B_m'}=\Phi'.$

Define $L_2=({\tilde \Phi'}\oplus L_2'\circ \psi_M\circ
\phi_{n,M}^{(1)}) \circ J.$ It is clear that
$$
[L_2]|_{{\cal P}_2\cup {\cal P}_{02}}=[L_2']|_{{\cal P}_2\cup
{\cal P}_{02}} =\alpha^{-1}|_{{\cal P}_2\cup {\cal P}_{02}}.
$$
Given the choice of $\Delta_2,$ we also have
$$
|\tau\circ L_2(a)-\tau(L_2'(a))|<\sigma_2/4
$$
for all $a\in A_{s.a}$ and $\tau\in T(A).$  In particular,
$$
\sup_{\tau\in T(A)}\{|\tau\circ L_2\circ \Phi_1(a)-\tau(a)|\}<\sigma_1/2
$$
for all $a\in {\cal H}_1.$
Since $\beta\circ \Phi_1^{\ddagger}({\bar w})={\bar w}$ for
all $w\in {\cal V}_1,$  by (ii') and 6.9, we have
$$
cel({\rm id}_A(h_0'(w^*))L_2(\Phi_1(h_0'(w)))
<8\pi +cel(g_w)+\ep/4< 8\pi+\ep/2\,\,\,\,{\rm in}\,\,\, U(q'Aq')
$$
for all $w\in {\cal V}_1.$ We also have, by 6.8, 6.9 and 6.10,
$$
cel({\rm id}_A(h_0'(w^*))L_2(\Phi_1(h_0'(w))))<2cel(w)+\pi/16
\,\,\,{\rm ( or}\,\,\, <
8\pi+2cel(w^k)/k+\pi/16 {\rm )}\,\,\,{\rm in}\,\,\,
U(q'Aq')
$$
for all $w\in {\cal V}_0$ (depends on if $[w]=0$ in $K_1(A)$ or
$[w]$ has order $k$).
It follows that (by adding $8\pi$)
$$
cel({\rm id}(h_0'(u^*))L_2(\Phi_1((h_0'(u))))<{\bf L}(u)
$$
for all $u\in U(A)\cap {\cal P}_2$ (in $q'Aq'$).

By applying 8.6, we obtain a unitary
$Z\in U(A)$ such that
$$
{\rm ad}Z\circ L_2\circ \Phi_1(a)\approx_{\ep/16} {\rm id}_A\,\,\,
{\rm on}\,\,\, {\cal F}_1.
$$

Therefore, by replacing $L_2$ by ${\rm ad} Z\circ L_2,$
 we obtain  the following ``approximate intertwining" diagram:
$$
\begin{array}{ccc}
B & \to_{\rm id_B} & B\\
\downarrow_{L_1} & \nearrow_{\Phi_1} & \downarrow_{L_2}\\
A& \to_{\rm id_A} & A.
\end{array}
$$
Since this process continues, we see that $L_1$ is
recursively ${\cal F}$-invertible (and $\Phi_1$ is
recursively ${\cal F}_1$-invertible----see 3.6 in \cite{LnC}).
It follows from an argument of Elliott
(see Theorem 3.6 in \cite{LnC}, for example)
that $A$ is isomorphic to $B.$
\QED

{\bf 10.5. Remark}
If $K_1(A)$ and $K_1(B)$ are torsion groups, then one can use the
``uniqueness theorem" 8.7. Since we do not need to control
exponential length in this case, section 7 is not needed.
Furthermore, we do not need to assume $B$
is AH. Consequently, we  do not need to assume the condition
on the $\partial_e(S_u(K_0(A))$ either. The whole proof is a
much shorter.

Since simple AH-algebra with very slow dimension growth have tracial
rank one or zero. We have the following.

\vspace{0.1in}

{\bf 10.6. Theorem} (\cite{EGL} and \cite{G2}) {\it Let $A$ and
$B$ be two unital simple AH-algebras with very slow dimension
growth and with torsion $K_1(A).$ Then $A$ is isomorphic to $B$ if
and only if
$$(K_0(A), K_0(A)_+, [1_A], K_1(A), T(A))\cong
(K_0(B), K_0(B)_+, [1_B], K_1(B), T(B))
$$
in the sense of 10.2.
}

{\bf 10.7} Let  $C$ be a stably finite non-unital \CA\, with
an approximate identity consisting of projections $\{e_n\}.$
Denote by $T(C)$ the set of traces $\tau$ on $A$ such that
$\sup_n\tau(e_n)=1.$ It is clear that $T(C)$ is the tracial
state space. Note that each tracial state extends to a tracial
state on ${\tilde C}.$
Therefore  $T({\tilde C})$ is the set of convex
combinations of $\tau\in T(C)$ and the tracial state which vanishes on $C.$
We also denote by $S_{u'}(K_0(C))$ the set of those order preserving
\hm s from $K_0(C)$ to ${\mathbb R}$ such that $\sup_ns([e_n])=1.$
Then each element in $S_u(K_0({\tilde C}))$ is the convex combination of $s\in
S_{u'}(K_0(C))$ and the state which vanishes on $j_*(K_0(C)),$
where $j: C\to {\tilde C}$ is the embedding.

\vspace{0.1in}

{\bf 10.8 Lemma} {\it Let $A$ be a unital separable simple \CA\,
with $TR(A)\le 1.$ Then there is a \CA\, $C=\lim_{n\to\infty}(C_n,
\phi_n),$ where $C_n\in {\cal I},$ satisfying the following:

{\rm (i)} each $C_n$ is a \SCA\, of $A$ and
$\{\phi_{n,\infty}(1_{C_n})\}$ forms an approximate identity for $C,$

{\rm (ii)} there is a sequence of \morp s $L_n: A\to C$ such that
$$
\lim_{n\to\infty}\|L_n(ab)-L_n(a)L_n(b)\|=0,\,\,\, a, b\in A,
$$

{\rm (iii)} there is an affine continuous (face-preserving)
isomorphism $r^{\sharp}: T(A)\to T(C)$ such that
 $$
r^{\sharp}(\tau)(\phi_{n, \infty}(b))=\lim_{k\to\infty}
\tau(\phi_{n,k}(b)) \,\,\,for\,\,\,all\,\,\, b\in C_n\,\,\,and
\,\,\,  \tau\in T(A),
$$

{\rm (iv)} there is an affine  continuous (face-preserving)
isomorphism $r_{\sharp}: S_{u'}(K_0(C))\to S_u(K_0(A))$ such that
$$
r_{\sharp}(s)([p])=\lim_{n\to\infty}\tau_s(\phi_{n,\infty}(L_n(p)))\,\,\,for\,\,\,all\,\,\,
s\in S_{u'}(K_0(C))\,\,\, and \,\,\, projection\,\,\,p\in A,
$$
where $\tau_s$ is the trace which induces $s.$ }

{\it Proof:} 
\,Let ${\cal F}_1\subset {\cal F}_2\subset\cdots {\cal F}_n\cdots$
be a sequence of finite subsets of $A$ such that $\cup_n {\cal
F}_n$ is dense in $A.$ Since $TR(A)\le 1,$ there is a
\SCA\, $C_1\subset A$ with
$C_1\in {\cal I}$ and $1_{C_1}=p_1$ such that

(1') $\|ap_1-p_1a\|<1/2$ for all $a\in {\cal F}_1$

(2') ${\rm dist}(p_1ap_1, C_1)<1/2$ for all $a\in {\cal F}_1.$

(3') $\tau(1-p_1)<1/4$ for all $\tau\in T(A).$

Let $1>\eta_1>0.$
By (3') and 4.7,
there is a projection $e_{(1,1)}\le p_1$ such that
$e_{(1,1)}$ is equivalent to $1-p_1.$
Since $\tau(p_1-e_{(1,1)})>1/2>\tau(1-p_1)$ for all $\tau\in T(A),$ by
4.7 again, we obtain
mutually orthogonal projections $e_{(1,1)},e_{(1,2)}$ such that
$e_{(1,i)}\le p_1,$ $[e_{(1,1)}]=[e_{(1,2)}]\ge [1-p_1].$ There are
$x_{(1,1)},x_{(1,2)}\in A$ such that $x_{(1,i)}^*x_{(1,i)}\ge 1-p_1$ and
$x_{(1,i)}x_{(1,i)}^*=e_{(1,i)}.$ Let ${\cal G}_1'$ be a finite set of
generators of $C_1$ and ${\cal G}_2={\cal F}_2\cup {\cal
G}_1'\cup\{x_{(1,i)},x_{(1,i)}^*, e_{(1,i)}: 1\le i\le 2\}.$ There is a
\SCA\, $C_2\subset A$ with $C_2\in {\cal I}$ and $1_{C_2}=p_2$
such that

(1'') $\|ap_2-p_2a\|<\eta_1/4$ for all $a\in {\cal G}_2,$

(2'') ${\rm dist}(p_2ap_2, C_2)<\eta_1/4$ for all $a\in {\cal
G}_2$ and

(3'') $\tau(1-p_2)<1/8$ for all $\tau\in T(A).$

By 2.1 (iii), with sufficiently
small $\eta_1,$ there is a \hm\, $\phi_1: C_1\to C_2$ such that
$$
\|\phi_1(b)-p_2bp_2\|<1/4
\,\,\,\rforal \,\,\,b\in {\cal F}_1\cup {\cal G}_1'.
$$
Put $q_2=\phi_1(1_{C_1}).$
With sufficiently small $\eta_1,$
since $x_{(1,i)}\in {\cal G}_2,$ we may also assume that
$2[p_2-q_2]\le [q_2]$ in $K_0(C_2).$

We continue in this fashion. Suppose that $C_n\subset A$ is a
unital \SCA\, which is in ${\cal I}$ has been constructed. If
$\tau(1-p_n)<1/2^{n+1}$ for all $\tau\in T(A),$  there  are
partial isometries $x_{(n,i)}\in A$ such that
$x_{(n,i)}^*x_{(n,i)}=e_{(n,i)}\le p_n,$ $e_{(n,i)}e_{(n,j)}=0$ if $i\not=j$ and
$[e_{(n,i)}]=[e_{(n,1)}]\ge [1-p_n],$ $1\le i\le 2^n.$
 Let ${\cal G}_n'$
be a finite set which contains a set of generators of $C_n$ and
$\phi_{i,n}(p_i),$ $i=1,2,...,n-1,$ where
$\phi_{i,n}=\phi_{n-1}\circ \phi_{n-2}\circ \cdots \circ \phi_i.$ (Note
that $C_i\subset A.$)  Let ${\cal G}_{n+1}={\cal
F}_{n+1}\cup {\cal G}_n'\cup\{e_{(n,i)}, x_{(n,i)}, x_{(n,i)}^*: 1\le i\le
2^n\}. $ Let $1>\eta_{n+1}>0$ be a positive number to be determined
(but it depends only on $C_n$ and ${\cal F}_n\cup {\cal G}_n'$).
Since $A$ has tracial topological rank one, there exist a \SCA\,
$C_{n+1}\subset A$ with $C_{n+1}\in {\cal I}$ and a projection $p_{n+1}$ with
$1_{C_{n+1}}=p_{n+1}$ such that

(1) $\|ap_{n+1}-p_{n+1}a\|<\eta_{n+1}/2^{n+1}$ for all $a\in {\cal G}_{n+1},$

(2) ${\rm dist}(p_{n+1}ap_{n+1}, C_{n+1})<\eta_{n+1}/2^{n+1}$
for all $a\in {\cal G}_{n+1}$ and

(3) $\tau(1-p_{n+1})<1/2^{n+2}$ for all $\tau\in T(A).$

We choose $\eta_{n+1}$ so small that there exists a \hm s
$\phi_n: C_n\to C_{n+1}$ such that
$$
\|\phi_n(b)-p_{n+1}bp_{n+1}\|<1/2^{n+1}
\rforal b\in {\cal F}_{n}\cup {\cal G}_n'.
$$
 Put $q_{n+1}=\phi_n(p_n).$
Since $x_{(n,i)}\in {\cal G}_{n+1},$ we may further assume that

(4) $2^n[p_{n+1}-q_{n+1}]\le [q_{n+1}]$ in $K_0(C_{n+1}).$

Set
$C=\lim_{n\to\infty}(C_n, \phi_n).$ Since $C_n$ are nuclear, there
is a \morp\, $L_n': A\to C_n$ such that
$$
\lim_{n\to\infty}\|L_n'(a)-p_nap_n\|=0
$$
for all $a\in A.$
Note that, by (1),
$$
\lim_{n\to\infty}\|L_n'(ab)-L_n'(a)L_n'(b)\|=0\rforal\,\,\, a, b\in A.
$$
Define $L_n=\phi_{n,\infty}\circ L_n'.$
It is clear that  $L_n$ satisfies (ii).
Put $\phi_{n,n+1}=\phi_n$ and for $k>n+1,$
$\phi_{n,k}=\phi_{k-1}\circ \cdots \circ \phi_n.$
Define $r^{\sharp}: T(A)\to T(C)$ as follows. For each $b\in C_n,$
define
$$
r^{\sharp}(\tau)(\phi_{n,
\infty}(b))=\lim_{k\to\infty}\tau(\phi_{n,k}(b)).
$$
Note that $\phi_{n,k}(b)\in C_k$ and $C_k\subset A.$
By (1)-(3) above and the definition of $\phi_n,$ the
right side converges. To see $r^{\sharp}$ is well defined, we let $c\in
C_m$ so that $\phi_{m,\infty}(c)=\phi_{n,\infty}(b).$ Then, for
any $\ep>0,$ there exists $N>\max\{n,m\}$ such that
$$
\|\phi_{n,k}(b)-\phi_{m,k}(c)\|<\ep
$$
for all $k\ge N.$ It follows that ($C_k\subset A$)
$$
|\tau(\phi_{n,k}(b))-\tau(\phi_{m,k}(c))|<\ep
$$
for all $\tau\in T(A)$ and $k\ge N.$ It follows that $r^{\sharp}$
is well defined.
It is then easy to see that $r^{\sharp}$ is an affine continuous
map. Define $r^{\sharp\,-1}: T(C)\to T(A)$ by
$$
r^{\sharp\,-1}(t)(a)=
\lim_{n\to\infty}t(L_n(a))=
\lim_{n\to\infty}t(\phi_{n,\infty}(L_n'(a)))\,\,\,\rforal
t\in T(C) \andeqn a\in A.
$$
It is well defined by (1)-(4) above, and, by (ii),
$r^{\sharp\,-1}(t)$ is a trace on $A.$  It is then clear that
$r^{\sharp\,-1}$ is an affine continuous map. It should be
noted that even if $a\in C_m$ (for $m<n$), $L_n'(a)\in C_n.$ Now
let $\tau\in T(A)$ and $a\in A.$
To show that $(r^{\sharp\,-1}\circ r^{\sharp})(\tau)(a)=\tau(a),$
we note
that
$$
(r^{\sharp\,-1}\circ r^{\sharp})(\tau)(a)=
\lim_{n\to\infty}r^{\sharp}(\tau)(\phi_{n,\infty}(L_n'(a)))
=\lim_{n\to\infty}(\lim_{k\to\infty}\tau(\phi_{n,k}(L_n'(a))))
$$
for all $a\in A$ and $\tau\in T(A).$
Let $\ep>0.$
Without loss of generality, we may assume that $a\in {\cal F}_n$
for some integer $n>0.$
Moreover, with sufficiently large $n,$ we may assume that
$1/2^n<\ep/4$
and
$$
\|L_n'(a)-p_nap_n\|<\ep/4.
$$
One estimates (with $k>n$)
$$
\|\phi_{n,k}(L_n'(a))-p_kp_{k-1}\cdots p_{n+1}p_nap_np_{n+1}\cdots p_{k-1}p_k\|
<\sum_{j=1}^{k-n}1/2^{n+j}+\ep/4<\ep/2.
$$
By (3), one has
$$
|\tau(p_kp_{k-1}\cdots p_{n+1}p_nap_np_{n+1}\cdots p_{k-1}p_k)-
\tau(a)|
<\ep/4
\rforal \,\,\,\tau\in T(A).
$$
It follows that
$$
|\tau(\phi_{n,k}(L_n'(a)))-\tau(a)|<\ep\rforal \tau\in T(A)
$$
if $k>n.$
Therefore

$
(5)\,\,\,\,\,\,\, \tau(a)=\lim_{n\to\infty}(\lim_{k\to\infty}\tau(\phi_{n,k}(L_n'(a))))
$
for all $a\in A$ and $\tau\in T(A).$

This also proves $r^{\sharp-1}\circ r^{\sharp}(\tau)(a)=\tau(a)$
for all $a\in A$ and $\tau\in T(A).$
Therefore
$r^{\sharp\,-1}\circ r^{\sharp}={\rm id}_{T(A)}.$

Suppose that
$t\in T(C)$ and $b\in C_n.$ Then
$$
r^{\sharp}\circ r^{\sharp\,-1}(t)(\phi_{n,
\infty}(b))=\lim_{k\to\infty} r^{\sharp\,-1}(t)(\phi_{n,k}(b))
=\lim_{k\to\infty}(\lim_{m\to\infty}t(\phi_{m,\infty}(L_m'(\phi_{n,k}(b)))).
$$
Fix $\ep>0.$ Choose $k>n$ such that $1/2^k<\ep/16.$ We may assume
that $\|b\|\le 1.$ For any $m>k,$ put $r_j=\phi_{j,m}(p_j),$
$j=k,...,m-1.$ Since $\phi_j(p_j)\le p_{j+1},$ $r_j\le r_{j+1}.$
By choosing a larger $k,$ applying (1) and (2) above, we may
assume that there is $c_1\in A$
such that (we view $\phi_{n,k}(b)\in C_k\subset  A$)
$$
r_jc_1=c_1r_j, k+1\le j\le m-1,\,\,\,\andeqn \|c_1-\phi_{n,k}(b)\|<\ep/8.
$$
We also have
$$
\|\phi_{n,m}(b)-p_mr_{m-1}\cdots r_k\phi_{n,k}(b)r_k\cdots r_{m-1}p_m\|
<\sum_{j=1}^{m-k}2^{k+j}<\ep/8.
$$
Put $c_2=p_{m}r_{m-1}\cdots r_kc_1.$ It then follows that
$$
c_3=L_{m}'(c_1)- c_2\le 2(p_{m}-r_k).
$$
Since each $C_j$ has stable rank one,
by (4),
there are $y_i\in C_m$ such that
$y_i^*y_i=p_{m}-r_k$ and $y_iy_i^*$ ($1\le i\le 2^m$)  are mutually orthogonal.
Let $z_i=\phi_{m, \infty}(y_i),$ $i=1,2,....,2^m.$
Then $z_i^*z_i=\phi_{m,\infty}(p_m-r_k)$ and
$z_iz_i^*$ ($1\le i\le 2^m$) are mutually orthogonal.
It follows that
$$
t(\phi_{m, \infty}(c_3))\le 2(1/2^m)<\ep/8
$$
for all $t\in T(C).$
On the other hand,
$$
\|[\phi_{m,\infty}(L_m'(\phi_{n,k}(b)))-\phi_{n,\infty}(b)]-
\phi_{m,\infty}(c_3)\|
$$
$$
\le \|[L_m'(\phi_{n,k}(b))-\phi_{n,m}(b)]-c_3\|\le
\ep/8+\ep/8+\|(L_m'(c_1)-c_2)-c_3\|=\ep/4
$$
It follows that
$$
|t(\phi_{m,\infty}(L_m'(\phi_{n,k}(b))))-t(\phi_{n,\infty}(b))|<\ep
$$
for all $t\in T(C)$ if $m>k.$
Thus

$
(6)\,\,\,\,\,\, t(\phi_{n,\infty}(b))=
\lim_{k\to\infty}(\lim_{m\to\infty}t(\phi_{m,\infty}
(L_m'(\phi_{n,k}(b)))))
$
for all $b\in C_n$ and $\tau\in T(C).$

It follows that $r^{\sharp}\circ r^{\sharp-1}={\rm id}_{T(C)}.$
Thus we have shown that $r^{\sharp}$ is an affine continuous surjective
map with an affine continuous inverse $r^{\sharp\,-1}.$ To see
$r^{\sharp}$  is face-preserving, let $\tau\in T(A),$ $t_1, t_2\in
T(C)$ and $0\le a\le 1$ for which
$$
r^{\sharp}(\tau)=at_1+(1-a)t_2.
$$
Let $\tau_1, \tau_2\in T(A)$ such that $r^{\sharp}(\tau_i)=t_i,$
$i=1,2.$  Then, since $r^{\sharp\,-1}$ is the inverse of
$r^{\sharp},$ we see that
$$
\tau=a\tau_1+(1-a)\tau_2.
$$

Fix a projection $p\in A$ and $s\in S_{u'}(K_0(C)).$
Here we will use the notation in 10.7 and 10.2. Then one
obtains a sequence of projections $e_n\in C_n$ such that
$$
\lim_{n\to\infty}\|p_npp_n-e_n\|=0.
$$
Or equivalently
$$
\lim_{n\to\infty}\|L_n'(p)-e_n\|=0.
$$
By (3) and (4)  above, $\{t(e_n)\}$ converges uniformly on $T(C).$
If $p\in M_K(A),$ then one can replace $C$ by $M_K(C)$ and $p_n$
by ${\rm diag}(p_n,...,p_n).$ Since $C$ is an inductive limit of
\CA s in ${\cal I},$  there exists $\sigma_s\in
T(C)$ such that $s([e])=\sigma_s(e)$ for any projection $e\in M_K(C).$
It follows that the following
$$
r_{\sharp}(s)([p])
=\lim_{n\to\infty}\sigma_s(\phi_{n,\infty}(L_n'(p)))=
\lim_{n\to\infty}s([\phi_{n,\infty}(e_n)])
$$
is independent of the choices of $\tau_s$ and
is well defined map from $S_{u'}(K_0(C))$ to $K_u(K_0(A)).$
 (Here we extend $L_n'$ and $\phi_{n,\infty}$ to
$M_K(A)$ and $M_K(C)$ in the obvious way.)  It is clear that
$r_{\sharp}$ is affine and continuous.
Let $e\in M_K(C_n)$ be a projection and $t\in S_u(K_0(A)).$
Since $A$ is a simple \CA\, with
$TR(A)\le 1,$
(by 10.2), there exists a $\tau_t\in T(A)$
such that $\tau_t$ induces $t.$

Define
$$
r_{\sharp}'(t)([\phi_{n,\infty}(e)])=
\lim_{k\to\infty}\tau_t(\phi_{n,k}(e))=\lim_{k\to\infty}t([\phi_{n,k}(e)]),
$$
where we view $C_n$ as a \SCA\, of $A.$ Then, by (1)-(4)
above, $r_{\sharp}'$ is a well defined affine continuous map. Now
let $p\in A$ be a projection and $t\in S_u(K_0(A)).$ By 10.2,
$t$ is induced by a trace $\tau_t\in T(A).$
One has, by (5),
$$
r_{\sharp}(r_{\sharp}'(t))([p])
=\lim_{n\to\infty}r_{\sharp}'(t)(\phi_{n, \infty}(L_n'(p)))
=\lim_{n\to\infty}(\lim_{k\to\infty}\tau_t(\phi_{n,k}(L_n'(p))))=\tau_t(p)=t([p]).
$$
It follows that $r_{\sharp}\circ r_{\sharp}'={\rm
id}_{S_u(K_0(A))}.$ On the other hand, let $e\in M_K(C_n)$ be a
projection and $s\in S_u(K_0(C)).$
Let $\sigma_s\in T(C)$ which induces $s.$ Then, by (6),
$$
r_{\sharp}'(r_{\sharp}(s))([\phi_{n,\infty}(e)])=\lim_{k\to\infty}
r_{\sharp}(s)([\phi_{n,k}(e)])
=\lim_{k\to\infty}(\lim_{m\to\infty}
\sigma_s(\phi_{m,\infty}(L_m'(\phi_{n,k}(e))))=\sigma_s(\phi_{n,\infty}(e))=
s([\phi_{n,\infty}(e)]).
$$
Thus $r_{\sharp}'\circ r_{\sharp}={\rm id}_{S_u(K_0(C))}.$ \QED

\vspace{0.1in}

{\bf 10.9 Lemma} {\it Let $A$ be a unital separable simple nuclear
\CA\, with $TR(A)\le 1.$  Then the map $\lambda: T(A)\to
S_u(K_0(A))$ maps $\partial_e(T(A))$ onto
$\partial_e(S_u(K_0(A))).$ Moreover, if $A$ is infinite
dimensional,\\ $K_0(A)/{\rm tor}(K_0(A))\not \cong {\mathbb Z}.$
In particular, there is a unital simple AH-algebra $B$ with no
dimension growth described in 10.1 such that
$$
 (K_0(A), K_0(A)_+, [1_A], K_1(A),
T(A)))=(K_0(B), K_0(B)_+, [1_B], K_1(B), T(B)).
$$
 }

{\it Proof:} We will apply Lemma 10.8. Let $C$ be the inductive
limit of \CA s in ${\cal I}$ as described in 10.8. By 1.11 in
\cite{V1}, the map from $T({\tilde C})$ to $S_u(K_0({\tilde
C}))$ maps extremal points onto extremal points. Let $t_0\in
T({\tilde C})$ be the trace such that $t_0(c)=0$ for all $c\in C$
and let $s_0\in S_u(K_0({\tilde C}))$ such that $s_0(x)=0$ for all
$x\in j_*(K_0(C)),$ where $j: C\to {\tilde C}$ is the embedding.
Note that $T({\tilde C})$ is the set of convex
combinations of $\tau\in T(C)$ and $t_0$ and $S_u(K_0({\tilde C}))$
is the set of convex combinations of $s\in S_{u'}(K_0(C))$ and
$s_0.$ Suppose that $\tau\in \partial_e(T(A)).$ Then, by 10.8, $
r^{\sharp}(\tau)\in
\partial_e(T(C))\subset \partial_e(T({\tilde C}).$
It follows that $r^{\sharp}(\tau)$ gives an extremal
state $s_{\tau}$ in $S_{u}(K_0({\tilde C}).$ It follows that
$s_{\tau}\in \partial_e(S_{u'}(K_0(C))).$
Note that $\lambda(\tau)=r_{\sharp}'(s_{\tau}).$ By 10.8,
this shows that $\lambda(\partial_e(T(A)))\subset \partial_e(S_u(K_0(A))).$
To see that $\lambda(\partial_e(T(A)))=\partial_e(S_u(K_0(A))),$
let $s\in \partial_e(S_u(K_0(A))).$
Set
$$
F=\{\tau\in T(A): \lambda(\tau)=s\}.
$$
It is clearly that $F$ is a closed and convex subset of $T(A).$ Furthermore it
is a face. By the Krein-Milman Theorem, it contains an extremal
point $t.$ Since $F$ is a face, $t\in \partial_e(T(A)).$ Thus
$\lambda(\partial_e(T(A)))=\partial_e(S_u(K_0(A))).$

To see $K_0(A)/{\rm tor}(K_0(A))\not\cong {\mathbb Z}$ when $A$ is
infinite dimensional, we note that $A$ has (SP) by 3.2. Since $A$
is simple, we obtain, for any integer $n>0,$ $n+1$ mutually
orthogonal nonzero projections (see for example 5.5)
$p_1,p_2,...,p_n$ and $q$ in $A$ for which $1=q+\sum_{i=1}^np_i,$
$[p_1]=[p_i]$ ($i=1,2,...,n$) and $[q]\le [p_1].$ This implies
that $K_0(A)/{\rm tor}(K_0(A))\not\cong {\mathbb Z}.$

The last statement follows from 10.1 and the above. \QED

Now by 10.4 and 10.9 we have the following:

\vspace{0.1in}

{\bf 10.10. Theorem}
{\it Let $A$ and $B$ be two unital separable nuclear simple \CA s
with $TR(A)\le 1$ and $TR(B)\le 1$ which satisfy the AUCT.
Then $A\cong B$ if and only if
$$
(K_0(A), K_0(A)_+, [1_A], K_1(A), T(A))
\cong (K_0(B), K_0(B)_+, [1_B], K_1(B), T(B))
$$
in the sense of 10.2.
}

\vspace{0.1in}

{\bf 10.11. Remark } Until very recently the only known
interesting examples of simple \CA s with $TR(A)=1$ are
constructed from inductive limits. In a recent paper (\cite{LO}),
it is shown  that separable simple \CA s with $TR(A)=1$ can be
constructed from crossed products. It is shown that certain
crossed products which satisfy some version of Rokhlin property
have tracial topological rank one.

\end{document}